\documentclass{amsproc}

\usepackage{amssymb}

\usepackage[cmtip,all]{xy}

\usepackage{amscd}
\usepackage{mathrsfs}
\usepackage{amsthm}
\usepackage{amsmath}
\usepackage{enumerate}
\usepackage{indentfirst}
\usepackage{xcolor}
\usepackage{comment}
\usepackage[obeyspaces,hyphens,spaces]{url}
\usepackage[hyperfootnotes=true]{hyperref}

\usepackage{etoolbox}

\newtheorem{theorem}{Theorem}[section]
\newtheorem{lemma}[theorem]{Lemma}

\theoremstyle{definition}
\newtheorem{definition}[theorem]{Definition}
\newtheorem{example}[theorem]{Example}

\theoremstyle{remark}
\newtheorem{remark}[theorem]{Remark}

\numberwithin{equation}{section}

\newtheorem{conjecture}[theorem]{Conjecture}

\newcommand{\fkm}{{\mathfrak m}}

\newcommand{\fkp}{{\mathfrak p}}

\newcommand{\fkS}{{\mathfrak S}}

\newcommand{\mS}{{\mathscr S}}
\newcommand{\Z}{\mathbb{Z}}

\newcommand{\Q}{\mathbb{Q}}
\newcommand{\Qlbar}{\overline{\mathbb{Q}}_\ell}
\newcommand{\Qbar}{\overline{\mathbb{Q}}}
\newcommand{\Qpbar}{\overline{\mathbb{Q}}_p}
\newcommand{\Flbar}{\overline{\mathbb{F}}_\ell}
\newcommand{\F}{\mathbb{F}}
\newcommand{\G}{\mathbb{G}}

\newcommand{\A}{\mathbb{A}}
\newcommand{\C}{\mathbb{C}}
\newcommand{\R}{\mathbb{R}}
\newcommand{\T}{\mathbb{T}}
\newcommand{\cA}{{\mathcal{A}}}
\newcommand{\cO}{{\mathcal{O}}}

\newcommand{\cI}{{\mathcal{I}}}
\newcommand{\cM}{{\mathcal{M}}}

\newcommand{\cS}{{\mathcal{S}}}

\newcommand{\cL}{{\mathcal{L}}}
\newcommand{\cU}{\mathcal{U}}
\newcommand{\cV}{\mathcal{V}}

\newcommand{\Fl}{\mathscr{F}\ell}
\newcommand{\HE}{\mathcal{HE}}

\newcommand{\Gal}{\operatorname{Gal}}

\newcommand{\Sh}{\operatorname{Sh}}

\newcommand{\Spec}{\mathrm{Spec}\,}
\newcommand{\Spa}{\mathrm{Spa}\,}

\newcommand{\Frob}{\tu{Frob}}
\newcommand{\Std}{\tu{Std}}

\newcommand{\tr}{\tu{tr}\,}
\newcommand{\Hom}{\tu{Hom}}

\newcommand{\End}{\tu{End}}

\newcommand{\Sp}{\tu{Sp}}
\newcommand{\SO}{\tu{SO}}

\newcommand{\GL}{\tu{GL}}
\newcommand{\GSO}{\tu{GSO}}
\newcommand{\GSp}{\tu{GSp}}

\newcommand{\Lie}{\tu{Lie}\,}
\newcommand{\diag}{\tu{diag}}

\newcommand{\ad}{\tu{ad}}

\newcommand{\Irr}{\tu{Irr}}
\newcommand{\Fpbar}{\overline{\F}_p}

\newcommand{\Res}{\mathrm{Res}}

\newcommand{\tu}{\textup}
\newcommand{\ol}{\overline}
\newcommand{\ra}{\rightarrow}
\newcommand{\isom}{\stackrel{\sim}{\ra}}
\newcommand{\m}{\mathfrak{m}}
\newcommand{\p}{\mathfrak{p}}
\newcommand{\HT}{\mathrm{HT}}

\numberwithin{equation}{section}

\begin{document}


\title[Langlands reciprocity for $\mathrm{GL}_n$: Shimura varieties and beyond]{Recent progress on Langlands reciprocity for $\mathrm{GL}_n$:\\ Shimura varieties and beyond}


\author{Ana Caraiani} 
\address{Department of
  Mathematics, Imperial College London,
  London SW7 2AZ, UK}
\email{caraiani.ana@gmail.com}
\thanks{A.C.~was supported in part by ERC Starting Grant 804176, by a Royal Society University Research Fellowship and by a Leverhulme Prize.}


\author{Sug Woo Shin}
\address{Department of Mathematics, UC Berkeley, Berkeley, CA 94720, USA / Korea Institute for Advanced Study, Seoul 02455, Republic of Korea}
\email{sug.woo.shin@berkeley.edu}
\thanks{S.W.S.~was partially supported by NSF grant DMS-2101688, NSF RTG grant DMS-1646385, and a Simons Fellowship no.1013886.}


\date{}

\begin{abstract}
    The goal of these lecture notes is to survey progress on the global Langlands reciprocity conjecture for $\GL_n$ over number fields from the last decade and a half. We highlight results and conjectures on Shimura varieties and more general locally symmetric spaces, with a view towards the Calegari--Geraghty method to prove modularity lifting theorems beyond the classical setting of Taylor--Wiles. 
    
\end{abstract}

\maketitle

\makeatletter
\patchcmd{\@dottedtocline}
  {\rightskip\@tocrmarg}
  {\rightskip\@tocrmarg plus 4em \hyphenpenalty\@M}
  {}{}
\makeatother

\setcounter{tocdepth}{2}
\tableofcontents

\section{Introduction}

The Langlands reciprocity conjecture for $\mathrm{GL}_n$ predicts a precise correspondence between automorphic representations of $\mathrm{GL}_n$ and $n$-dimensional $\ell$-adic Galois representations over each global field. 
Over global function fields, this is a celebrated theorem by L.~Lafforgue. While the conjecture is still open over number fields, Shimura varieties have been at the center of our approach to the problem ever since Eichler, Shimura, Deligne, Ihara, and Langlands analyzed the special case of modular curves. This led to the birth of the Langlands--Kottwitz--Rapoport (LKR) method, consisting of a web of conjectures concerning the geometry of Shimura varieties and the endoscopic classification of automorphic representations, to study the $\ell$-adic cohomology of Shimura varieties. The LKR method almost reached the current form at the conjectural level in the 1990s, though there still remains much to be investigated.

In the early 2010s, new perspectives on Shimura varieties and their cohomology emerged in connection with Langlands reciprocity conjecture. Calegari--Geraghty proposed that many new cases of modularity would be accessible if we better understood the torsion cohomology of locally symmetric spaces, of which (the underlying manifolds of) Shimura varieties are examples; in particular if we had ``torsion'' Galois representations attached to the torsion cohomology classes. The existence of such Galois representations appeared to be a mere dream apart from some computational evidence in low rank, but Scholze's revolutionary use of perfectoid Shimura varieties realized the torsion Galois representations. These developments stimulated further progress on the cohomology of Shimura varieties and locally symmetric spaces with torsion coefficients, as well as the modularity lifting theorems themselves, with a view towards implementing the Calegari--Geraghty method in some key cases, for example for elliptic curves over CM fields or for abelian surfaces over totally real fields. 

Our goal in this survey is to give a gentle introduction to the aforementioned conjectures and results with an emphasis on the role of Shimura varieties. In \S~\ref{ss:conj-and-main-thms} we review the Langlands reciprocity conjecture for $\mathrm{GL}_n$ with both characteristic 0 and torsion coefficients. The rest of \S~\ref{s:construction} is devoted to one direction of the conjecture, from the automorphic side to the Galois side. After a brief review of the LKR method and its consequence on the so-called conjugate self-dual case (\S~\ref{ss:perfectoid-Shimura}), we proceed to the perfectoid geometry of Shimura varieties and its application to the construction of torsion Galois representations (\S~\ref{ss:construction-torsion-Galois}).

In \S~\ref{s:calegari-geraghty} we give a brief introduction to the Calegari--Geraghty method for proving modularity lifting theorems in the setting of the Betti cohomology of locally symmetric spaces for $\GL_n$ over a general number field $F$. This is a generalization of the famous Taylor--Wiles method. Firstly, in \S~\ref{ss:prerequisites to CG}, we discuss the prerequisites to the Calegari--Geraghty method - what we need to understand about torsion in the cohomology of these locally symmetric spaces in order to implement the method. We state several conjectures and discuss the state-of-the-art results towards them. Secondly, in \S~\ref{ss:CG sketch}, we aim to give the reader a flavour of the method by contrasting it with the original Taylor--Wiles method. We sketch both methods in special cases in order to highlight the role of the various conjectures on locally symmetric spaces.   

In \S~\ref{s:vanishing}, we revisit Shimura varieties and the emerging web of conjectures concerning their cohomology with torsion coefficients. In \S~\ref{ss:geometry of pi_HT}, we discuss the geometry of the Hodge--Tate period morphism, which has been instrumental in computing a suitably generic part of the cohomology of these Shimura varieties and which is likely to play a fundamental role in obtaining a full understanding. In \S~\ref{ss:vanishing}, we discuss several different approaches towards computing the generic part of the cohomology of Shimura varieties PEL type A and point towards generalizations. In \S~\ref{ss:boundary}, we come back full circle and discuss applications of these results to modularity via the boundary of the Borel--Serre compactification of these Shimura varieties.  

Finally, we confess that we have barely touched on exciting new directions in the Langlands program such as the derived and categorical aspects as well as further new geometry over local nonarchimedean fields. The reader is urged to refer to \cite{FengHarrisIHES}, \cite{EmertonGeeHellmannIHES}, and \cite{FarguesScholzeIHES} in these proceedings to learn about them.

\subsection{Notation and conventions}\label{ss:notation}

Given a mathematical object $V$ over a ring $R$ (or $\Spec R$) such as a module or a scheme, let $V_{R'}$ denote the base change of $V$ along a ring homomorphism $R\ra R'$.

When $k$ is a field, let $\ol k$ denote its separable closure. For a finite extension $k'/k$ write $\Res_{k'/k}$ for the functor carrying algebraic groups over $k'$ to those over $k$. When $k'/k$ is a Galois extension, write $\Gamma_{k'/k}$ for the Galois group, and set $\Gamma_k:=\Gamma_{\ol k/k}$. If $k$ is a global field, $\A_k$ will mean the ring of ad\`eles. Given $S$ a finite set of places of $k$, denote by $\A_k^S$ the ring of ad\`eles away from $S$ (with the $S$-components removed from $\A_k$). When $S$ is the set of infinite places of $k$, we also write $\A_k^\infty$. If $k$ is a global field and $v$ is a place of $k$, write $W_{k_v}$ for the local Weil group over $k_v$. Upon choosing a $k$-algebra embedding $i_v:\ol k\hookrightarrow \ol k_v$, we have induced embeddings of $W_{k_v}$ and $\Gamma_{k_v}$ into $\Gamma_k$. The conjugacy class of each induced embedding is independent of the choice of $i_v$. As this is harmless for representation-theoretic considerations, we will usually keep the choice of $i_v$ implicit. If $v$ a finite place, let $\Frob_v$ denote a geometric Frobenius element at $v$ in $W_{k_v}$ or $\Gamma_{k_v}$ (or in $\Gamma_k$ via $i_v$).

Let $G$ be a connected reductive group over a global field $F$.
We fix a quasi-split inner form $G^*$ of $G$ over $F$ together with an inner twist $G_{\ol F}\simeq G^*_{\ol F}$ and a pinning for $G^*$ defined over $F$. Accordingly, as in \cite[\S1.2]{KottwitzShelstad}, we fix the dual group $\hat G=\hat G^*$ with an $L$-action as well as a $\Gamma_F$-equivariant bijection between the dual root datum of $G$ and the root datum of $\hat G$. This determines the $L$-group $\hat G\rtimes W_F$ as in \emph{loc.~cit.}, which is useful in the study of endoscopy, but when considering $\ell$-adic Galois representations valued in the $L$-group, we take the more convenient definition $^L G:=\hat G\rtimes \Gal(F'/F)$ for a finite Galois extension $F'$ of $F$ over which $G$ becomes split; this makes sense as the $L$-action factors through $\Gal(F'/F)$. We will be only brief about $C$-groups, but the $C$-group of $G$ can be formed using the Galois group for the same $F'$ over $F$. (It turns out that the choice of $F'$ does not matter for our purpose.)
The pinning for $G^*$ determines a hyperspecial subgroup $K^*_v\subset G^*(F_v)$ at every finite place $v$ where $G^*$ is unramified (e.g., see \cite[\S4.1]{WaldspurgerLendoscopietordue}). Thereby, via the inner twist, a hyperspecial subgroup $K_v\subset G(F_v)$ is determined at every $v$ where $G$ is unramified.

When $\rho$ is a finite dimensional representation of a group, write $\rho^{\tu{ss}}$ for its semisimplification. If $\rho$ is a continuous representation, so is $\rho^{\tu{ss}}$.

By convention, the algebraic groups $\Sp_{2g}$ and $\GSp_{2g}$ are defined over $\Z$ by the symplectic form on $\Z^g\times \Z^g$ given by $((a_i),(b_j))\mapsto \sum_{i=1}^g (a_i b_{g+i}-a_{g+i}b_i)$. 

Write $\C_p$ for the completion of $\Qpbar$, and $\cO_{\C_p}$ for the ring of integers in $\C_p$.

\subsection*{Acknowledgments}
We are grateful to the organizers Pierre-Henri Chaudouard, Wee Teck Gan, Tasho Kaletha, and Yiannis Sakellaridis of the IHES summer school, July 11--29, 2022, for inviting us to participate and contribute this article, which is an outgrowth of the lecture notes by the authors. Section~\ref{s:vanishing} of this article expands the discussion on the same topic in AC's contribution to the Proceedings of the ICM 2022, therefore there is some overlap with~\cite{CaraianiICM}. 

We are grateful to the Hausdorff Institute of Mathematics for hosting the trimester program ``The Arithmetic of the Langlands program'', during which part of this article was written. We thank Linus Hamann, Teruhisa Koshikawa, Arno Kret, James Newton, Peter Scholze, Matteo Tamiozzo and Mingjia Zhang for valuable discussions that have informed our understanding of the material in this article. We also thank Frank Calegari, Tasho Kaletha, Alberto Minguez, Matteo Tamiozzo, and Yihang Zhu for their comments on an earlier draft. Finally we appreciate many helpful comments and suggestions from the anonymous referees.

\section{Construction of Galois representations with torsion coefficients}\label{s:construction}

\subsection{Conjecture and Main Theorems}\label{ss:conj-and-main-thms}

We recall Buzzard--Gee's formulation of the global Langlands correspondence for reductive groups in the automorphic-to-Galois direction. Then we specialize to the case of general linear groups to state the main theorem of this section, with both characteristic zero and torsion coefficients.
Throughout this section, let $G$ be a connected reductive group over a number field $F$. We adopt the convention of \S~\ref{ss:notation} for the pinning and $L$-group $^L G$ of $G$. In particular, $^L G$ is the semi-direct product of the dual group $\hat G$ and a finite Galois group over $F$.

\subsubsection{On $L$-normalization vs $C$-normalization}\label{sss:L-vs-C}

Buzzard--Gee \cite{BuzzardGee} introduce the notion of \emph{$L$-algebraicity}, resp.~\emph{$C$-algebraicity}, for automorphic representations $\pi$ of $G(\A_F)$, according as the infinitesimal characters of the archimedean components of $\pi$ are integral, resp.~integral after shifting by the half sum of positive roots of $G$. (These conditions will be made explicit for $G=\GL_n$ below.)

They go on to formulate precise conjectures on the corresponding Galois representations, building upon Clozel's work \cite{ClozelAnnArbor}. In the $C$-algebraic case, the Galois representations are valued in the so-called $C$-group $^C G$ of $G$, which can be thought of either as the $L$-group of a natural $\G_m$-extension of $G$ \cite[\S5.3]{BuzzardGee}, or as a certain semi-direct product of $^L G$ with $\G_m$, cf.~\cite[\S1.1]{ZhuSatake}.

In this article, we are mainly interested in automorphic representations contributing to the cohomology of locally symmetric spaces, which are always $C$-algebraic \cite[\S7.2]{BuzzardGee} (but not necessarily $L$-algebraic). Hence the associated Galois representations are to be valued in $^C G(\Qlbar)$; compare the conjectures in \S~3.2 and \S~5.3 of \cite{BuzzardGee}. As such, Conjecture \ref{conj:GLC-ell-adic-G} below should be stated with $^C G$ in place of $^L G$ in general. However we made the deliberate choice to use $^L G$ everywhere in favor of simplicity, as this is harmless for $G=\GL_n$, our main case of interest. In fact, $^C G$ is isomorphic to the \emph{direct} product of $^L G$ and $\G_m$ for $\GL_n$ and many other reductive groups (for instance if $\rho_{\tu{ad}}$ lifts to a cocharacter of $\hat G$; see \cite[\S1.1, Example 2 (3)]{ZhuSatake}). For such groups, we can project from $^C G$ to $^L G$ and forget about the $\G_m$-factor in $^C G$ (since the Galois representation into $\G_m(\Qlbar)$ should be given by the $\ell$-adic cyclotomic character, no information is lost). Thereby $^L G$ can be used as the target group of Galois representations.
ß
Here is an alternative perspective. Whenever $G$ admits a twisting element in the sense of \cite[\S5.2]{BuzzardGee}, which is the case for $G=\GL_n$, one can always twist a $C$-algebraic representation~$\pi$ by a character to make it $L$-algebraic. Since $L$-algebraic automorphic representations are expected to correspond to $^L G$-valued Galois representations,  Conjecture \ref{conj:GLC-ell-adic-G} should be fine as written as long as the twisting element is incorporated into the statement.

\subsubsection{The Buzzard--Gee conjecture}\label{sss:BG-conj}

Let  $\ell$ be a prime. Fix a field isomorphism $\iota:\C\isom \Qlbar$ (disregarding topology). 
 Let $S$ be a finite set of places of $F$ including all places above $\ell$ or $\infty$ as well as the places where $G$ is ramified. Write $S^\infty\subset S$ for the subset of finite places. 
 Let $K_v\subset G(F_v)$ be a family of hyperspecial subgroups at finite places $v\notin S$, chosen at almost all places following \S~\ref{ss:notation}. Set $F_\infty:=F\otimes_{\Q} \R$. Let $K_\infty\subset G(F_\infty)$ be a maximal compact subgroup, and $K_\infty^0$ its neutral component in the archimedean topology. Set $K'_\infty:=K^0_\infty A_G(\R)^0$, where $A_G$ is the maximal $\Q$-split subtorus in $\Res_{F/\Q}Z_G$.
 
Define $\mathcal A_{\tu{ac}}(G,F)$ to be the set of $C$-algebraic
cuspidal automorphic representation of $G(\A_F)$. Write $\mathcal A_{\tu{rac}}(G,F)$ for the subset consisting of $\pi\in\mathcal A_{\tu{ac}}(G,F)$ that is regular, i.e., $\pi_v$ has regular infinitesimal character at all infinite places $v$ of~$F$. Define $\mathcal A_{\tu{coh}}(G,F)$ to be the set of $L^2$-discrete automorphic representations $\pi$ of $G(\A_F)$ which are cohomological, i.e., there exists an irreducible algebraic representation $\xi$ of $(\Res_{F/\Q}G)\times_{\Q}\C$ such that the relative Lie algebra cohomology $H^*(\Lie G(F_\infty),K'_\infty,\pi_\infty\otimes \xi)$ is non-vanishing in some degree. Write  $\mathcal A_{\tu{coh,c}}(G,F)$ for the subset of cuspidal members in $\mathcal A_{\tu{coh}}(G,F)$. We have the inclusions
\begin{equation}\label{eq:three-As}
    \mathcal A_{\tu{coh,c}}(G,F) \subset \mathcal A_{\tu{rac}}(G,F) \subset \mathcal A_{\tu{ac}}(G,F).
\end{equation}
The second inclusion is strict unless $G$ is a torus. The first inclusion turns out to be an equality for $G=\GL_n$ by \cite[\S3]{ClozelAnnArbor}; for general $G$, the equality should follow from \cite[Thm.~1.8]{Salamanca-Riba}. A superscript $S$ over each set will mean the subset consisting of $\pi$ such that $\pi_v$ is unramified for every $v\notin S$.

Let $k\in \{\Flbar,\Qlbar,\C\}$ with the discrete topology on $\Flbar$ and the usual topology on $\Qlbar,\C$.
Here is a variant of $\mathcal A_{\tu{coh}}(G,F)$ with $k$-coefficients.  Let $K_{S^\infty}\subset \prod_{v\in S^\infty} G(F_v)$ be a sufficiently small open compact subgroup, which is allowed to vary. Take $K=K_{S^\infty}\times \prod_{v\notin S} K_v\subset G(\A_F^\infty)$.
Consider the locally symmetric space
\begin{equation}\label{eq:loc-sym-space}
    Y_K:= G(F)\backslash G(\A_F)/K K'_\infty.
\end{equation}
The double coset algebra $\T^S:=\Z[K^S\backslash G(\A_F^S)/K^S]$, referred to as the (abstract) Hecke algebra (away from $S$), acts on the Betti cohomology: 
\begin{equation}\label{eq:TS-on-Betti}
    \T^S~ \circlearrowleft ~\bigoplus_{i\ge 0} 
    \Big( \varinjlim_{K_{S^\infty}} H^i(Y_K,k) \Big).
\end{equation}
Let $\mathcal{HE}^S_{\tu{coh}}(G,F)_{k}$ denote the set of eigencharacters $\T^S\ra k$ appearing in \eqref{eq:TS-on-Betti}, which are also known as \emph{Hecke eigensystems} or \emph{Hecke eigencharacters}. We often denote an eigencharacter by its kernel $\fkm\subset \T^S$, which is a maximal ideal, understanding that the eigencharacter also comes with an embedding $\T^S/\fkm\hookrightarrow k$.
For $S\subset S'$ there is a natural inclusion $\mathcal {HE}^S_{\tu{coh}}(G,F)_{k}\subset \mathcal {HE}^{S'}_{\tu{coh}}(G,F)_{k}$ compatibly with the obvious map $\T^{S'}\ra \T^{S}$. 

We can define a more general set $\mathcal {HE}^S_{\tu{coh}}(G,F)'_{k}$ by replacing the coefficient field $k$ in \eqref{eq:TS-on-Betti} with local systems arising from irreducible algebraic representations of $\Res_{F/\Q}G$ (or a suitable integral model thereof, if $\tu{char}(k)>0$) on $k$-vector spaces. 
When $k=\C$, we have a natural map
\begin{equation}\label{eq:A-to-HE}
    \mathcal{A}^S_{\tu{coh}}(G,\C)\ra \mathcal {HE}^S_{\tu{coh}}(G,F)'_{\C}. 
\end{equation}
The map assigns to each $\pi$ the Hecke eigensystem arising from $\pi^S$. This map is an injection for $G=\GL_n$ by the strong multiplicity one theorem but not in general. When $k=\Flbar$, it turns out that 
$\mathcal {HE}^S_{\tu{coh}}(G,F)'_{k}=\mathcal {HE}^S_{\tu{coh}}(G,F)_{k}$, namely no new Hecke eigensystems are obtained from non-constant coefficients. This ultimately comes from the fact that a continuous representation of an $\ell$-adic reductive group on a finite-dimensional $\Flbar$-space is trivial on an open subgroup. (To apply this fact, note that the places of $F$ over $\ell$ are in $S$ and that the limit is taken over open compact subgroups at such places in \eqref{eq:TS-on-Betti}.)

Define $\mathcal{G}({}^L G,F)_{k}$ to be the set of weak equivalence classes of continuous semisimple representations $\rho:\Gamma_F\ra {}^L G(k)$ that are unramified at almost all places of $F$, where $\rho_1$ and $\rho_2$ are considered weakly equivalent if the semisimple parts $\rho_1(\Frob_v)_{\tu{ss}}$ and
$\rho_2(\Frob_v)_{\tu{ss}}$ are $\hat G(k)$-conjugate at every finite place $v\nmid \ell$ where both $\rho_1$ and $\rho_2$ are unramified. When $k=\Qlbar$, we introduce the subset $\mathcal{G}_{\tu{dR}}({}^L G,F)_{\Qlbar}$ consisting of $\rho\in\mathcal{G}({}^L G,F)_{\Qlbar}$ 
such that $\rho|_{\Gamma_{F_v}}$ is a de Rham representation at every $v|\ell$. (An $L$-group valued representation is said to be de Rham if its composition with some, thus every, faithful representation of $^L G$ is de Rham in the usual sense.) As in the automorphic case, if we write $\mathcal{G}^S$ in place of $\mathcal{G}$, it denotes the subset of (weak equivalence classes of) $\rho$ which are unramified outside $S$.

We are ready to state a coarse form of the Buzzard--Gee conjecture (cf.~the conjectures in \cite[\S3.2, \S5.3]{BuzzardGee}, and their refinements \cite[\S4]{JohanssonThorne} in the cohomological case); this is statement (1) below. An analogous conjecture with $\Flbar$-coefficients was suggested by Ash \cite{Ash1,Ash2}, at least for $G=\GL_n$; see statement (2) below.
(Recall from \S~\ref{sss:L-vs-C} that $C$-groups should be used to state the conjecture in general, but that we are glossing over this point.)

\begin{conjecture}\label{conj:GLC-ell-adic-G}
The following are true.
\begin{enumerate}
 \item There exists a map (which depends on the choice of $\iota$)
$$\tu{GLC}_{G,\Qlbar}:\mathcal A^S_{\tu{ac}}(G,F) \ra \mathcal G^S_{\tu{dR}}({}^L G,F)_{\Qlbar},\qquad \pi=\otimes'_v \pi_v \mapsto \rho_{\pi,\iota},$$
such that $\iota\pi_v$ corresponds to $(\rho_{\pi,\iota}|_{\Gamma_{F_v}})^{\tu{ss}}$ via the unramified local Langlands correspondence at all $v\notin S$.
\item There exists a map
$$\tu{GLC}_{G,\Flbar}:\mathcal{HE}^S_{\tu{coh}}(G,F)_{\Flbar} \ra \mathcal G^S({}^L G,F)_{\Flbar},\qquad \fkm \mapsto \ol\rho_{\fkm},$$
such that $\ol\rho_{\fkm}$ is unramified at $v\notin S$ and such that $\tu{char}(\ol\rho_{\fkm}(\Frob_v))$ is explicitly described in terms of the Hecke eigenvalues at $v$ encoded in $\fkm$ for every $v\notin S$.
\end{enumerate}
\end{conjecture}
If true, part (1) of the conjecture assigns to each $\pi$ a (weakly) compatible system of Galois representations as the prime $\ell$ and the isomorphism $\iota$ vary in the following sense: for $\iota:\C\isom \Qlbar$ and $\iota':\C\isom \ol{\Q}_{\ell'}$, the conjugacy classes $(\rho_{\pi,\iota}|_{\Gamma_{F_v}})^{\tu{ss}}$ and $(\rho_{\pi,\iota'}|_{\Gamma_{F_v}})^{\tu{ss}}$ are conjugate via $\iota'\iota^{-1}$ and moreover defined over $\Qbar$.

In fact the conjecture by Buzzard--Gee is more precise than part (1) above in that the Hodge--Tate cocharacter and the image of complex conjugations are predicted for $\rho_{\pi,\iota}$. 
On the other hand, it is nontrivial to formulate the correct local-global compatibility at ramified places in both (1) and (2). Already in (1) (the problem is worse in (2)), one reason is that the local $L$-packets and $A$-packets are only conjectural in general; another problem is the so-called CAP (cuspidal but associated with parabolic subgroups) automorphic representations. However the situation is somewhat favorable for $G=\GL_n$ as we will see below.

The maps $\tu{GLC}_{G,\Qlbar}$ and $\tu{GLC}_{G,\Flbar}$ do not uniquely characterize the isomorphism classes of $\rho_{\pi,\iota}$ and $\ol\rho_{\fkm}$ even if the isomorphism classes of $\rho_{\pi,\iota}|_{\Gamma_{F_v}}$ are prescribed at all places $v$ of $F$. (See \cite{LarsenConjugacy,LarsenConjugacy2} for related discussions.) This is why we phrase the conjecture in terms of weak equivalence classes of Galois representations.
Moreover, the maps $\tu{GLC}_{G,\Qlbar}$ and $\tu{GLC}_{G,\Flbar}$ are not expected to be either injective or surjective. When the coefficient is $\Qlbar$,
the failure of injectivity is related to the presence of local packets;
the problem with surjectivity amounts to describing the image of $\tu{GLC}_{G,\Qlbar}$, which should incorporate the phenomenon of global $A$-packets among other things.

\begin{remark}
The preceding discussion suggests that the characteristic 0 analogue of $\mathcal{HE}^S_{\tu{coh}}(G,F)_{\Flbar}$ should be $\cA^S_{\tu{coh,c}}(G,F)$ or $\mathcal{HE}^S_{\tu{coh}}(G,F)'_{k}$ with $k\in \{\C,\Qlbar\}$, which is less general than $\mathcal A^S_{\tu{ac}}(G,F)$. 
Then it is natural to ask whether $\mathcal A^S_{\tu{ac}}(G,F)$ admits mod $\ell$ automorphic analogues, such that the latter accommodates ``mod $\ell$'' of non-regular members of $\mathcal A^S_{\tu{ac}}(G,F)$ in particular. In the special case when $G$ is part of a Shimura datum, a partial answer may be given by the coherent cohomology (see the introductory paragraphs in \S \ref{s:calegari-geraghty}) with mod $\ell$ coefficients. Other than that, the authors have no idea.
\end{remark}

\subsubsection{The $\GL_n$-case}

Now we specialize to the case of $G=\GL_n$, and give more precise statements and results. Let $k\in \{\Flbar,\Qlbar,\C\}$ as before.
In this case, we write $\mathcal{A}(n,F)$, $\mathcal{HE}^S_{\tu{coh}}(n,F)_k$,  $\mathcal{G}(n,F)_k$ for $\mathcal{A}(\GL_n,F)$, $\mathcal{HE}^S_{\tu{coh}}(\GL_n,F)_k$,  $\mathcal{G}(\GL_n,F)_k$, and likewise with more superscripts and subscripts.  Each member of $\mathcal{G}(n,F)_k$ is represented by an $n$-dimensional representation $\Gamma_F\ra \GL_n(k)$, and each weak equivalence is the same as an isomorphism class in the usual sense.

For each $\pi\in \mathcal{A}(n,F)$ and each embedding $\tau:F\hookrightarrow \C$, which determines an archimedean place $v$ of $F$ together with $F_v\hookrightarrow \C$, the infinitesimal character of $\pi_v$ may be viewed as a multi-set of $n$ complex numbers by the Harish-Chandra isomorphism, to be denoted by $\tu{HC}_\tau(\pi)$. Then $\pi$ is $C$-algebraic if and only if $\tu{HC}_{\tau}(\pi)-\frac{n-1}{2}$ is integral (i.e., $a\in\Z+\frac{n-1}{2}$ for all $a\in \tu{HC}_{\tau}(\pi)$) for every $\tau$.

We start by specifying the normalization of local Langlands in part (1) of Conjecture \ref{conj:GLC-ell-adic-G}. Let $\tu{LL}$ denote the unramified local Langlands correspondence for $\GL_n$ over $F_v$ with $\C$-coefficients, with $v$ non-archimedean, carrying Satake parameters to the eigenvalues of geometric Frobenius elements.\footnote{That is, if continuous characters $\chi_{i,v}:F^\times_v\ra \C^\times$ and $\sigma_{i,v}:W_{F_v}\ra\C^\times$ correspond via local class field theory (normalized to match uniformizers of $F_v$ with geometric Frobenius elements of $W_{F_v}$) for $i=1,...,n$, then the unique irreducible unramified subquotient of the normalized parabolic induction of $\otimes_{i=1}^n \chi_{i,v}$ from $\prod_{i=1}^n \GL_1(F_v)$ to $\GL_n(F_v)$ is sent to the $n$-dimensional representation $\oplus_{i=1}^n\sigma_{i,v}$ of $W_{F_v}$.}
Then our requirement on $\pi\mapsto \rho_{\pi,\iota}$ is that the semisimplification of $\iota^{-1}\rho_{\pi,\iota}|_{W_{F_v}}$ is isomorphic to $\tu{LL}(\pi_v\otimes|\det|^{(1-n)/2})$. In this case, the Hodge--Tate weights of $\rho_{\pi,\iota}$ relative to $\tau:F\hookrightarrow \Qlbar$ 
(defined in \cite[\S2.2.2]{PatrikisMAMS}, for example)
should be given by the multi-set $\tu{HC}_{\iota^{-1} \tau}(\pi)+\frac{n-1}{2}$.

Next we make part (2) of Conjecture \ref{conj:GLC-ell-adic-G} concrete.
At each finite place $v$ of $F$ not above $\ell$, write $q_v$ for the residue field cardinality at $v$, and $\varpi_v$ for a uniformizer of $F_v$. The square root $q_v^{1/2}\in \R_{>0}^\times$ is viewed as an element of $\Qlbar^\times$ via $\iota$, or as an element of $\Flbar^\times$ via the reduction map (since $q_v^{1/2}$ is an $\ell$-unit).
When $v\notin S$, define $T_{v,i}$ to be the double coset operator 
$K_v\left(\begin{smallmatrix}
    \varpi_v I_{i} & 0 \\ 0 & I_{n-i}
\end{smallmatrix}\right)K_v$, which is independent of the choice of $\varpi_v$ and viewed as an element of $\T^S$.
There is a ``universal'' Hecke polynomial
$$
H_v(x):=x^n-T_{v,1}x^{n-1} + q_v T_{v,2}x^{n-2}+\cdots + (-1)^n q_v^{n(n-1)/2}T_{v,n}\in \T^S[q_v^{1/2}][x].
$$
The Hecke polynomial for each Hecke eigencharacter $\fkm\in \mathcal {HE}^S_{\tu{coh}}(n,F)_{k}$ at $v$ is defined by applying it to the coefficients of $H_v$:
$$H_v(x)_{/\fkm}\in (\T^S/\fkm)(q_v^{1/2})[x]\subset k(q_v^{1/2})[x].$$
When $k=\Flbar$, the condition at each $v$ in Conjecture \ref{conj:GLC-ell-adic-G} (2) can now be made concrete: it is the equality 
\begin{equation}\label{eq:Frobpoly=Heckepoly}
    \tu{char}(\ol\rho_{\fkm}(\Frob_v))=H_v[x]_{/\fkm}.
\end{equation}
When $k=\Qlbar$, if $\fkm_{\pi}$ denotes the image of $\pi\in \mathcal{A}^S_{\tu{coh}}(n,F)$ under \eqref{eq:A-to-HE} via  $\iota:\C\simeq \Qlbar$, then the equality
\eqref{eq:Frobpoly=Heckepoly} is equivalent to the condition in Conjecture \ref{conj:GLC-ell-adic-G} (1) that $\iota \pi_v$ corresponds to $(\rho_{\pi,\iota}|_{\Gamma_{F_v}})^{\tu{ss}}$.
Thus the two parts of Conjecture \ref{conj:GLC-ell-adic-G} are consistent through the operation of taking mod $\ell$ on the Galois side.

In the conjecture for $G=\GL_n$, 
the representations $\rho_{\pi,\iota}$ and $\ol\rho_{\fkm}$ are uniquely characterized by the Brauer--Nesbitt theorem and the Chebotarev density theorem. The map $\tu{GLC}_{n,\Qlbar}$, if exists, must be injective by the strong multiplicity one theorem for $\GL_n$. The conjectural map $\tu{GLC}_{n,\Qlbar}$ should also be injective in view of \eqref{eq:Frobpoly=Heckepoly}.

For each real place $v$ of $F$, write $c_v\in \Gamma_{F}$ for a complex conjugation at $v$ (well-defined up to conjugacy).
Every Galois representation $\rho$ in the image of $\mathcal{A}^S_{\tu{rac}}(n,F)$ under $\mathrm{GLC}_{n,\Qlbar}$ is believed to be \emph{odd} in the sense that $|\tr \rho(c_v)|\in\{0,\pm1\}$ at every real place $v$ of $F$. (Some non-regular $\pi\in\mathcal{A}^S_{\tu{ac}}(n,F)$ also correspond to odd Galois representations, e.g., if $\pi$ comes from a holomorphic cuspidal eigenform of weight 1.) Similarly every $\ol\rho$ in the image under $\mathrm{GLC}_{n,\Flbar}$ is expected to be odd, i.e., $|\tr \ol\rho(c_v)|\in \{0,\pm1\}$ at every real $v$. See \S \ref{sss:self-dual-ref} below for known results. Our definition of oddness does not coincide with the one often found in the literature (e.g., \cite[Def.~1.2]{FKP}) if $n>2$; see Remark \ref{rem:odd} below for a discussion.

We make further observations on Conjecture \ref{conj:GLC-ell-adic-G}.  In this case, the image of $\mathcal{A}^S_{\tu{ac}}(n,F)$ under $\tu{GLC}_{n,\Qlbar}$ is believed to be exactly the irreducible representations in $\mathcal{G}^S_{\tu{dR}}(n,F)_{\Qlbar}$. (The dictionary in the $\GL_n$-case is that cuspidality corresponds to irreducibility through $\tu{GLC}_{n,\Qlbar}$.)
At least when $F$ is totally real, 
one could guess that an irreducible $\ol\rho\in\mathcal{G}^S(n,F)_{\Flbar}$ belongs to the image of $\tu{GLC}_{n,\Flbar}$ exactly when $\ol\rho$ is odd. This is reasonable for $n=2$ as predicted by the Serre modularity conjecture for Hilbert modular forms. However the situation is unclear for $n>2$ as it does not seem known whether $\ol\rho$ lifts to some $\rho\in \mathcal{G}^S_{\tu{dR}}(n,F)_{\Qlbar}$ even at the expense of ``raising the level''; this case is not covered by lifting results in the literature such as \cite[Thm.~A]{FKP}.

\begin{remark}\label{rem:odd}
One may partly extend from $\GL_n$ to a general $G$ the expectation that the Galois representations arising from $\mathcal{A}^S_{\tu{rac}}(n,F)$ should be odd. 
We propose the following notion of oddness for $\rho\in \mathcal G^S_{\tu{dR}}({}^L G,F)_{\Qlbar}$. For each real place $v$ of $F$, write $M_{F_v}$ for a Levi subgroup of a fundamental parabolic subgroup of $G_{F_v}$. (The latter is a parabolic subgroup minimal for the property of containing a fundamental maximal torus of $G_{F_v}$, cf.~\cite[\S 4.1]{BorelWallach}.) Then 
$M_{F_v}$ is unique up to $G(F_v)$-conjugacy. We want to say that $\rho$ is \emph{odd} if $\rho(c_v)$ at every real place $v$ is $\hat G(\Qlbar)$-conjugate to a conjugacy class in ${}^L M_{F_v}(\Qlbar)$ satisfying Gross's oddness condition for $M_{F_v}$ as stated in  \cite[Def.~1.2]{FKP} (which applies when the coefficient field is either $\Qlbar$ or $\Flbar$). If $G_{F_v}$ contains an elliptic maximal torus then $M_{F_v}=G_{F_v}$, so our oddness coincides with that of \cite[Def.~1.2]{FKP}, but not in general, e.g., for $\GL_n$ with $n>2$. In the latter case, our oddness is equivalent to the condition that $|\tr \rho(c_v)|\in\{0,\pm1\}$ at real places $v$. It seems reasonable to guess that whenever $\pi\in \mathcal{A}_{\tu{coh,c}}(G,F)$ has a local component that is generic (which should be equivalent to having generic components at almost all places by a version of the Shahidi conjecture), the corresponding Galois representation should be odd in our sense.
We leave the reader to formulate the analogues for the case of torsion coefficients.
\end{remark}

To state results towards the above conjecture, we introduce some terminology.
Let $F$ be a totally real field, or a CM field (a totally imaginary quadratic extension of a totally real field). Such an $F$ admits a complex conjugation $c$ that is independent of how $F$ is embedded into $\C$. Write $F^+$ for the fixed subfield of $F$ by $c$, and $N_{F/F^+}$ for the norm map (which is the identity map if $F$ is totally real).
By definition $\pi\in \mathcal A_{\tu{ac}}(n,F)$ is said to be \emph{essentially conjugate self-dual} if $\pi\circ c \simeq \pi^\vee\otimes (\chi\circ N_{F/F^+})$ for a Hecke character $\chi:(F^+)^\times\backslash \A_{F^+}^\times \ra \C^\times$. 
(In fact an additional parity condition used to be imposed on $\chi$ when $F=F^+$, but the condition was shown to be superfluous by \cite{PatrikisSign}.)
Write $\tilde{\mathcal A}_{\star}(n,F)$ for the subset of essentially conjugate self-dual representations in $\mathcal A_{\star}(n,F)$ for $\star\in \{\tu{ac},\tu{rac}\}$. 
Let $\mathcal A^1_{\star}(n,F)$ denote the subset of conjugate self-dual members, namely $\pi$ such that $\pi\circ c\simeq \pi^\vee$.
If $F$ is totally real, the word ``conjugate'' is usually omitted everywhere.

Here is the main theorem of this section. (See the end of \S~\ref{ss:csd}, resp.~\S~\ref{ss:construction-torsion-Galois}, for known results on local-global compatibility at $v\in S$ and other comments.) Part (1) is due to Harris--Lan--Taylor--Thorne \cite{HLTT}. Shortly thereafter, Scholze \cite{Sch15} proved both parts by a different method based on perfectoid Shimura varieties and Hodge--Tate period morphisms. Boxer's thesis \cite{BoxerThesis} suggests an alternative path for both parts.

\begin{theorem}\label{thm:construction-main}
Let $F$ be a totally real or CM field, and $\iota$ as above.
\begin{enumerate}
 \item There exists an injection
$$\tu{GLC}_{n,\Qlbar}:\mathcal A^S_{\tu{rac}}(n,F) \ra \mathcal G^S_{\tu{dR}}(n,F)_{\Qlbar},\qquad \pi \mapsto \rho_{\pi,\iota},$$
such that $\iota\pi_v$ corresponds to $(\rho_{\pi,\iota}|_{\Gamma_{F_v}})^{\tu{ss}}$ via the unramified LLC at all $v\notin S$.
\item There exists a map
$$\tu{GLC}_{n,\Flbar}:\mathcal{HE}^S_{\tu{coh}}(n,F)_{\Flbar} \ra \mathcal G^S(n,F)_{\Flbar},\qquad \fkm \mapsto \ol\rho_{\fkm},$$
such that  $\tu{char}(\ol\rho_{\fkm}(\Frob_v))=H_v(x)_{/\fkm}$ for every $v\notin S$.
\end{enumerate}
\end{theorem}

From now, restrict $F$ to be a totally real or CM field. 
We follow Scholze's proof as it has further applications as presented in \S~\ref{s:vanishing}.
The proof breaks up into two main parts, carried out in the following order:
\begin{enumerate}
   \item[(A)] Prove Theorem \ref{thm:construction-main} (1) for $\pi\in\tilde{\mathcal A}_{\tu{rac}}(n,F)$.
   \item[(B)] Prove Theorem \ref{thm:construction-main} (2) and its analogue for mod $\ell^m$-coefficients, $m\ge 1$.
\end{enumerate}
Part (A) was completed by a combination of works by Clozel, Kottwitz, et al., spanning over nearly two decades, before Theorem \ref{thm:construction-main} was proved in general; see \S~\ref{sss:self-dual-ref} for a guide to references. The argument will be reviewed in \S~\ref{ss:csd}. Explanation of Part (B) will take up \S~\ref{ss:perfectoid-Shimura}--\S~\ref{ss:construction-torsion-Galois} below. 
Once this is done, Theorem \ref{thm:construction-main} (1) in full follows from (2) by realizing $\pi$ in the cohomology and ``gluing'' Galois representations modulo powers of $\ell$. See the proof of \cite[Cor.~5.4.2]{Sch15} for details. 

\begin{remark}
Before we restrict ourselves to the regular or cohomological case from now on, let us mention several results on Conjecture \ref{conj:GLC-ell-adic-G} in the \emph{non-regular} case. The basic principle is that we can often prove the conjecture when $\pi$ contributes to the so-called coherent cohomology of Shimura varieties (even if $\pi$ is not cohomological, i.e., not appearing in the Betti cohomology). When $G=\GL_2$ over $\Q$, such $\pi$ come from weight 1 cuspforms. Deligne--Serre \cite{DeligneSerre} proved the conjecture in this case. For $\GL_2$ over a totally real field, the problem is about Hilbert cuspforms of partial weight 1. Here the result is due to Rogawski--Tunnell, Ohta, and Jarvis. (See \cite{Jarvis} and the references therein.) In semisimple rank greater than 1, there has been much progress since Taylor's work \cite{TaylorGSp4} on $\tu{GSp}_4$ over $\Q$; see the last paragraph of \S~\ref{ss:construction-torsion-Galois} for recent results. All of these results rest on variants and generalizations of Hasse invariants to be discussed in \S~\ref{ss:construction-torsion-Galois}.

However we stress that Conjecture \ref{conj:GLC-ell-adic-G} is still unsettled even when $G=\GL_2$ over $\Q$. The remaining $\pi$ (which do not show up in either coherent or Betti cohomology of modular curves) correspond to Maass cuspforms. Except when $\pi$ arise from Hecke characters over real quadratic fields via automorphic induction (the conjecture is reduced to class field theory for such $\pi$), very little is known.
\end{remark}

\subsection{The conjugate self-dual case}\label{ss:csd}

This subsection sketches the proof of the following special case for essentially conjugate self-dual representations.

\begin{theorem}\label{thm:csd-case}
Theorem \ref{thm:construction-main} (1) is true if $\pi\in\tilde{\mathcal A}_{\tu{rac}}(n,F)$.
\end{theorem}

In fact there is a standard argument, cf.~the proofs of \cite[Thm.~1.1, 1.2]{BLGHT} and \cite[Prop.~4.3.1]{CHT}, to reduce the proof of Theorem \ref{thm:csd-case} to the following case, so we will assume it throughout \S~\ref{ss:csd}:
\begin{itemize}
    \item $F$ is a CM field, and $\pi\in \tilde{\mathcal A}^1_{\tu{rac}}(n,F)$.
\end{itemize}

\subsubsection{Setup}\label{sss:setup}
Before getting to the proof of Theorem \ref{thm:csd-case}, we explain the basic strategy to make progress towards constructing the map $\tu{GLC}_{G,\Qlbar}$ in an idealized setup, based on the Langlands--Kottwitz--Rapoport method for the cohomology of Shimura varieties. 

Let $(G,X)$ be a Shimura datum \cite[\S2.1.1]{DeligneShimuraCorvallis}; see \S~2.3.1 of Caraiani's article in \cite{AWS2017perfectoid} for further explanation. Thus $G$ is a connected reductive group over $\Q$, and $X$ is a Hermitian symmetric domain with a transitive $G(\R)$-action. From $X$ we obtain a conjugacy class $[\mu_X]$ of cocharacters $\G_{m,\C}\ra G_{\C}$. Then $[\mu_X]$ is defined over a number field $E:=E(G,X)\subset \C$ called the reflex field.
We fix a $\Gamma_{\Q}$-pinning $(\hat B, \hat T, \{\hat X_{\alpha}\})$. 
From the inverse of $[\mu_X]$ we obtain a $\hat B$-dominant character $-\mu\in X^*(\hat T)$, which is then fixed by the $\Gamma_E$-action. 
So there exists a representation
$$r_{-\mu}:{}^L G_{E}=\hat G \rtimes \Gamma_E \ra \GL(V_{-\mu}),$$
characterized uniquely up to isomorphism by the two conditions: $r_{-\mu}|_{\hat G}$ is an irreducible representation of highest weight $-\mu$, and $\Gamma_E$ acts trivially on the highest weight space.

We have a projective system of canonical models of Shimura varieties $\Sh=(\Sh_K)$, which are quasi-projective smooth varieties over $E$ and labeled by neat open compact subgroups $K\subset G(\A^\infty)$. The system $\Sh$ is equipped with the ``Hecke'' action defined over $E$ by  $G(\A^\infty)$. Each $\Sh_K(\C)$ is a complex manifold, identified with $Y_K$ in \eqref{eq:loc-sym-space} by construction (except the technical point that $K'_\infty$ may be replaced with a larger subgroup of $G(\R)$ in which $K'_\infty$ has finite index). Each $\xi\in \Irr((\Res_{F/\Q}G)_{\C})$ gives rise via $\iota$ to lisse $\Qlbar$-sheaf $\cL_{\xi,K}$ over $\Sh_K$ compatibly with respect to the transition maps between $(\Sh_K)$. The \'etale cohomology with compact support
\begin{equation}\label{eq:Hixic}
    H^i_{\xi,c}:=\varinjlim_{K} H^i_c(\Sh_{K}\times_E \ol{E},\cL_{\xi,K}),\qquad i\in\Z_{\ge 0},
\end{equation}
is a $G(\A^\infty)\times \Gal(\ol E/E)$-module,  and 
$$[H_{\xi,c}]:=\textstyle\sum_{i\ge 0} (-1)^i [H^i_{\xi,c}]$$
is a virtual $G(\A^\infty)\times \Gal(\ol E/E)$-module (viewed in the Grothendieck group; see \cite[pp.23--25]{HarrisTaylor}). 

\subsubsection{The Langlands--Kottwitz--Rapoport (LKR) method}
\label{sss:LKR-method}
This is a systematic way to compute $[H_{\xi,c}]$ at each rational prime $p\neq \ell$. We restrict ourselves to the case where Shimura varieties have good reduction modulo $p$; see \S\ref{sss:complements-further-refs} below for the case of bad reduction. 

Let $\fkp$ be a prime of $E$ above $p$. Assume that $G_{\Q_p}$ is an unramified group and fix a hyperspecial subgroup $K_p\subset G(\Q_p)$. 
As (sufficiently small) open compact subgroups $K^p\subset G(\A^{\infty,p})$ vary, we have a projective system $\Sh_{K_p}:=\varprojlim_{K^p} \Sh_{K_pK^p}$ over $E$ equipped with a $G(\A^{\infty,p})$-action.
By taking the limit in \eqref{eq:Hixic} over $K^p$, we can define $H^i_{\xi,c,K_p}$, which is the $K_p$-fixed subspace of $H^i_{\xi,c}$, as well as a virtual $G(\A^{\infty,p})\times \Gal(\ol E/E)$-module $[H_{\xi,c,K_p}]$. Let $\fkp$ be a prime of $E$ above $p$.
The LKR method aims to show that the local Galois action at $\fkp$ is unramified and to describe 
$$\mbox{the virtual}~G(\A^{\infty,p})\times \Frob_{\fkp}^{\Z}\mbox{-module}~[H_{\xi,c,K_p}]$$ 
via automorphic representations.
The method consists of the following steps.

\begin{enumerate}
\item Construct an integral model $\mathscr S_{K_p}$ of $\Sh_{K_p}$ over $\cO_{E,\fkp}$, which is canonical in the sense of \cite{KisinModels}.
    \item Establish a trace formula for the action of $G(\A^{\infty,p})\times \Frob_{\fkp}^{\Z}$ on $[H_{\xi,c,K_p}]$ by studying the structure of points on $\mathscr S_{K_p}$ modulo $\fkp$.
    \item Stabilize the trace formula.
 \item Compare the outcome with the stabilization of the Arthur-Selberg trace formula to obtain a representation-theoretic description of $[H_{\xi,c,K_p}]$.
\end{enumerate}

Part (1) was done for PEL-type Shimura varieties by Kottwitz \cite{KottwitzPoints} and for abelian-type Shimura varieties by Kisin \cite{KisinModels} ($p>2$) and Kim--Madapusi Pera \cite{KimMadapusiPera} $(p=2)$.
Until recently, the best result on (2) and (3) has been due to Kottwitz \cite{KottwitzAnnArbor,KottwitzPoints} in the case of PEL-type Shimura varieties.
This has been used to carry out (4) in some important special cases, which were enough to prove Theorem \ref{thm:csd-case}. 

Now the first two steps are complete more generally for Shimura varieties of abelian type in \cite{KSZ} following \cite{KisinLR}; the key is the proof of a version of the Langlands--Rapoport (LR) conjecture regarding the structure of points in (2). (The work of Kottwitz \cite{KottwitzPoints}, and its recent extension to the case of Hodge type by~Lee \cite{LeeLKmethod} bypasses the LR conjecture. However the LR conjecture seems essential to deal with abelian-type Shimura varieties.) The outcome has been employed to prove new cases of Conjecture \ref{conj:GLC-ell-adic-G} (2) for the groups $\GSp_{2n}$ and certain quasi-split forms of $\GSO_{2n}$ over totally real fields. Further applications are being worked out in \cite{KSZ2}.

\subsubsection{Idealized situation}\label{sss:idealized}
Kottwitz gave a conjectural description of the intersection cohomology version of $[H_{\xi,c}]$ (thus also $[H_{\xi,c,K_p}]$) in \cite[\S10]{KottwitzAnnArbor} but it is quite complicated. To convey the basic principle, we restrict our attention to the $\pi^\infty$-isotypic part of the cohomology, where $\pi\in \mathcal{A}_{\tu{coh}}(G,\Q)$, and make the following simplifying hypotheses:
\begin{itemize}
    \item[(a)] $G^{\ad}$ is anisotropic over $\Q$,
    \item[(b)] $\pi$ is essentially tempered at every place,
    \item[(c)] we pretend that endoscopy disappears (e.g., conjugacy classes are stable conjugacy classes, and $L$-packets are singletons).
\end{itemize}
Hypothesis (c) is satisfied if $G$ is an inner form of (restriction of scalars of) $\GL_2$ but typically false unless we restrict attention to only a certain subset of $\pi$.
Hypothesis (a) implies that Shimura varieties and their canonical integral models are proper over the base, so that there is no need for compactification; thus we freely write $H_\xi$ instead of $H_{\xi,c}$. The same hypothesis also implies that each $H^i_{\xi,c}$ is a semisimple $G(\A^{\infty})$-module via Matsushima's formula, since the $L^2$-automorphic spectrum is semisimple.
Hypothesis (b) is supposed to ensure that $\pi$ contributes to cohomology only in the middle degree, namely $H^{d}_{\xi}$ where $d:=\dim \Sh$.
Consider the assignment
\begin{equation}\label{eq:pi-isotypic}
  \pi\quad\mapsto\quad R_{\pi,\mu,\iota}:=\Hom_{G(\A^{\infty})}(\pi^\infty,H^{d}_{\xi}).  
\end{equation}
If we had $\rho_{\pi,\iota}$ (whose construction is our goal), then our expectation, i.e., the Kottwitz conjecture, is essentially that
\begin{equation}\label{eq:Rpi-expectation}
R_{\pi,\mu,\iota}= r_{-\mu}\circ \rho_{\pi,\iota}
\end{equation}
up to twists having to do with normalizations. This leads to:

\medskip

\paragraph{\textbf{An approach to Conjecture \ref{conj:GLC-ell-adic-G}}}
Suppose $F=\Q$. Given $G,\iota$ and $\pi$, choose (possibly several) Shimura data and study the resulting representations $R_{\pi,\mu,\iota}$. From this, construct $\rho_{\pi,\iota}$ and show the desired properties of $\rho_{\pi,\iota}$.

\begin{remark}
 If $F\neq \Q$ in Conjecture \ref{conj:GLC-ell-adic-G}, then the group $G$ over $F$ in the conjecture will differ from the reductive groups used in Shimura data, because the latter are always over $\Q$. Often the latter is chosen to be a form of $\Res_{F/\Q} G$.
\end{remark}

Let us illustrate the mechanism for unitary groups in relation to Theorem \ref{thm:csd-case}. Let $F$ be a CM field. We fix an algebraic closure $\ol \Q\subset \C$. Put $F^+:=F^{c=1}$ for the totally real subfield fixed by $c$. Choose a CM type $\Phi$, namely a representative for the coset space $\Gal(F/F^+)\backslash\Hom_{\Q}(F,\ol \Q)$, so that we have the identification
$$F\otimes_{\Q} \R = \textstyle\prod_{\tau\in \Phi} \C,\qquad s\otimes t\mapsto (\tau(s)t).$$
Consider a unitary group $U$ over $F^+$ for the $n$-dimensional Hermitian space $W$ over $F$ (with respect to $F/F^+$) whose signature at each $\tau\in \Phi$ is $(a_\tau,b_\tau)$ with $a_\tau+b_\tau=n$.
We identify $U\times_{F^+} \R = U(a_\tau,b_\tau)$, where the latter is the real unitary group for the form $$((z_i),(w_j))\in \C^n\times \C^n\quad \mapsto \quad \textstyle\sum_{i=1}^{a_\tau} z_i\ol w_j- \textstyle\sum_{i=a_\tau+1}^{n} z_i \ol w_i.$$
Now define $G:=\Res_{F^+/\Q} U$. Consider the $\R$-morphism $h:\Res_{\C/\R}\G_m \ra G_{\R}=\prod_{\tau\in \Phi} U(a_\tau,b_\tau)$, which on $\R$-points sends $z\in \C^\times$ to the block diagonal matrix $\diag((z/\ol z)I_{a_{\tau}},I_{b_\tau})$.
Setting $X$ to be the $G(\R)$-conjugacy class of $h$, we obtain a Shimura datum $(G,X)$.
In this example, we have 
$$\hat G=\prod_{\tau\in \Phi} \GL_n(\C), \qquad \mu(z)=(\diag(z I_{a_\tau},I_{b_\tau})),$$
if we take the upper triangular Borel subgroup in the pinning of $\hat G$. So
\begin{equation}\label{eq:V-mu-unitary}
    V_{-\mu}|_{\hat G}=\otimes_{\tau\in \Phi} \wedge^{a_{\tau}}(\Std^\vee),
\end{equation}
where $\Std^\vee$ denotes the dual of the standard representation of $\GL_n$.
We leave it as an exercise for the reader to describe the $\Gamma_{\Q}$-action on $\hat G$ and compute the reflex field $E$.
In light of \eqref{eq:V-mu-unitary}, an important case is when $U$ has signatures 
$$(a_{\tau_0},b_{\tau_0})=(1,n-1)~\mbox{for some}~\tau_0\in \Phi
\quad\mbox{and}\quad (a_{\tau},b_{\tau})=(0,n)~\mbox{for}~\tau\in \Phi\backslash\{\tau_0\}.$$
The corresponding Shimura varieties have dimension $n$.
In this case, view $F$ as a subfield of $\ol \Q$ (thus also $\C$) via $\tau_0$. Then $E=F$ (unless $E=\Q$ and $n=2$, but this case is excluded by (a)). The representation $V_{-\mu}|_{\hat G}$ is $\Std^\vee$ on the $\tau_0$-component and trivial on the other components. This extends to a representation of $^L G_F = \hat G\rtimes \Gamma_F$ by making $\Gamma_F$ act trivially. So the expected output from the cohomology of Shimura varieties is that, under hypotheses (a), (b), and (c), 
$$\pi\in \mathcal{A}_{\tu{coh}}(G,\Q)=\mathcal{A}_{\tu{coh}}(U,F^+) \quad\mapsto\quad R_{\pi,\mu,\iota} = \rho_{\pi,\iota}^\vee ~\in \mathcal{G}_{\tu{dR}}(n,F)$$
up to twists. 

\begin{remark}
When Theorem \ref{thm:csd-case} was proved, it was based on Kottwitz's result  for PEL-type Shimura data \cite{KottwitzPoints} in the case of certain unitary \emph{similitude} groups. We have chosen to avoid similitude groups to keep the notation and exposition simpler.
The price to pay is that the Shimura datum above is not of PEL or Hodge type but only of abelian type. So the recent work \cite{KSZ} is needed to apply the LKR method; the overall argument mostly stays the same other than that.
\end{remark}

Going back to Theorem \ref{thm:csd-case}, recall that our main interest lies in $\Pi\in \tilde{\mathcal A}^1_{\tu{rac}}(n,F)$. The connection between the two sets $\mathcal{A}_{\tu{coh}}(U,F^+)$ and $\tilde{\mathcal A}^1_{\tu{rac}}(n,F)$ is provided by an automorphic base change due to Clozel and Labesse:
\footnote{For instance, see \cite[Cor.~5.3, Thm.~5.4]{LabesseBaseChange} for both base change and descent results, which are proven under mild technical hypotheses. In practice, one can overcome the hypothesis to obtain Theorem \ref{thm:csd-case}. A general base change and descent result follows from the endoscopic classification for unitary groups in \cite{Mok,KMSW}, but these works are conditional on some papers to be written. (See the introduction of \cite{KMSW} for details.)}
$$\tu{BC}:\mathcal{A}_{\tu{coh}}(U,F^+)\ra \tilde{\mathcal A}^1_{\tu{ra,iso}}(n,F),$$
where the target is by definition the set of isobaric sums $\boxplus_{i\in I} \Pi_i$ such that $\Pi_i\in \tilde{\mathcal A}^1_{\tu{ra}}(n_i,F)$ with $\sum_i n_i=n$. Conversely, they also prove that every $\Pi\in \tilde{\mathcal A}^1_{\tu{rac}}(n,F)$ is the image of some $\pi\in \mathcal{A}_{\tu{coh}}(U,F^+)$ under $\tu{BC}$, namely $\Pi$ admits a ``descent'' $\pi$ (such a $\pi$ is typically far from unique), provided that $U$ is quasi-split at all finite places. In fact there is a parity obstruction (e.g., see \cite[\S2]{ClozelGalois}) for finding such a $U$ with signatures $(1,n-1),(0,n),...,(0,n)$ at $\infty$ if $n$ is even (and $[F^+:\Q]$ is even at the same time), leading to complication in \S~\ref{sss:reality} (b) and (c) below, but let us ignore it until we revisit the issue. 

Thus the construction of $\rho_{\Pi,\iota}$ for Theorem \ref{thm:csd-case} would go as follows:
\begin{equation}\label{eq:ideal-construction}
 \Pi\in \tilde{\mathcal A}^1_{\tu{rac}}(n,F) ~~
\stackrel{\tu{descent}}{\rightsquigarrow}~~\pi\in \mathcal{A}_{\tu{coh}}(U,F^+)  ~~
\stackrel{[H_{\xi}]}{\rightsquigarrow} ~~
\rho^\vee_{\pi,\iota} =: R_{\Pi,\iota} ~\in \mathcal{G}_{\tu{dR}}(n,F).   
\end{equation}
A key point is to compute the restriction of $R_{\Pi,\iota}$ to local Galois groups, at the places $v$ of $F$ where $\Pi$ is unramified. 
The unramified components of $\pi$ can be described from those of $\Pi$ as this is part of the automorphic base change and descent. From this information, the LKR method discussed above gets us to understand $R_{\Pi,\iota}$ at each $v$ to complete the argument.

\subsubsection{Reality}\label{sss:reality}

The reader is cautioned that the argument sketched in the preceding paragraph is not rigorous but idealistic. In fact (a) is somewhat harmless since it is satisfied if $F^+\neq \Q$; even if $F^+=\Q$ there is a standard argument \cite{SorensenPatching} to reduce to the case $F^+\neq \Q$. 
However we do not know (b) a priori, which amounts to a generalized Ramanujan conjecture. In other words, (b) cannot be used as an ingredient. (Rather, it turns out to be a consequence, cf.~\cite{CaraianiLGCellnotp}.) More seriously, (c) cannot be implemented since endoscopy does exhibit itself for unitary groups. 

These difficulties can be overcome with more effort, at the expense of modifying the ideal version \eqref{eq:ideal-construction}.
At first, Theorem \ref{thm:csd-case} was proven by \cite{KottwitzLambdaAdic,ClozelGalois,HarrisTaylor} under a technical assumption that $\Pi$ is essentially square-integrable at a finite place. The assumption is made to avoid endoscopy (i.e., to implement (c)) via the Jacquet--Langlands transfer to the unit group $D^\times$ of a division algebra and by using a unitary group that is a form of $D^\times$. When Ng\^o proved the fundamental lemma \cite{NgoFL} (also thanks to earlier contributions by Cluckers--Loeser, Hales, Waldspurger and others), it became possible to embrace endoscopy and prove Theorem \ref{thm:csd-case} in full generality. This is done in \cite{ClozelHarrisLabesse,ShinGalois,ChenevierHarris}.
Now we explain more precise ideas on how to work without hypotheses (a), (b), and (c).

\subsubsection*{(a) Reducing to the anisotropic case}

As mentioned above, the only nontrivial case is when $F^+=\Q$ (so $F$ is imaginary quadratic over $\Q$); then $G=U$ is a unitary group of $\Q$-rank one so condition (a) fails. 

Suppose that the theorem was proved in all cases for $F^+\neq \Q$.
We want to construct $R_{\Pi,\iota}\in \mathcal{G}_{\tu{dR}}(n,F)$ starting from $ \Pi\in \tilde{\mathcal A}^1_{\tu{rac}}(n,F)$.
The reduction step is based on the fact that there is a ``huge'' supply of real quadratic fields $F^+_0$. 
For each $F^+_0$, set $F_0:=F^+_0 F$. Then we have a diagram
$$
\xymatrix{
\Pi\in \tilde{\mathcal A}^1_{\tu{rac}}(n,F)~ \ar@{-->}[r]^-{?} \ar[d]_-{\tu{base change}~}
& ~R_{\Pi,\iota}\in \mathcal{G}_{\tu{dR}}(n,F)\ar[d]^-{\Res} \\
\Pi_0\in \tilde{\mathcal A}^1_{\tu{rac}}(n,F_0)~\ar[r]_-{\tu{assumed}}
& ~R_{\Pi_0,\iota}\in \mathcal{G}_{\tu{dR}}(n,F_0),
}
$$
where the vertical arrows are the Arthur--Clozel base change and the restriction of Galois representations along $F_0/F$, respectively.
If $R_{\Pi,\iota}$ exists with the correct restriction to local Galois groups, then it should make the diagram commutative. Thus the problem is whether we can construct $R_{\Pi,\iota}$ from the data of $R_{\Pi,\iota}|_{\Gamma_{F_0}}$ for the (infinite) family of real quadratic fields $F^+_0$. This algebraic problem was solved by Blasius--Ramakrishnan \cite[\S4.3]{BlasiusRamakrishnan}; see \cite{SorensenPatching} for a generalization.

The above argument is flexible. For instance, one can analogously reduce to the case when $[F^+:\Q]$ is even by making similar quadratic base changes.

\subsubsection*{(b) Concentration in the middle degree and temperedness}

We begin with a slogan: tempered automorphic representations should contribute to the $\ell$-adic cohomology of Shimura varieties only in the middle degree. Although we cannot use the principle since the relevant case of the Ramanujan conjecture is a priori unknown, the idea is that an approximation towards Ramanujan can still do the job thanks to the Weil conjectures. Moreover the Ramanujan conjecture for relevant cuspidal automorphic representations can be shown as a by-product.

To be more precise, let $\Pi$ and $\pi$ be as in \eqref{eq:ideal-construction}. 
We want to show that the $\pi^\infty$-isotypic part of the cohomology is concentrated in degree $d=\dim \Sh$, so that $R_{\pi,\mu,\iota}$ in \eqref{eq:pi-isotypic} is equal to the alternating sum of the $\pi^\infty$-isotypic part of the cohomology over all degrees (up to the sign $(-1)^d$). The alternating sum is amenable to computation by the fixed-point formula in the context of the LKR method, so we would then have a grip on the local Galois action of $R_{\pi,\mu,\iota}$. 

As we are in the anisotropic case thanks to (a), mod $\fkp$ Shimura varieties are proper (and smooth, since the level subgroup $K_p$ at $p$ is hyperspecial), so that the Weil conjectures apply to the Frobenius eigenvalues on the cohomology. In particular there is no cancellation between different degrees in the alternating sum. Now the LKR method relates the $\Frob_{\fkp}$-eigenvalues in the $\pi^\infty$-isotypic part of cohomology to the Satake parameters of $\pi$ at $\fkp$. Now the Jacquet--Shalika bound on the Satake parameters of $\Pi$ show that they are not far from essentially tempered. 
On the other hand, the $\Frob_{\fkp}$-eigenvalues in degree $i\neq d$ are ``far'' from those in degree $d$ by the Weil conjectures.
From this one can conclude that the $\Frob_{\fkp}$-eigenvalues in the $\pi^\infty$-isotypic part must appear in degree $d$, and further that $\pi_{\fkp}$ is essentially tempered (a fortiori its Satake parameters are Weil numbers of appropriate weight). 

If we start from $ \Pi\in \tilde{\mathcal A}^1_{\tu{rac}}(n,F)$ (rather than $\pi$) then we would prove the Ramanujan conjecture for $\Pi$ by descending to $\pi$ as in the first arrow of \eqref{eq:ideal-construction} and run the preceding argument. This can be made to work if $n$ is odd, but there is a parity obstruction if $n$ is even; see the paragraph above \eqref{eq:ideal-construction}. For even $n$ (when $[F^+:\Q]$ is also even), the Ramanujan conjecture for $\tilde{\mathcal A}^1_{\tu{rac}}(n,F)$ was completed in \cite{ClozelRamanujan,CaraianiLGCellnotp} by switching to the signature  $(2,n-2),(0,n),...,(0,n)$ or $(1,n-1),(1,n-1),(0,n),...,(0,n)$ at $\infty$ to avoid the obstruction and by making a more elaborate argument.

\subsubsection*{(c) Working in the endoscopic setting}
This is the most serious problem. A starting point is the stabilization of the fixed point formula for $\Sh$. This can be done thanks to the fundamental lemma (which also implies that the Langlands--Shelstad transfer exists by \cite{WaldspurgerFLimpliesTransfer}).

The situation is somewhat favorable when $n$ is odd. Then the construction goes as in \eqref{eq:ideal-construction} as there is no obstruction for finding $U$ and descending $\Pi$ to $\pi$. In this case the main problem is to run the trace formula portion of the LKR method to compute the local Galois action of $R_{\Pi,\iota}=\rho^\vee_{\pi,\iota}$. (At primes where $\pi$ is ramified, one uses a variant method following \cite{HarrisTaylor,ShinGalois}.) This involves the stabilized formula for $\Sh$, a comparison with the twisted trace formula for $\GL_n$, and the strong multiplicity one theorem to single out the $\Pi$-isotypic part in the twisted trace formula.

The case of even $n$ is complicated due to the obstruction mentioned in part (b). Here switching to a different signature at $\infty$ as in the last paragraph of (b) is not a good idea in view of \eqref{eq:Rpi-expectation} and \eqref{eq:V-mu-unitary}. Instead, we go from $\Pi\in \tilde{\mathcal A}^1_{\tu{rac}}(n,F)$ to an isobaric sum $\Pi\boxplus \chi$ on $\GL_{n+1}(\A_F)$ by making an auxiliary choice of conjugate self-dual Hecke character $\chi$ of $\GL_1(\A_F)$. Then we consider a unitary group $U$ in $n+1$ variables which is quasi-split at all finite places and signature $(1,n),(0,n+1),...,(0,n+1)$; there is no obstruction since $n+1$ is odd. This unitary group gives rise to Shimura varieties as before. Then the recipe \eqref{eq:ideal-construction} is modified as
$$\Pi\boxplus \chi  ~~
\stackrel{\tu{descent}}{\rightsquigarrow}~~\pi\in \mathcal{A}_{\tu{coh}}(U,F^+)  ~~
\stackrel{[H_{\xi}]}{\rightsquigarrow} ~~
\rho^\vee_{\pi,\iota} =: R_{\Pi,\iota}.$$
One can think of the first arrow as taking \emph{advantage} of endoscopy. The local Galois action at $\fkp$ of $\rho^\vee_{\pi,\iota}$ can be computed using a similar method as for the case of odd~$n$. However the computation shows that the Galois representation $\rho^\vee_{\pi,\iota}$ sometimes corresponds to $\Pi$ and sometimes to $\chi$; we definitely want the former to obtain the correct Galois representation $R_{\Pi,\iota}$. It turns out that the archimedean components of $\Pi$ and $\chi$ control which case should occur. As long as $\Pi$ satisfies a mild regularity condition at $\infty$, it is shown in \cite{ClozelHarrisLabesse,ShinGalois} that $\chi$ can be chosen such that $R_{\Pi,\iota}$ corresponds to $\Pi$ rather than $\chi$. 

In summary, if $n$ is even then we are done as long as $\Pi$ satisfies a mild regularity hypothesis. The remaining case is covered by a congruence argument due to Chenevier--Harris \cite{ChenevierHarris} utilizing the eigenvarieties for definite unitary groups. 

\subsubsection{Complements and further references}\label{sss:self-dual-ref}

There are several expository articles on the Langlands--Kottwitz--Rapoport method with different emphases, including \cite{BlasiusRogawskiMotives,ClozelLKsurvey,MilnePoints,GenestierNgo,ZhuLKmethod}. In the case of bad reduction (when the level subgroup is parahoric but not hyperspecial), one can start from the surveys \cite{RapoportBadReduction,HainesTestFunction}. See also \cite{ShinICM} for a discussion of the LKR method and a variety of related approaches, results, and further references.

The proof of Theorem \ref{thm:csd-case} along the lines of \cite{HarrisTaylor,ShinGalois} is surveyed in \cite{ShinParisBook}, where some aspects of \S~\ref{sss:setup}--\S~\ref{sss:reality} are more detailed and more references are given.
The argument to prove Theorem \ref{thm:csd-case} may be partly simplified by using \cite{ScholzeShin} (which circumvents Mantovan's formula and the works \cite{ClozelHarrisLabesse,ShinGalois}) or \cite{FintzenShin} (which further avoids endoscopy and reduces to the case covered by \cite{ClozelGalois}).

For a stronger version of Theorem \ref{thm:csd-case}, including the local-global compatibility at \emph{all} finite places, refer to \cite[Thm.~2.1.1]{BLGGTpotentialautomorphy} and the references therein due to many. This theorem is complemented by \cite{TaylorComplexConj,TaibiComplexConj,CaraianiLeHung},  describing the image of complex conjugation. In fact \cite{CaraianiLeHung} proves more generally that the $\ell$-adic and torsion Galois representations in Theorem \ref{thm:construction-main} are odd in the sense of Remark \ref{rem:odd}.

From Theorem \ref{thm:csd-case}, one can deduce many cases of part (1) of Conjecture \ref{conj:GLC-ell-adic-G} for classical groups $G$ over totally real fields $F^+$. The additional work consists in establishing a Langlands functoriality for the standard embedding of $^L G$ into a general linear group, and showing that the Galois representations in Theorem \ref{thm:csd-case} satisfy a certain sign condition so as to yield Galois representations valued in the $L$-group (or $C$-group) of $G$. The sign has been pinned down by  \cite{BellaicheChenevierSign}. As for the functoriality, given $\pi\in \mathcal A_{\tu{ac}}(G,F^+)$ we want to find a (not necessarily cuspidal) automorphic representations $\Pi$ of $\GL_N(\A_F)$ for a suitable $N$ such that the map $\pi\mapsto \Pi$ is compatible with the embedding of $^L G$ via the unramified local Langlands at almost all places. Such a map is obtained by \cite{Arthur,Mok} in a very precise form when $G$ is quasi-split, conditionally on some works in preparation. There is a simpler argument to prove a weaker functoriality that is enough for the application to Conjecture \ref{conj:GLC-ell-adic-G}. Either way, the result for $G$ requires $\Pi$ to be regular and algebraic (up to a character twist) in order to apply Theorem \ref{thm:csd-case}. This is the case if $\pi$ is regular in $\mathcal A_{\tu{ac}}(G,F^+)$, and if $\pi$ satisfies a mild extra regularity in the even orthogonal case. See \cite{ShinWeakTransfer} for further details regarding this paragraph.

As in \cite{KretShinGSp,KretShinGSO}, the LKR method can be applied to obtain Conjecture \ref{conj:GLC-ell-adic-G}~(1) when $G$ is $\GSp_{2n}$ or $\GSO_{2n}$ (general symplectic or special orthogonal group) over a totally real field under technical hypotheses. The argument also makes use of the known case of Conjecture~\ref{conj:GLC-ell-adic-G}~(1) for $\Sp_{2n}$ and $\SO_{2n}$in the last paragraph.

\subsection{Perfectoid Shimura varieties}\label{ss:perfectoid-Shimura}

As a preparation for Theorem \ref{thm:construction-main} (2), we need some fundamental results on perfectoid Shimura varieties. We assume familiarity with basics of adic spaces and perfectoid spaces; see \cite{AWS2017perfectoid, berkeley-lectures} for expositions and more references therein. From here until the end of \S~\ref{s:construction} we are in the setting $p=\ell$. This makes our notation close to that of \cite{Sch15}.

Let $\mathcal{X}$ be an adic space and $(\mathcal{X}_i)$ be a filtered projective system of adic spaces with qcqs transition maps. We write 
\begin{equation}\label{eq:perfectoid-limit}
    \mathcal{X}\sim \varprojlim_i \mathcal{X}_i
\end{equation} 
to mean that $\mathcal{X}$ is equipped with maps $f_i:\mathcal{X}\ra \mathcal{X}_i$ compatible with the system $(\mathcal{X}_i)$ such that (i) $(f_i)$ induce a homeomorphism $|\mathcal{X}|\simeq \varprojlim_i |\mathcal{X}_i|$ and (ii) over an affinoid cover of $\mathcal{X}$, the map on the level of Huber rings has dense image. (See \cite[Def.~2.4.1]{ScholzeWeinstein} for the precise definition.) In this article \eqref{eq:perfectoid-limit} will be considered only when $\mathcal{X}$ is perfectoid over a perfectoid field. Then $\mathcal{X}$ is \emph{unique} up to a unique isomorphism \cite[Prop.~2.4.5]{ScholzeWeinstein}.
(Given $(\mathcal{X}_i)$ as above, the limit always exists in the category of diamonds. The point of the definition is that the diamond is actually representable by a perfectoid space.) When a group $\Delta$ acts on each of $\mathcal{X}$ and the system $(\mathcal{X}_i)$, we say that \eqref{eq:perfectoid-limit} is $\Delta$-equivariant if $(f_i)$ are $\Delta$-equivariant. 

\begin{example}\label{ex:perfectoid-limi}
Define $R^+_r:=\Z_p^{\tu{cyc}}\langle t\rangle$, the $p$-adic completion of $\Z_p^{\tu{cyc}}[t]$, and $R_r:=R_r^+[1/p]$ for $r\ge 1$ so that $\Spa(R_r,R^+_r)$ is a closed unit disc. Let $R^+$ denote the $p$-adic completion of $\varinjlim_r R^+_r$ via the ``relative Frobenius'' transition maps $R^+_r\ra R^+_{r+1}$ given by $f(t)\mapsto f(t^p)$, and define $R$ similarly with $R_r$ in place of $R^+_r$. (The transition maps become relative Frobenius maps modulo $p$.)  
Then $\Spa(R,R^+)$ is affinoid perfectoid, known as the ``perfectoid closed unit disc'', and
$$ \Spa(R,R^+) \sim \varprojlim_{r\ge 1} \Spa(R_r,R_r^+).$$
\end{example}

Let $(G,X)$ be a Shimura datum of Hodge type. This means that $(G,X)$ embeds in the Siegel Shimura datum associated with the symplectic similitude group $\GSp_{2g}$ for some $g\in \Z_{\ge 1}$.\footnote{In fact \cite{Sch15} uses a variant of Hodge-type Shimura data in the context of \emph{connected} Shimura data; see \cite[\S4.1]{Sch15}. We ignore this point in our exposition.} 
In this case the Shimura varieties $\Sh_K$ over $\C$ for neat open compact subgroups $K\subset G(\A^{\infty})$ are equipped with closed embeddings into Siegel modular varieties of genus $g$ at suitable levels. (We could work over the reflex field but it is enough for our purpose to stay over a bigger field like $\C$ or $\C_p$.) The Satake--Baily--Borel (a.k.a.~minimal) compactification \cite{BailyBorel} of $\Sh_K$ is denoted by $\Sh_K^*$. The latter is a projective normal variety over $\C$ which is typically singular. Write $\cS_K,\cS_K^*$ for the adic space over~$\C_p$ associated with the base change of $\Sh_K,\Sh_K^*$ to $\C_p$ via a fixed isomorphism $\C\simeq \Qpbar$ and the embedding $\Qpbar\hookrightarrow \C_p$. In fact the canonical model of $\Sh_K$ over the reflex field $E$ extends uniquely to a model of $\Sh^*_K$ over $E$, thus $\cS_K$ and $\cS_K^*$ are defined over $E$ as well.

When a hyperspecial subgroup $K_p^0\subset G_{\Q_p}$ is given (in particular $G_{\Q_p}$ is an unramified group), we have a canonical integral model $\mS_{K_p^0K^p}$ over $\cO_{E,\fkp}$ for $\Sh_{K_p^0K^p}$ at each prime $\fkp$ of $E$ above $p$, cf.~\S\ref{sss:LKR-method}. 
(One needs not assume that $G_{\Q_p}$ is unramified for the analogue of Theorem \ref{thm:perfectoid-Sh} in Remark \ref{rem:HT-period-Hodge-type} to hold true.) In that case, write $\fkS_{K_p^0K^p}$ for the $p$-adic completion of $\mS_{K_p^0K^p}$. Let $\cS^\circ_{K_p^0K^p}$ denote its generic fiber, which is an open subspace of $\cS_K$. We refer to $\cS^\circ_{K_p^0K^p}$ as the good reduction locus; when $(G,X)$ is a Siegel datum, for example, $\Qpbar$-points of $\cS^\circ_{K_p^0K^p}$ is the locus where the parametrized abelian varieties have good reduction.
For an open compact subgroup $K_p\subset K_p^0$, define $\cS^\circ_{K_pK^p}$ to be the preimage of $\cS^\circ_{K_p^0K^p}$ in $\cS_{K_pK^p}$.

The datum $(G,X)$ determines the so-called Hodge  cocharacter $\mu:\G_m\ra G_{\C}$, canonical up to $G(\C)$-conjugacy. The cocharacter $\mu$ gives rise to a parabolic subgroup 
$$P_\mu:=\big\{g\in G: \lim_{t\ra 0} \tu{Ad}(\mu(t)) g~\tu{exists} \big\} \subset G_{\C},$$
well-defined up to conjugacy, 
and to a flag variety $G/P_\mu$ \cite[\S2.1]{CaraianiScholzeCompact}. (The parabolic subgroup $P_\mu$ is opposite to the parabolic subgroup used to construct the flag variety that receives the Borel embedding of the complex Shimura variety.)
This flag variety is a priori defined over $\C$, but in fact it has a model over the reflex field of the Shimura datum. Write $\Fl_G$ for the adic space over $\C_p$ associated with this flag variety. The dependence of $\Fl_G$ on $\mu$ is clear but we omit $\mu$ from the notation for simplicity. When $(G,X)$ is a Siegel Shimura datum, $\Fl_G$ parametrizes maximal isotropic subspaces of the underlying symplectic space.

\begin{theorem}\label{thm:perfectoid-Sh}
Let $(G,X)$ be a Siegel Shimura datum. 
Then there exists a perfectoid space $\cS^*_{K^p}$ over $\C_p$ with a $G(\Q_p)$-action, as well as a $G(\A^{\infty,p})$-action on the projective system $\{\cS^*_{K^p}\}$ as $K^p$ varies, such that we have $G(\Q_p)$-equivariantly
\begin{equation}\label{eq:perfectoid-Sh}
    \cS^*_{K^p}\sim \varprojlim_{K_p\subset G(\Q_p)} \cS^*_{K_p K^p}.
 \end{equation} 
Moreover $\cS^*_{K^p}$ is equipped with the Hodge--Tate period morphism
$$\pi_{\tu{HT}}: \cS^*_{K^p}\ra \Fl_{G}$$
enjoying the following properties.
\begin{enumerate}
    \item $\pi_{\tu{HT}}$ is $G(\Q_p)$-equivariant. As $K^p$ varies, the morphism is also $G(\A^{\infty,p})$-equivariant with respect to the trivial action on $\Fl_{G}$.
    \item There exists a finite cover $(\cU_i)$ of $\Fl_G$ by open affinoid subsets such that $\cV_i:=\pi_{\tu{HT}}^{-1}(\cU_i)$ are affinoid perfectoid satisfying: (i) each $\cV_i$ is the preimage of an affinoid subset of $\cS_{K_pK^p}$ for some open compact $K_p\subset G(\Q_p)$, (ii) a suitable analogue of \eqref{eq:perfectoid-Sh} with $\cV_i$ in place of $\cS^*_{K^p}$ holds true.
    \item There is an ample line bundle $\omega_{\Fl}$ over $\Fl_G$ whose pullback $\omega^*_{K^p}:=\pi_{\tu{HT}}^*\omega_{\Fl}$ is $G(\Q_p)$-equivariantly isomorphic to an automorphic ample line bundle.
\end{enumerate}

\end{theorem}

\begin{remark}\label{rem:HT-period-Hodge-type}
Assume now that $(G,X)$ is a general Shimura datum of Hodge type. The analogous statements are true except that \eqref{eq:perfectoid-Sh} has to be slightly modified by introducing certain \emph{ad hoc} minimal compactifications, which are good enough for applications. 
The relation \eqref{eq:perfectoid-Sh} still holds if it is stated for Shimura varieties without compactification, and it holds for the usual minimal compactifications if $\sim$ in \eqref{eq:perfectoid-Sh} is weakened to an isomorphism of the associated diamonds. For more details and references, see the discussion between Theorem 3.1 and Example 3.2 in \cite{CaraianiICM}. 

\end{remark}

\begin{remark}
It is fruitful to study $\pi_{\tu{HT}}$ with respect to the Newton stratification of $\Fl_G$ to gather further information. Informally speaking, the geometric fibers of $\pi_{\tu{HT}}$ are perfectoid Igusa varieties, and are constant over each individual Newton stratum. This fact is shown and exploited to great effect in \cite{CaraianiScholzeCompact,CaraianiScholzeNonCompact}. We  discuss this more in \S~\ref{ss:geometry of pi_HT}. 
\end{remark}

\begin{proof}[Sketch of proof]
The rest of \S~\ref{ss:perfectoid-Shimura} is devoted to sketching the proof of Theorem \ref{thm:perfectoid-Sh} in the \emph{Siegel case}.
So $G=\GSp_{2g}$ until the end of \S~\ref{ss:perfectoid-Shimura}.
The case of Hodge type can be deduced from this case with further work; refer to \cite[\S4]{Sch15} and \cite[\S2]{CaraianiScholzeCompact}. See \cite{Shen} for an analogous 
result in the case of abelian type. 

Recall the convention for the groups $\Sp_{2g}$ and $\GSp_{2g}$ over $\Z$ from \S~\ref{ss:notation}. Write $\tu{Id}_g$ for the $g\times g$ identity matrix.
For $r\in \Z_{\ge 0}$, we introduce the three congruence subgroups of $K_p^0:=\GSp_{2g}(\Z_p)$:
$$\Gamma(p^r)\subset \Gamma_1(p^r)\subset \Gamma_0(p^r),$$
given by the conditions that $\gamma$ mod $p^m$ is congruent to the forms\footnote{The definition of $\Gamma_0(p^r)$ in \cite{Sch15} also requires that $\det \gamma\equiv 1~\mbox{mod}~p^m$, but we ignore this point.}
$$\begin{pmatrix}
\tu{Id}_g & 0 \\ 0 & \tu{Id}_g
\end{pmatrix}, \quad \begin{pmatrix}
\tu{Id}_g & * \\ 0 & \tu{Id}_g
\end{pmatrix}, \quad \begin{pmatrix}
* & * \\ 0 & *
\end{pmatrix}, \quad \mbox{respectively}.$$

\subsubsection*{Step 0. Hasse invariants}
It is helpful to recall some key facts about the Hasse invariant, denoted by $\tu{Ha}$,
in the setting of modular forms and modular curves, following \cite[\S2.1]{KatzPadic}.
\begin{itemize}
\item[(H1)] $\tu{Ha}$ is a Katz modular form of level 1 and weight $p-1$ over $\F_p$. In particular it is a section of an ample line bundle.
\item[(H2)] $\tu{Ha}$ lifts to a classical modular form $\widetilde{\tu{Ha}}$ if $p\ge 5$. (The lift can be taken to be the Eisenstein series $E_{p-1}$ of level 1 and weight $p-1$.) A power of $\tu{Ha}$ lifts to a classical modular form for $p\le 3$.
\item[(H3)] $\tu{Ha}$ is invertible on the ordinary locus and zero on the supersingular locus. 
\item [(H4)]The $q$-expansion of $\tu{Ha}$ is 1.
\end{itemize}

In general, given an abelian scheme $f:A\ra S$ of dimension $g$ such that $p=0$ in $S$, consider the line bundle $\omega_{A/S}:=\wedge^g (\pi_*\Omega_{A/S})$ over $S$. Put $A^{(p)}:=A\times_{S,\Frob_S}S$, where $\Frob_S:S\ra S$ is the absolute Frobenius. Then the Verschiebung morphism $A^{(p)}\ra A$ induces a map $\omega_{A/S}\ra \omega_{A^{(p)}/S}\simeq \omega_{A/S}^{\otimes p}$, hence a section
$\tu{Ha}(A/S)\in \Gamma(S,\omega_{A/S}^{\otimes (p-1)})$. The Hasse invariant $\tu{Ha}$ over $\mS_{K^0_pK^p,\F_p}$ is obtained by applying this construction to the universal abelian scheme; when $g=1$, we recover the classical Hasse invariant above.
Property (H3) generalizes to the fact that $\tu{Ha}$ is invertible exactly on the ordinary locus. (If $g>1$, there are several intermediate strata between the ordinary and supersingular loci.) 

Even though $\tu{Ha}$ is defined only in characteristic $p$, it can be used to study the $p$-adic formal scheme $\fkS_{K^0_pK^p}$ and its adic generic fiber $\cS^\circ_{K^0_pK^p}$. 
For this, one lifts $\tu{Ha}$ to  $\widetilde{\tu{Ha}}$ over $\fkS_{K^0_pK^p}$, which is always possible locally on the base. This weaker version of (H2) is enough to define $\fkS_{K^0_pK^p}(\epsilon)$ and $\cS_{K^0_pK^p}(\epsilon)$ for rational numbers $0\le \epsilon<1/2$, which are the loci where $\widetilde{\tu{Ha}}$ divides $p^\epsilon$. The locus $\cS_{K^0_pK^p}(\epsilon)$ of $\cS^\circ_{K^0_pK^p}$
should be thought of as the $\epsilon$-neighborhood of the ordinary locus $\cS_{K^0_pK^p}(0)$. Loosely speaking, the points of $\cS_{K^0_pK^p}(\epsilon)$ correspond to abelian varieties that are ``not too supersingluar'' in the sense quantified by $\epsilon$. 

The above construction extends to minimal compactifications. It is classical that $\tu{Ha}$ is defined over  $\mS^*_{K^0_pK^p,\F_p}$ if $g=1$. For $g\ge 2$, $\tu{Ha}$ is extended to $\mS^*_{K^0_pK^p,\F_p}$ by Hartog's extension principle (cf.~\cite[Lem.~3.2.10]{Sch15}; see Remark 3.2.11 therein for an analogy between $\C$-analytic and $p$-adic settings), using the normality of $\mS^*_{K^0_pK^p,\F_p}$ and the fact that the boundary has codimension~$g$. By means of lifts $\widetilde{\tu{Ha}}$ locally over $\fkS^*_{K^0_pK^p}$, we can define\footnote{Actually they are defined over $\Z_p^{\tu{cyc}}$ and $\Q_p^{\tu{cyc}}$ respectively, but we do not keep track of optimal base rings in these notes.}
$$\fkS^*_{K^0_pK^p}(\epsilon)\ra \fkS^*_{K^0_pK^p}~\mbox{over}~\cO_{\C_p} \quad\mbox{and}\quad \cS^*_{K^0_pK^p}(\epsilon) \ra \cS^*_{K^0_pK^p}~\mbox{over}~\C_p.$$
The former is an admissible blow-up in Raynaud's sense, and the latter is an open immersion. For open compact subgroups $K_p\subset K^0_p$, set $\cS_{K_pK^p}(\epsilon)$ to be the preimage of $\cS_{K^0_pK^p}(\epsilon)$ under the projection map $\cS_{K_pK^p}\ra \cS_{K^0_p K^p}$. Define $\cS^*_{K_pK^p}(\epsilon)$ likewise.

\subsubsection*{Step 1. Anti-canonical tower at  $\Gamma_0(p^\infty)$-level}
We want to show that the projective limit $\varprojlim_{K_p}\cS^*_{K_pK^p}$ is perfectoid. In light of Example \ref{ex:perfectoid-limi}, we would hope that the transition maps induce relative Frobenius modulo $p$. While this is overly optimistic, the situation is better if we restrict ourselves to the \emph{anti-canonical} tower over $\cS^*_{K_pK^p}(\epsilon)$ as $K_p$ goes through $\Gamma_0$-level structures. In such a restricted setup, (anti-)canonical subgroups are available to help us show that the mod $p$ picture of such a tower is given by relative Frobenius maps. 

Let $r\in \Z_{\ge 1}$ and $\epsilon\in \Q$ with $0<\epsilon<1/2$.
Let $A$ be an abelian scheme of dimension $g$ over a $p$-adically complete $\cO_{\C_p}$-algebra $R$. 
If $A\times_R R/pR$ is not too supersingular in terms of $r$ (measured by the Hasse invariant), then $A$ 
possesses a unique flat subgroup scheme $C_r\subset A[p^r]$ of order~$p^{rg}$ which becomes the kernel of the $r$-th power Frobenius on $A$ modulo $p^{1-\epsilon}$. (See \cite[Cor.~3.2.6]{Sch15} for the precise statement.) Such a $C_r$ is called the weak canonical subgroup of level~$p^r$. If $A$ contains~$C_1$, i.e., a canonical subgroup of level $p$, a level $p^r$ \emph{anti-canonical subgroup} of $A$ is defined to be a totally isotropic subgroup $D\subset A[p^r]$  of order $p^{rg}$ such that $D[p]\cap C_1=\{0\}$.

Let $\cS_{\Gamma_0(p^{r})K^p}(\epsilon)_{\tu{anti}}$ denote the open subspace of $\cS_{\Gamma_0(p^{r})K^p}(\epsilon)$ on which the universal subgroup of $\cA[p^r]$ coming from the $\Gamma_0(p^r)$-level structure is the anti-canonical subgroup. 
The following commutative diagram is a key:
$$\xymatrix{
\cS_{\Gamma_0(p^{r+1})K^p}(\epsilon)_{\tu{anti}} \ar[r]^-{\sim} \ar[d] & \cS_{K_p^0 K^p}(p^{-r-1}\epsilon) \ar[d] & (A,D_{r+1}) \ar@{|->}[r] \ar@{|->}[d] & A/D_{r+1} \ar@{|->}[d] \\
\cS_{\Gamma_0(p^{r})K^p}(\epsilon)_{\tu{anti}} \ar[r]^-{\sim}  & \cS_{K_p^0 K^p}(p^{-r}\epsilon)
& (A,D_r)=(A,D_{r+1}[p^r]) \ar@{|->}[r] & A/D_r,
}$$
where the diagram on the right gives the moduli-theoretic description of the maps in terms of principally polarized abelian varieties, together with anti-canonical subgroups in the left column. The two vertical maps are natural projections. The horizontal maps are isomorphisms as the inverse maps can be given as the quotient maps by suitable canonical subgroups, cf.~\cite[Thm.~3.2.15.(ii)]{Sch15}. The projective system in the left column as $r\ra\infty$ is the anti-canonical tower. 

The use of anti-canonical subgroups is significant in two ways. First, $\epsilon$ does not shrink as $r\ra\infty$ in the anti-canonical tower. (In a similar construction via canonical subgroups, $\epsilon$ tends to $0$.) Second, the right vertical map $A/D_{r+1}\mapsto A/D_r$ is a lift of relative Frobenius. Analogously with Example \ref{ex:perfectoid-limi}, this essentially tells us that a perfectoid limit $\cS_{\Gamma_0(p^{\infty})K^p}(\epsilon)_{\tu{anti}}$ of the anti-canonical tower exists, satisfying $ \cS_{\Gamma_0(p^{\infty})K^p}(\epsilon)_{\tu{anti}}\sim \varprojlim_{r} \cS_{\Gamma_0(p^{r})K^p}(\epsilon)_{\tu{anti}}$.

The anti-canonical tower extends to the tower of $\cS^*_{\Gamma_0(p^{r})K^p}(\epsilon)_{\tu{anti}}$ over the compactifications; Hartog's extension principle is used if $g\ge 2$ (cf.~\cite[Lem.~3.2.10]{Sch15}). Essentially the same argument as above shows the existence of a perfectoid limit
$$\cS^*_{\Gamma_0(p^{\infty})K^p}(\epsilon)_{\tu{anti}}\sim \varprojlim_{r} \cS^*_{\Gamma_0(p^{r})K^p}(\epsilon)_{\tu{anti}}.$$
Moreover it is affinoid perfectoid. The intuitive reason is that each $\cS^*_{\Gamma_0(p^{r})K^p}(\epsilon)_{\tu{anti}}$ is isomorphic to $\cS^*_{K^p}(p^{-r}\epsilon)$, and the latter is affinoid as it is the $p^{-r}\epsilon$-neighborhood given in terms of a section of an ample line bundle $\omega^{\otimes p^r(p-1)}$, namely a lift of the $p^r$-th power of $\tu{Ha}$, cf.~(H1) above. (Such a lift exists globally for $r\gg 1$ by ampleness of $\omega$.) A precise argument is made from computing the tilts in characteristic $p$ \cite[Cor.~3.2.20]{Sch15}.

\subsubsection*{Step 2. Going up to $\Gamma(p^\infty)$-level}
We have shown that $\varprojlim_r \cS^*_{\Gamma_0(p^r) K^p}(\epsilon)_{\tu{anti}}$ is perfectoid. Since $\Gamma_0(p^r)$ do not form a neighborhood basis at $1$ in $G(\Q_p)$ but $\Gamma(p^r)$ do, we would like to go from the $\Gamma_0(p^\infty)$-level all the way up to the $\Gamma(p^\infty)$-level. This is done in two steps via $\Gamma_1(p^\infty)$.

Going from $\Gamma_1(p^\infty)$ to $\Gamma(p^\infty)$ is relatively easy. Since 
$$\cS^*_{\Gamma(p^r) K^p}(\epsilon)_{\tu{anti}}\ra \cS^*_{\Gamma_1(p^r) K^p}(\epsilon)_{\tu{anti}}$$
are finite \'etale for all $r\ge 1$, one can appeal to almost purity, which means that an adic space that is finite \'etale over an affinoid perfectoid space is also affinoid perfectoid.
The hard part is to pass from $\Gamma_0(p^\infty)$ to $\Gamma_1(p^\infty)$. Since 
$$\cS_{\Gamma_1(p^r) K^p}(\epsilon)_{\tu{anti}}\ra \cS_{\Gamma_0(p^r) K^p}(\epsilon)_{\tu{anti}}$$
is finite \'etale, there is no problem with obtaining a perfectoid limit over the interior. The main difficulty is how to extend to the boundary since the map between compactifications is ramified at the boundary. The idea is to tilt the interior to characteristic $p$, extend it to an affinoid perfectoid space including the boundary, and prove that this extension is the correct tilt. The proof involves Hartog's extension principle in characteristic 0 and a $p$-adic version of Riemann's Hebbarkeitssatz (``bounded~ functions have removable singularities'') in characteristic $p$.

The output is an affinoid perfectoid space $$\cS^*_{K^p}(\epsilon)_{\tu{anti}}\sim\varprojlim_r \cS^*_{\Gamma(p^r) K^p}(\epsilon)_{\tu{anti}}.$$

\subsubsection*{Step 3. Topological Hodge--Tate morphism}
Motivated by the definition of \eqref{eq:perfectoid-limit}, we define the topological space $$|\cS^*_{K^p}|:=\varprojlim_{K_p\subset G(\Q_p)} |\cS^*_{K_pK^p}|,$$
which is naturally equipped with a continuous $G(\Q_p)$-action. One of our goals is to upgrade $|\cS^*_{K^p}|$ to a perfectoid space. As a preparation and also as a step towards the Hodge--Tate period map, cf.~Step 5 below, we construct a continuous $G(\Q_p)$-equivariant map of topological spaces $$|\pi_{\tu{HT}}|: |\cS^\circ_{K^p}| \ra |\Fl_G|,$$
described pointwise as follows.
Each $x\in |\cS^\circ_{K^p}|$ can be represented by a $\Spa(\kappa,\kappa^+)$-point for a complete nonarchimedean extension $\kappa\supset \C_p$, which  corresponds to a principally polarized abelian variety $A_x$ over $\kappa$ with a symplectic isomorphism $T_p A_x \simeq \Z_p^{2g}$ (matching symplectic forms up to scalars). Then $x$ is sent to the Hodge--Tate filtration
$$0\ra \Lie A_x (1) \ra T_p A_x\otimes_{\Z_p} C_x \simeq C_x^{2g} \ra  (\Lie A^\vee_x)^\vee\ra 0, $$
which determines a point of $|\Fl_G|$. Continuity of $|\pi_{\tu{HT}}|$ essentially follows from the fact that the Hodge--Tate filtration works well in families. In fact, the same construction defines $|\pi_{\tu{HT}}|$ on a larger domain than the good reduction locus $\cS^\circ_{K^p},$
namely $|\cS^*_{K^p}|$ minus the boundary, the point being that the Hodge--Tate filtration belongs in characteristic zero.

\subsubsection*{Step 4. Spreading the perfectoid property}
Fix $0\le \epsilon<1/2$.
Recall from Step 2 that the open subset $|\cS_{K^p}^*(\epsilon)_{\tu{anti}}|$ of $|\cS_{K^p}^*|$ is affinoid perfectoid. In fact 
$$|\cS_{K^p}^*(\epsilon)|=G(\Z_p)\cdot |\cS_{K^p}^*(\epsilon)_{\tu{anti}}|$$
so there is a perfectoid space $\cS_{K^p}^*(\epsilon)$ whose underlying space is $|\cS_{K^p}^*(\epsilon)|$. It is covered by finitely many translates of $\cS_{K^p}^*(\epsilon)_{\tu{anti}}$. (Finiteness comes from continuity of the $G(\Z_p)$-action and the fact that an open subgroup of $G(\Z_p)$ has finite index.)
Denote by $\cS_{K^p}(\epsilon)$ the open subspace of $\cS_{K^p}^*(\epsilon)$ away from the boundary.

To promote $|\cS_{K^p}^*|$ to the desired perfectoid space in Theorem \ref{thm:perfectoid-Sh}, it is thus enough to show that finitely many translates of $\cS_{K^p}^*(\epsilon)$ under the $G(\Q_p)$-action cover $|\cS_{K^p}^*|$.
The starting point is the observation that the ordinary locus $\cS_{K^p}(0)$ maps to $\Fl_G(\Q_p)$ under $|\pi_{\tu{HT}}|$. 
Now it is important to take $\epsilon>0$ to work with a strict neighborhood of the ordinary locus. Then
one can show the existence of an open subset $U\subset \Fl_G$ containing $\Fl_G(\Q_p)$ such that $$|\pi_{\tu{HT}}|^{-1}(U)\subset |\cS^*_{K^p}(\epsilon)|.$$
On the other hand, $G(\Q_p)U=\Fl_G$ as this is true for every open subset containing $\Fl_G(\Q_p)$; the proof uses an explicit open affinoid covering $\Fl_G=\cup_J \Fl_{G,J}$ coming from the Pl\"ucker embedding, indexed by cardinality $g$ subsets  $J\subset\{1,2,...,2g\}$. Since $\Fl_G$ is quasi-compact, a finite subset $\{g_i\}_{i\in I}\subset G(\Q_p)$ can be chosen such that
$$\Fl_G = \cup_{i\in I} g_i U .$$
The two displayed facts together with the $G(\Q_p)$-equivariance of $|\pi_{\tu{HT}}|$ imply that $|\cS^*_{K^p}|=\cup_{i\in I} g_i |\pi_{\tu{HT}}|^{-1}(U)$. (The boundary is taken care of by a topological argument.) Hence
$$|\cS^*_{K^p}|=\cup_{i\in I} g_i |\cS^*_{K^p}(\epsilon)|.$$

\subsubsection*{Step 5. Finishing the construction of $\pi_{\tu{HT}}$}

Now that we have a perfectoid space $\cS_{K^p}$, the topological map $|\pi_{\tu{HT}}|$ can be upgraded to a map of adic spaces
$$\pi_{\tu{HT}}: \cS^*_{K^p}\ra \Fl_G$$
by largely repeating the same idea for constructing $|\pi_{\tu{HT}}|$ and extending the map to the boundary by a $p$-adic Hebbarkeitssatz.

\subsubsection*{Step 6. Further properties of $\pi_{\tu{HT}}$}

We briefly comment on verifying properties (2) and (3) of $\pi_{\tu{HT}}$ in Theorem \ref{thm:perfectoid-Sh}.
For (2), it turns out that the covering  $(\Fl_{G,J})$ of $\Fl_G$ from Step 4 works.
The argument again utilizes the outcome of Step 2 and the $G(\Q_p)$-equivariance of $\pi_{\tu{HT}}$. 
For (3), if $W_{\Fl_G}\subset \cO^{2g}_{\Fl_G}$ denotes the universal totally isotropic subspace, then it follows from the construction of $\pi_{\tu{HT}}$ that there is a $G(\Q_p)$-equivariant isomorphism of locally free rank $g$ modules
$$\Lie \cA_{K^p}\simeq \pi_{\tu{HT}}^* W_{\Fl_G}$$
away from the boundary. Property (3) is obtained by taking the dual of the top exterior power and extending to the boundary.
\end{proof}

\subsection{Construction of torsion Galois representations}\label{ss:construction-torsion-Galois}

\begin{theorem}\label{thm:construction-torsion}
Theorem \ref{thm:construction-main} (2) holds true. We also have a suitable analogue for mod $p^m$ coefficients for $m\ge 1$, in terms of Chenevier's determinants \cite{ChenevierDet}.
\end{theorem}

A more precise statement is as follows. Denote by $\T^S_{\GL_n}$ the abstract Hecke algebra for $\GL_n$ over $F$ away from $S$ as in \S~\ref{sss:BG-conj}. Let $K'\subset \GL_n(\A_F^\infty)$ be an open compact subgroup as in \S \ref{sss:BG-conj}. Then $H^*(Y_{\GL_n,K'},\Z/p^m\Z)$ is naturally a $\T^S_{\GL_n}$-module. Write 
$$\T^S_{\GL_n}(K',m):=\tu{image}(\T^S_{\GL_n}\ra \End_{\Z/p^m\Z}(H^*(Y_{\GL_n,K'},\Z/p^m\Z))).$$
Then the theorem asserts that there exists an $n$-dimensional determinant $D$ of $\Gamma_F$ which is unramified away from $S$ and valued in $\T^S_{\GL_n}(K',m)/J$ for a nilpotent ideal $J\subset \T^S_{\GL_n}(K',m)$ of bounded nilpotence exponent, such that the characteristic polynomial of $D$ at each $\Frob_v$, $v\notin S$, is given by $H_v(x)/J$. (Throughout this subsection, we use the representation-theoretic notion of determinants from \cite{ChenevierDet}. See therein for a comparison with similar notions such as pseudocharacters.) A technical remark is that it is advantageous to deal with Hecke algebras and Hecke modules mod $p^m$ rather than Hecke eigensystems mod $p^m$. 

The existence of such a $D$ subsumes Theorem \ref{thm:construction-main} (2). Indeed, $\fkm\in \HE^S_{\tu{coh}}(n,F)_{\Flbar}$ corresponds to a morphism $f:\T^S_{\GL_n}(K',1)/J\ra \Fpbar$, so $f\circ D$ is an $n$-dimensional determinant of $\Gamma_F$ valued in $\Fpbar$ with the correct Frobenius characteristic polynomials away from~$S$. Since a determinant valued in an algebraically closed field determines a unique isomorphism class of semisimple  $n$-dimensional representations (and vice versa) by \cite[Thm.~A]{ChenevierDet}, we obtain the desired $\ol\rho_{\fkm}$ as in Theorem \ref{thm:construction-main} (2). As done in \cite{Sch15} one can keep track of topology to check that $D$ and $f \circ D$ are continuous; it follows from this that $\ol \rho_{\fkm}$ is a continuous representation.

\begin{proof}[Sketch of proof]
We concentrate on the case where $F$ is a CM field and employ unitary Shimura varieties.
If $F$ is totally real, then the argument below can be adapted to the setting of Siegel modular varieties for the group $G=\Res_{F/\Q}\Sp_{2n}$, cf.~\cite[\S5.1]{Sch15}.

Recall that $F^+$ denotes the maximal totally real subfield of $F$. Let $U_{2n}$ denote the quasi-split unitary group in $2n$ variables relative to the quadratic extension $F/F^+$.
For $G:=\Res_{F^+/\Q} U_{2n}$, there is a standard choice for $X$ to make $(G,X)$ a Shimura datum.
One has the tower of connected Shimura varieties $\Sh_{K}$ for $G$ over $\Q(\mu_\infty)$ indexed by neat open compact subgroups $K=\prod_{v\neq \infty} K_v\subset G(\A^\infty)$.
(This differs from the convention in the proof part of \S~\ref{ss:perfectoid-Shimura}, where $G$ was $\GSp_{2g}$ and $\Sh_K$ denoted the Siegel modular varieties.) 
Let $S$ be a finite set of places as in \S~\ref{sss:BG-conj}, which in particular contains $\{p,\infty\}$, and keep $K_v$ hyperspecial at $v\notin S$.
Write $\cI$ for the ideal sheaf of the boundary of $\cS^*_{K_pK^p}$, and $\cO^+$ for the $+$-part of the structure sheaf of $\cS^*_{K_pK^p}$. Set $\cI^+:=\cI\cap \cO^+$ over $\cS^*_{K_pK^p}$. We still denote by $\cI^+$ the pullback of $\cI^+$ to $\cS^*_{K^p}$. This sheaf may be thought of as the sheaf of cuspforms. Define
$$\tilde H^i(\Sh_{K^p},\Z/p^m\Z):=\varinjlim_{K_p} H^i(\Sh_{K_p K^p},\Z/p^m\Z),\quad i\ge 0.$$
The tilde indicates that it is the mod $p^m$ version of the \emph{completed} cohomology.

\subsubsection*{(0) Overview}

The basic strategy is to reduce to the case of \S~\ref{ss:csd} where Galois representations are already constructed, via $p$-adic congruences. 
The proof progresses by studying various cohomology spaces in order:
$$H^*(Y_{\GL_n,K'},\Z/p^m\Z) ~\stackrel{(1)}{\rightsquigarrow}~ \tilde H^*(\Sh_{K^p},\Z/p^m\Z) ~\stackrel{(2)}{\rightsquigarrow}~
H^*(\cS^*_{K^p},\cI^+/p^m)$$
$$ ~\stackrel{(3)}{\rightsquigarrow}~ H^0(\cV,\cI^+/p^m)
~\stackrel{(4)}{\rightsquigarrow}~ 
H^*(\cS^*_{K_pK^p},\omega^{\otimes mk}_{K_pK^p}\otimes \cI),
$$
where $\cV$ is an open affinoid perfectoid subspace of $\cS^*_{K^p}$ (see (3) below), and $\omega_{K_pK^p}$ is the automorphic ample line bundle at a finite level $K_p\subset G(\Q_p)$ whose pullback to the infinite level is as in property (3) of Theorem \ref{thm:perfectoid-Sh}.
The five cohomology groups are computed on their respective topological spaces (the first two of which are real; the next three are $p$-adic). It is important to compute the second and third groups also in the \'etale topology to make the transitions (2) and (3).

Step (1) is a reduction of the problem about locally symmetric spaces for $\GL_n$ over $F$ (which are not Shimura varieties if $n>1$) to one about Shimura varieties associated with $G=\Res_{F^+/\Q} U_{2n}$. The next steps (2)--(4)  establish $p$-adic congruences, namely that the mod $p^m$ (and $p$-adically completed) cohomology of the Shimura varieties can be $p$-adically approximated by classical cuspforms. 

In summary, if we keep a careful track of Hecke algebras along the way, then 
the above steps reduce the construction of Galois representations (or rather determinants of the Galois group) for $H^*(Y_{\GL_n,K'},\Z/p^m\Z)$ to that for $H^0(\cS^*_{K_p K^p},\omega_{K_pK^p}^{\otimes mk}\otimes \cI)$.
The latter consists of classical cuspidal automorphic forms on $U_{2n}$ over $F^+$, so we do have associated Galois representations from the conjugate self-dual case treated in \S~\ref{ss:csd}. This completes the proof.
Now we give more details on each step below.

\subsubsection*{(1) Borel--Serre compactification}\label{sss:Borel-Serre}

The Borel--Serre compactification $\Sh^{\tu{BS}}_K$ is a compact real manifold with corners such that $\Sh_K\hookrightarrow\Sh^{\tu{BS}}_K$ is a homotopy equivalence. 
(See \cite[\S5.2]{Sch15}, cf.~the original work \cite{BorelSerre} and a survey \cite[\S4]{GoreskyClay}.) 
Write $\partial \Sh^{\tu{BS}}_K$ for the complement of $\Sh_K$. There is a long exact sequence
$$
\cdots \ra H^i_c(\Sh_K,\Z/p^m \Z) \ra H^i(\Sh^{\tu{BS}}_K,\Z/p^m\Z) \ra H^i(\partial \Sh^{\tu{BS}}_K,\Z/p^m\Z) \ra \cdots,
$$
where the middle term is identified with $H^i(\Sh_K,\Z/p^m\Z)$.

The boundary strata in $\partial \Sh^{\tu{BS}}_K$ are indexed by the conjugacy classes of $\Q$-rational proper parabolic subgroups of $U_{2n}$.
From the basic fact that $\Res_{F/F^+}\GL_n$ is a Levi subgroup of such a parabolic subgroup, it follows that the ``parabolic induction'' (for modules over Hecke algebras, as we are working at fixed levels) of $H^*(Y_{\GL_n,K'},\Z/p^m\Z)$ contributes to $H^*(\partial \Sh^{\tu{BS}}_K,\Z/p^m\Z)$.
According to the long exact sequence, such a contribution to each $H^i(\partial \Sh^{\tu{BS}}_K,\Z/p^m\Z)$ is captured by $H^i(\Sh_K,\Z/p^m\Z)$ and $H^{i+1}_c(\Sh_K,\Z/p^m\Z)$.
So we are essentially reduced to attaching mod $p^m$ Galois representations to $H^*(\Sh_K,\Z/p^m\Z)$ and $H^*_c(\Sh_K,\Z/p^m\Z)$.\footnote{More precisely, we construct a pseudo-determinant valued in the Hecke algebra corresponding to each cohomology group. A priori we get a $2n$-dimensional pseudo-determinant as $\Sh_K$ are associated with $U_{2n}$. It takes an extra twisting argument to obtain the desired $n$-dimensional pseudo-determinant from this. See \cite[\S5.3]{Sch15}, cf.~\cite[\S7]{HLTT}.} Thanks to Poincar\'e duality, it is enough to deal with $H^*_c(\Sh_K,\Z/p^m\Z)$, or 
$$\tilde H^i_c(\Sh_{K^p},\Z/p^m\Z)$$ 
after taking the limit over $K_p$. 

\subsubsection*{(2) A comparison isomorphism}
Write $j:\Sh_{K_pK^p}\hookrightarrow \Sh^*_{K_pK^p}$ for the open immersion. The main point of (2) is the following isomorphisms which are equivariant for the Hecke action away from $p$, where the tensor products are over $\Z/p^m\Z$:
$$ \tilde H^*_c(\Sh_{K^p},\Z/p^m\Z)\otimes \cO_{\C_p}/p^m ~\simeq~ \tilde H^*(\Sh^*_{K^p},j_! \Z/p^m\Z) \otimes \cO_{\C_p}/p^m$$
$$\stackrel{a}{\simeq}~ H^*_{\tu{\'et}}(\cS^*_{K^p},j_!\cO^+/p^m) ~\simeq~ H^*_{\tu{\'et}}(\cS^*_{K^p},\cI^+/p^m).$$
The first isomorphism is formal. The second map is a comparison \cite[Thm.~3.13]{ScholzeCDM} for \'etale cohomology between a proper (not necessarily smooth) scheme over $\C_p$ and its adification with constructible sheaf coefficients. (The theorem for $\F_p$-sheaves therein is easily extended to the case of $\Z/p^m\Z$-sheaves.) Actually it is only an \emph{almost} isomorphism, denoted by the symbol $\stackrel{a}{\simeq}$, but we ignore the technical point as it does not affect our goal. 

\subsubsection*{(3) Hecke-stable  affinoid perfectoid \v{C}ech covering}

Now it truly matters to pass to infinite-level at $p$.
We employ the affinoid perfectoid covering $(\cV_i)$ of $\cS^*_{K^p}$ from Theorem \ref{thm:perfectoid-Sh} to compute $H^*_{\tu{\'et}}(\cS^*_{K^p},\cI^+/p^m)$ as the \v{C}ech cohomology of $\cS^*_{K^p}$ as a topological space.
Here we appeal to cohomological vanishing of affinoid perfectoid spaces as in \cite[Prop.~6.14, 7.13]{ScholzePerfectoid}, and also the fact that the boundary of the perfectoid space $\cS^*_{K^p}$ is strongly Zariski closed in each $\cV_i$. 

The outcome is that $H^*_{\tu{\'et}}(\cS^*_{K^p},\cI^+/p^m)$ can be computed by a \v{C}ech complex whose terms consist of $$H^0(\cV,\cI^+/p^m),$$
where $\cV$ stands for finite intersections made from the covering $(\cV_i)$. Since each $\cV_i$ is stable under the Hecke correspondences away from $S$ in view of part (1) of Theorem \ref{thm:perfectoid-Sh}, so is $\cV$. Thus the Hecke algebra action away from $S$ is encoded by each $H^0(\cV,\cI^+/p^m)$.

\subsubsection*{(4) Fake Hasse invariants}

The point of this step is to extend the sections of $\cI^+/p^m$ on $\cV$ to the whole Shimura variety and also lift to characteristic 0, at the expense of multiplying the coefficient sheaf by a power of an automorphic line bundle. This is achieved by introducing fake Hasse invariants, mimicking the classical Hasse invariants.
The name ``fake'' suggests that, unlike the classical Hasse invariants, they do not arise as sections of vector bundles over typical integral models of Shimura varieties at finite level. 

In the case of modular forms, when a mod $p^m$ modular form is multiplied by a sufficiently divisible power of the Hasse invariant (or its lift to characteristic 0), the poles outside the ordinary locus are removed while the $q$-expansion mod $p^m$ is unchanged. (This argument appears in the proof of \cite[Thm.~4.5.1]{KatzPadic}, for instance.) Moreover the resulting modular form lifts to characteristic 0 by a cohomological vanishing theorem for ample line bundles.

Fake Hasse invariants are designed to play similar roles. They are sections of $\omega_{K^p}$ pulled back via $\pi_{\tu{HT}}$ from explicit sections of $\omega_{\Fl}$, via Theorem \ref{thm:perfectoid-Sh} (3). 
One can choose ``exotic'' formal integral models of $\cS^*_{K^p}$ and $\omega_{K^p}$ at finite level $K_p\subset G(\Q_p)$ where these sections mod $p^m$ are defined. (See footnote 19 of \cite[p.1030]{Sch15} for ``exotic''.)
The key properties of fake Hasse invariants are as follows, cf.~(H1)--(H4) in Step 0 of the proof of Theorem \ref{thm:perfectoid-Sh}, besides abundance of fake Hasse invariants thanks to Theorem \ref{thm:perfectoid-Sh} (2) (i.e., the Hodge--Tate period map is ``affinoid'').
\begin{itemize}
\item[(i)] They are sections of ample line bundles.
    \item[(ii)] Multiplying them commutes with the prime-to-$p$ Hecke action.
    \item[(iii)] The formal integral model of $\cV$ at level $K_p$ is the invertible locus of a suitable fake Hasse invariant.
\end{itemize}
With these inputs, one can descend to  finite levels $K_p$ at $p$ and mimic the classical arguments, by multiplying a sufficiently divisible power of fake Hasse invariants, to extend sections from $\cV_{K_p}$ to $\cS^*_{K_pK^p}$ and lift the sections from mod $p^m$ to characteristic 0 coefficients, compatibly with the Hecke actions away from $p$. This completes Step (4).
\end{proof}

\subsubsection*{Complements and further references}\label{sss:complements-further-refs}
Regarding \S\S~\ref{ss:perfectoid-Shimura}--\ref{ss:construction-torsion-Galois}, we recommend the survey articles 
\cite{ScholzeICM2014, ScholzeICM2018, ScholzeJapanese,MorelBourbaki, WeinsteinBAMS} as well as Caraiani's chapter in \cite{AWS2017perfectoid}. 

Varma \cite{VarmaThesis} has strengthened Theorem \ref{thm:construction-main} (1) by showing the local-global compatibility with characteristic 0 coefficients (up to semisimplification) at all primes away from $\ell$, not just at unramified places. See \cite[\S\S3--5]{10authors} for some local-global compatibility results at all primes including those above $\ell$ (the ``$\ell=p$ case'', which is most subtle), where Galois representations are valued in suitable Hecke algebras; this requires understanding of the torsion setting.

Scholze's method was adapted by Pilloni--Stroh \cite[\S3]{PilloniStrohCoherent} to prove the analogue of Theorem \ref{thm:construction-main} (1) for automorphic representations contributing to the \emph{coherent cohomology} of certain Shimura varieties of Hodge type via analogues of Hasse invariants.  Since such representations can be \emph{non-regular} (e.g., the archimedean components can be non-degenerate limits of discrete series), their theorem applies to new cases that are not covered by Theorem \ref{thm:construction-main} (1). Getting around perfectoid Shimura varieties, Boxer and  Goldring--Koskivirta \cite{BoxerThesis,GoldringKoskivirta} (independently of each other) develop different analogues of Hasse invariants and their applications; see therein for further references. The results of \cite{PilloniStrohCoherent} are also obtained by \cite{GoldringKoskivirta}.

\section{The Calegari--Geraghty method for Betti cohomology}\label{s:calegari-geraghty}

The goal of this section is to give a rough outline of the Calegari--Geraghty method in the Betti setting for the group $\mathrm{GL}_n$ and to discuss its prerequisites. The Calegari--Geraghty method is a vast extension of the Taylor--Wiles method for proving modularity lifting, first proposed in~\cite{CG18}.
There are two versions of the Calegari--Geraghty method in the literature, both introduced in~\cite{CG18}: one that applies to the Betti cohomology of locally symmetric spaces and one that applies to the coherent cohomology of Shimura varieties. 
In both cases, it is essential to apply the Taylor--Wiles patching method to complexes rather than to individual cohomology groups. 

The remarkable feature of the Betti version is that it applies to locally symmetric spaces that are not Shimura varieties - this is the case on which we will focus. The coherent version applies to coherent cohomology groups of Shimura varieties, particularly in non-regular weight. This version has been spectacularly applied in the work of Boxer--Calegari--Gee--Pilloni to prove the potential modularity of abelian surfaces over totally real fields~\cite{BCGP}. Another essential ingredient for this result comes from higher Hida theory, which was introduced by Pilloni~\cite{pilloni-HH} and is now being further developed by Boxer and Pilloni (see~\cite{BoxerPilloni} for the modular curve case of higher Hida theory). 
There are many beautiful ideas here - see, for example,~\cite{Calegari-CDM} for a survey focused on this topic.    

\subsection{Conjectures on the cohomology of locally symmetric spaces}\label{ss:prerequisites to CG}

For this section, we let our connected reductive group $G$ be $\GL_n$ and we work over an arbitrary number field $F$. We let $\ell$ be a prime and $S$ be a finite set of places of $F$ containing all the ramified places and all the places above $\ell$ and $\infty$. For a sufficiently small compact open subgroup $K\subset G(\A^\infty_F)$ as in \S~\ref{sss:BG-conj}, we have the locally symmetric space 
\[
Y_K = G(F)\backslash G(\A_F)/K K'_\infty,
\] 
which is a smooth Riemannian manifold. This locally symmetric space does not arise from a Shimura variety except in very special cases. 

\begin{example}\leavevmode
\begin{enumerate}
\item Let $n=2$ and $F$ be a totally real field. If $F=\Q$, the locally symmetric space $Y_K$ arises from a modular curve defined over $\Q$. If $[F:\Q]>1$, the locally symmetric space $Y_K$ is not a Shimura variety itself. However, it is closely related to one, since it is essentially a torus bundle over a Hilbert modular variety. Hilbert modular varieties are Shimura varieties of abelian type. 
\item Let $n=2$ and $F$ be an imaginary quadratic field. Then $G(F_\infty) = \GL_2(\C)$ and $K'_\infty = \mathrm{U}_2(\R)\cdot \R_{>0}$. The quotient $G(F_\infty)/K'_{\infty}$ can be identified with the $3$-dimensional hyperbolic space $\mathbb{H}^3$. Therefore, the locally symmetric space $Y_K$ is an arithmetic hyperbolic $3$-manifold as it can be obtained from a disjoint union of finitely many quotients of $\mathbb{H}^3$ by congruence subgroups of $\GL_2(\cO_F)$. In particular, $Y_K$ does not admit a complex structure. The manifolds $Y_K$ are also called \emph{Bianchi manifolds}.
\item More generally, if $n\geq 2$ and $F$ is not totally real, or if $n\geq 3$ and $F$ is an arbitrary number field, the locally symmetric spaces $Y_K$ cannot be directly related to Shimura varieties. This follows from a classification of groups that admit a Shimura datum.   
\end{enumerate}
\end{example}

\medskip 

\subsubsection{Existence of Galois representations} The Calegari--Geraghty method relies on a deep understanding of the Betti cohomology of the locally symmetric spaces $Y_K$, including a formulation of the Buzzard--Gee conjecture for torsion classes in the cohomology. To make this more precise, we introduce certain local systems on $Y_K$. 

We will use the following notation throughout this subsection. Let $E/\Q_{\ell}$ be a finite extension which will be our field of coefficients, with ring of integers $\cO$, uniformizer $\varpi$, and residue field $k$. We always assume that $E$ is large enough, for example so that it contains the image of every embedding of $F$ into $\overline{\Q}_{\ell}$. 

We let $\lambda$ be a highest weight vector for $\mathrm{Res}_{F/\Q}\GL_n$. In other words, we take $\lambda \in (\Z^n)^{\Hom_{\Q}(F,E)}$ such that $\lambda_{\tau, 1}\geq \lambda_{\tau, 2}\geq \dots \geq \lambda_{\tau, n}$ for every embedding $\tau\in \Hom_{\Q}(F, E)$. We let $\sigma_{\lambda}$ denote the irreducible algebraic representation of $\mathrm{Res}_{F/\Q}G$ over $\overline{\Q}_{\ell}$ of highest weight $\lambda$. Since $E$ is assumed large enough, we may assume that $\sigma_{\lambda}$ is defined over $E$. We let $\sigma_{\lambda}^\circ$ denote a $\prod_{v\mid \ell}\GL_n(\cO_{F,v})$-stable lattice in $\sigma_{\lambda}$. 

We assume that the compact open subgroup $K$ is of the form $\prod_{v} K_v$, where $v$ runs over finite places of $F$ and $K_v\subset \GL_n(\cO_{F,v})$ is a compact open subgroup. We let $K$ act on $\sigma_{\lambda}^\circ$ via the projection to $\prod_{v\mid \ell}K_v$, pull back $\sigma^\circ_{\lambda}$ to a $K$-equivariant sheaf on $G(F)\backslash G(\A_F)/K'_{\infty}$ and the descend to a local system of $\cO$-modules $\cV_{\lambda}$ on $Y_K$. We denote $\cV_{\lambda}[1/\ell]$ by $V_{\lambda}$. 

From now on, we will be interested in understanding the cohomology groups 
\begin{equation}\label{eq:integral cohomology groups}
H^i\left(Y_K, \cV_{\lambda}\right), 
\end{equation}
which can be equipped with an action of the abstract Hecke algebra $\mathbb{T}^S$. Perhaps the most conceptual way of constructing the Hecke action in this setting is by observing that the pullback of $\sigma^\circ_{\lambda}$ to $G(F)\backslash G(\A_F)/K'_{\infty}$ is actually a $G(\A_{F}^S)\times K_S$-equivariant sheaf. The shadow of the $G(\A_{F}^S)$-action on this equivariant sheaf at infinite level induces an action of the Hecke algebra $\mathbb{T}^S$ on the finite level complex $R\Gamma(Y_K, \cV_{\lambda})$, viewed as an object in an appropriate derived category. Taking cohomology, one obtains an action of $\mathbb{T}^S$ on the $H^i\left(Y_K, \cV_{\lambda}\right)$. See~\cite[\S2]{newton-thorne} for more details.  

The cohomology groups in~\eqref{eq:integral cohomology groups} are finitely generated $\cO$-modules: this follows from the existence of the Borel--Serre compactification of $Y_K$, which is a compact manifold with corners that is homotopy equivalent to $Y_K$. However, these cohomology groups can and do contain torsion. For example, in the Bianchi case the torsion is conjectured to grow exponentially with the index of $K$ in $\GL_2(\widehat{\cO}_F)$ - see~\cite{BergeroVenkatesh} for more details. In what follows, the case of torsion classes will be the most subtle. We define 
\begin{equation}\label{eq:faithful Hecke algebra}
\mathbb{T}^S(K,\lambda):=\mathrm{Im}\left(\mathbb{T}^S\to \mathrm{End}_{\cO}\left(\bigoplus_{i=0}^{\dim_{\mathbb{R}} Y_K} H^i(Y_K, \cV_{\lambda})\right)\right)
\end{equation}
and formulate several conjectures about $\mathbb{T}^S(K,\lambda)$ and about the cohomology groups it acts on. 

To make precise statements, we introduce explicit generators for the spherical Hecke algebra at a place $v\notin S$ of $F$. We assume that, for such a place $v$, we have $K_v=\GL_n(\cO_{F_v})$. We let $\varpi_v$ denote a uniformizer of $F_v$ and $q_v$ denote the cardinality of the residue field $k_v:=\cO_{F,v}/\varpi_v$. For each $i=1,\dots, n$, we let $T_{v,i}$ denote the double coset operator in the spherical Hecke algebra of $\GL_n(F_v)$ given by
\[
T_{v,i} := \GL_n(\cO_{F_v})\mathrm{diag}[\varpi_v,\dots,\varpi_v, 1,\dots 1]\GL_n(\cO_{F_v}),  
\]
where $\varpi_v$ occurs $i$ times on the diagonal. 

\begin{conjecture}\label{conj:existence of Galois representations 1} 
Let $\mathfrak{m}\subset \mathbb{T}^S(K,\lambda)$ be a maximal ideal; its residue field is a finite extension $k'$ of the residue field $k$ of $E$. 
Then there exists a unique continuous, semi-simple Galois representation 
\[
\bar{\rho}_{\mathfrak{m}}: \Gamma_F\to \GL_n(k')
\]
characterized by the fact that, for any finite $v\notin S$, $\bar{\rho}_{\m}\mid_{\Gamma_{F_v}}$ is unramified and the characteristic polynomial of $\bar{\rho}_{\mathfrak{m}}(\Frob_v)$ is equal to the image of the polynomial 
\[
X^n - T_{v,1}X^{n-1} + \dots + (-1)^iq_v^{i(i-1)/2}T_{v,i}X^{n-i} + \dots + (-1)^n q_v^{n(n-1)/2} T_{v,n}. 
\]
modulo $\m$. 
\end{conjecture}

\begin{remark} This is a version of the second part of Conjecture~\ref{conj:GLC-ell-adic-G} in the case when $G=\GL_n$, formulated for more general coefficients. This can be seen from the Hochschild--Serre spectral sequence and from the fact that one can trivialize the local system $\cV_{\lambda}/\varpi$ by raising the level $\prod_{v\mid \ell}K_v$. 
Once again, the normalization is compatible with the $C$-algebraic version of the Buzzard--Gee conjecture.
\end{remark}

We say that a maximal ideal $\m\subset \mathbb{T}^S(K,\lambda)$ is \emph{non-Eisenstein} if the Galois representation $\bar{\rho}_{\m}$ exists and is absolutely irreducible. (This is a stronger condition than the notion of non-Eisenstein ideal going back to the work of Mazur~\cite{Mazur}. See~\cite[\S3.8]{CalegariVenkatesh} for a discussion of various related notions for $n=2$.) 

\begin{conjecture}\label{conj:existence of Galois representations 2}
Let $\m\subset \mathbb{T}^S(K,\lambda)$ be a non-Eisenstein maximal ideal. Then there exists a unique continuous Galois representation 
\[
\rho_{\m}: \Gamma_F\to \GL_n\left(\mathbb{T}^S(K,\lambda)^{\widehat{\ }}_{\m}\right)
\]
characterized by the fact that, for any finite place $v\notin S$ of $F$, $\rho_{\m}\mid_{\Gamma_{F_v}}$ is unramified and the characteristic polynomial of $\rho_{\m}(\Frob_v)$ is equal to the image of the polynomial 
\[
X^n - T_{v,1}X^{n-1} + \dots + (-1)^iq_v^{i(i-1)/2}T_{v,i}X^{n-i} + \dots + (-1)^n q_v^{n(n-1)/2} T_{v,n}. 
\]
in $\mathbb{T}^S(K,\lambda)^{\widehat{\ }}_{\m}[X]$. 
\end{conjecture}

We now discuss the status of Conjectures~\ref{conj:existence of Galois representations 1} and~\ref{conj:existence of Galois representations 2}. For a general number field, they are wide open. From now on, we restrict to the case when $F$ is a totally real or imaginary CM field. In this case, Conjectures~\ref{conj:existence of Galois representations 1} and~\ref{conj:existence of Galois representations 2} are essentially known and their proof, due to Scholze, has been discussed in \S~\ref{ss:construction-torsion-Galois}. The Galois representation $\rho_{\m}$ $\ell$-adically interpolates both the Galois representations attached to torsion classes occurring in $H^*(Y_K,\cV_{\lambda})_{\m}$ as well as the characteristic $0$ automorphic Galois representations first constructed in~\cite{HLTT}. 

We emphasise that Conjecture~\ref{conj:existence of Galois representations 2} is only proved up to a nilpotent ideal. More precisely, by work of Newton--Thorne~\cite{newton-thorne}, there exists a nilpotent ideal $J\subset \mathbb{T}^S(K,\lambda)^{\widehat{\ }}_{\m}$ with $J^4=1$ such that there exists a Galois representation 
\[
\rho_{\m}:\Gamma_F\to \GL_n\left(\mathbb{T}^S(K,\lambda)^{\widehat{\ }}_{\m}/J\right)
\]
with the right characteristic polynomials at places $v\not\in S$ of $F$. This does not usually cause trouble for applications to modularity lifting theorems, since the key point there is to understand the support of a certain patched module inside a Galois deformation ring. Nilpotent ideals do not affect the support. However, nilpotent ideals can cause trouble for more subtle questions, such as those concerning the Bloch--Kato conjecture, as in the appendix to~\cite{CGH}. 
If the prime $\ell$ splits completely in the totally real or imaginary CM field $F$, then Conjecture~\ref{conj:existence of Galois representations 2} is known as stated by the refined construction in~\cite{arizona}.  

\medskip 

\subsubsection{Local-global compatibility} Continue to assume that $F$ is a totally real or imaginary CM field and assume further that $\m\subset \mathbb{T}^S(K,\lambda)$ is a non-Eisenstein maximal ideal. Let $R_{\bar{\rho}_{\m}}$ denote the global deformation ring of $\bar{\rho}_{\m}$, which parametrizes deformations of $\bar{\rho}_{\m}$ that are unramified outside $S^\infty$. The Galois representation $\rho_{\m}$ is itself such a deformation, so by the universal property of $R_{\bar{\rho}_{\m}}$
we obtain a map 
\begin{equation}\label{eq:Galois to Hecke}
R_{\bar{\rho}_{\m}}\to \mathbb{T}^S(K,\lambda)^{\widehat{\ }}_{\m}/J.
\end{equation}
This is an incarnation of the usual map from a Galois deformation ring to a Hecke ring that shows up in both the classical Taylor--Wiles method and in the extension due to Calegari--Geraghty. The ultimate goal is to prove that this map is an isomorphism, or at least that it is an isomorphism on the underlying reduced rings. However, in order for the map~\eqref{eq:Galois to Hecke} to have a chance of being an isomorphism, we first need to refine it by imposing appropriate local conditions at the primes in $S$. The next conjecture concerns local-global compatibility for the Galois representation $\rho_{\m}$, which will help us refine this map by showing that it factors through a quotient of $R_{\bar{\rho}_{\m}}$ with the appropriate local conditions. 

Let $p$ be a prime and let $v\mid p$ be a place of $F$, possibly contained in $S^\infty$. If $\ell=p$, we set $\lambda_v= (\lambda_{\tau})_{\tau}$, where $\tau\in \Hom_{\Q}(F, E)$ runs over embeddings that induce the place $v$.  Local-global compatibility can roughly be formulated as the following question: how does the restriction $\rho_{\m}|_{\Gamma_{F_v}}$ depend on the level $K_v$ at which $\m$ occurs, on the Hecke eigenvalues at $v$, and, if $\ell =p$, on $\lambda_v$? 
\footnote{If $F$ is a totally real field, the description of $\rho_{\m}(c)$ for a choice of complex conjugation $c\in \Gal(\overline{F}/F)$ can also be interpreted as a form of local-global compatibility -- at the infinite places.}

For the characteristic $0$ automorphic Galois representations constructed in~\cite{HLTT}, this can be formulated more precisely by requiring compatibility with a suitable normalization of the classical local Langlands correspondence for $\GL_n/F_v$. In the case $\ell \not = p$, this compatibility up to semi-simplification is a theorem of Varma~\cite{VarmaThesis}. In the case $\ell = p$, one needs to first show that the automorphic Galois representations in question are de Rham at $v$ (or, equivalently, potentially semi-stable at $v$), then identify their Hodge--Tate numbers in terms of the sets $\lambda_v$, and finally establish the compatibility with classical local Langlands. 
In~\cite{A'Campo}, A'Campo recently proved that the automorphic Galois representations in question are indeed de Rham at primes $v\mid \ell$ under certain technical assumptions and identified their Hodge--Tate weights. 

The formulation of local-global compatibility for the integral Galois representations $\rho_{\m}$ is more subtle, particularly when $\ell=p$. This is because, outside a restricted family of cases such as the Fontaine--Laffaille and ordinary ones, it is not (a priori) clear how to formulate integral $p$-adic Hodge-theoretic conditions that should be satisfied by $\rho_{\m}|_{\Gamma_{F_v}}$. However, it turns out that we can formulate the local-global compatibility conjecture at $v\mid \ell$ using the potentially semi-stable Galois deformation rings constructed by Kisin in~\cite{kisin-PST}. This version of the conjecture first appears in the case $n=2$ in~\cite{gee-newton} and it is particularly well-suited for applications to modularity lifting theorems. 

For simplicity, we restrict to the crystalline case. Let $\bar{\rho}_{\m, v}$ denote the local representation $\bar{\rho}_{\m}|_{\Gamma_{F_v}}$. 
We let $R^{\square}_{\bar{\rho}_{\m, v}}$ denote its framed unrestricted local deformation ring. We let $R^{\square, \mathrm{crys}}_{\bar{\rho}_{\m, v}}(\lambda_v)$ denote its quotient whose characteristic $0$ points parametrize crystalline lifts of $\bar{\rho}_{\m,v}$ with Hodge--Tate weights equal to $\{\lambda_{\tau, n}, \lambda_{\tau, n-1}+1,\dots, \lambda_{\tau, 1}+ n-1\}$ at each embedding $\tau: F\hookrightarrow E$ inducing the place $v$. This quotient was constructed by Kisin in~\cite{kisin-PST}, first after inverting $\ell$ and then integrally by taking Zariski closure. 

\begin{conjecture}\label{conj:integral LGC} Assume that $K_v = \GL_n(\cO_{F_v})$, that $\m\subset \mathbb{T}^S(K,\lambda)$ is a non-Eisenstein maximal ideal, and that Conjecture~\ref{conj:existence of Galois representations 2} holds for $\m$. Consider the composition 
\begin{equation}\label{eq:galois to hecke factoring}
\xymatrix{R^{\square}_{\bar{\rho}_{\m,v}}\ar[r]\ar@{->>}[dr] & R^{\square}_{\bar{\rho}_{\m}}\ar[r] & \mathbb{T}^S(K,\lambda)^{\widehat{\ }}_{\m}
\\ \  & R^{\square, \mathrm{crys}}_{\bar{\rho}_{\m,v}}(\lambda_v)\ar@{-->}[ur]  &\ }
\end{equation}
where the first horizontal map is induced by restriction and the second horizontal map is induced from the existence of the Galois representation $\rho_{\m}$. This composition factors through $R^{\square, \mathrm{crys}}_{\bar{\rho}_{\m,v}}(\lambda_v)$. 
\end{conjecture}

\begin{remark}\label{rem:alternate integral LGC} We give an equivalent formulation of the local-global compatibility conjecture at a place $v\mid \ell$ when we are in the Fontaine--Laffaille setting. Concretely, we assume that $K_v = \GL_n(\cO_{F_v})$ as above, but also that $\ell$ is unramified in $F$ and that the vectors $\lambda_v$ are all in the Fontaine--Laffaille range, in other words
\[
0\leq \lambda_{\tau, n}\leq \dots \leq \lambda_{\tau, 1} \leq \ell-n-1
\]
for every $\tau: F\hookrightarrow E$ that induces $v$. Then it is equivalent to conjecture that $\rho_{\m}|_{\Gamma_{F_v}}$ is in the image of the Fontaine--Laffaille functor and has Hodge--Tate weights equal to $\{\lambda_{\tau, n}, \lambda_{\tau, n-1}+1,\dots, \lambda_{\tau, 1}+ n-1\}_{\tau}$ with $\tau: F\hookrightarrow E$ running over embeddings that induce $v$. One can state a similar explicit local-global compatibility conjecture in the ordinary case.  
\end{remark}

\begin{remark}\label{rem:potentially ss LGC}
It should be possible to formulate a version of Conjecture~\ref{conj:integral LGC} in the potentially crystalline case. For this, one could replace the algebraic representation $\sigma_{\lambda}$ with a locally algebraic representation $\sigma_{\tau, \lambda}:= \sigma_{\tau}\otimes \sigma_{\lambda}$ of $\GL_n(\cO_{F_v})$. Here, $\sigma_{\tau}$ should be taken to be a smooth representation of $\GL_n(\cO_{F_v})$ over a finite-dimensional $E'$-vector space, which should correspond to an inertial type $\tau: I_{F_v}\to \GL_n(E')$ under the inertial local Langlands correspondence, for some finite extension $E'/E$. For more details, see~\cite[\S\S3--4]{CEGGPS}, which describes the inertial local Langlands correspondence for $\GL_n$ in the potentially crystalline case. One can further generalize to the potentially semi-stable case by handling the monodromy operator carefully. 
\end{remark}

To make progress on Conjecture~\ref{conj:integral LGC}, we need a new way to access the Galois representations $\rho_{\m}$. The reason is that the construction of $\rho_{\m}$, which is described in \S~\ref{ss:construction-torsion-Galois}, involves increasing the level $\prod_{v\mid \ell} K_v$ as well as losing track of the weight $\lambda$ by trivializing the local systems $\cV_{\lambda}/\varpi^m$. 

At the moment, all progress on Conjecture~\ref{conj:integral LGC} has relied crucially on a series of increasingly more powerful vanishing theorems for the cohomology of unitary Shimura varieties with integral coefficients~\cite{CaraianiScholzeCompact, CaraianiScholzeNonCompact, koshikawa}. We discuss the general strategy in detail in \S~\ref{s:vanishing}. Based on this strategy, one can obtain partial results towards Conjecture~\ref{conj:integral LGC} in the Fontaine--Laffaille case, in the formulation of Remark~\ref{rem:alternate integral LGC}. This is done in~\cite[\S4]{10authors}. One can also obtain partial results in the ordinary case - see~\cite[\S5]{10authors}. 

The recent paper~\cite{CaraianiNewton} proves Conjecture~\ref{conj:integral LGC} in the case when $F$ is an imaginary CM field under some additional technical assumptions. Below, we give a flavour of the main assumptions; see Theorem 4.2.15 of \emph{loc.~cit.} for the precise statement.

\begin{enumerate}
    \item There is an assumption on the image of the residual representation $\bar{\rho}_{\m}$, which is needed in order to appeal to known vanishing results for the cohomology of unitary Shimura varieties with torsion coefficients. With a better understanding of these cohomology groups of Shimura varieties, this assumption can be progressively weakened.
    
    \item There is also an assumption that there are ``enough'' places above $\ell$ in the totally real subfield $F^+$ of $F$. This is needed in order to play these places against each other: one focuses on one place $\bar{v}\mid \ell$ of $F^+$ where one wishes to prove local-global compatibility, using the remaining places $\bar{v}'\mid \ell$, $\bar{v}'\not = \bar{v}$ as auxiliary ones. At auxiliary places, one employs congruences to vary the weight. This assumption excludes, for example, the case when $F$ is an imaginary quadratic field and the $Y_K$ are Bianchi manifolds.  
\end{enumerate}

\noindent It should be possible to adapt the method of~\cite{CaraianiNewton} to the potentially semi-stable case, under similar technical assumptions as described above. This is the subject of the upcoming PhD thesis of Bence Hevesi. 
 
\medskip

\subsubsection{Vanishing conjectures for locally symmetric spaces}

Conjectures~\ref{conj:existence of Galois representations 2} and~\ref{conj:integral LGC} are natural generalizations of results that hold true in the classical Taylor--Wiles setting. For example, they hold true in the conjugate self-dual case, when one can implement the classical Taylor--Wiles method starting from algebraic automorphic forms on a definite unitary group, as is done for example in~\cite{CHT}. 

We now discuss a third conjecture, which roughly says that the non-degenerate part of the cohomology $H^*(Y_K, \F_{\ell})$ is concentrated in a restricted range of degrees. This is the most novel prerequisite to the Calegari--Geraghty method and it illustrates some of the differences between general locally symmetric spaces and Shimura varieties. 

To state the conjecture, we define the following numerical invariants:
\[ 
l_0:= \mathrm{rk}\ \GL_n(F\otimes_{\Q}\R) - \mathrm{rk}\ K'_\infty\ \mathrm{and} 
\ q_0:=\frac{1}{2}(\dim_{\R} Y_K - l_0),
\]
which turn out to be non-negative integers. Concretely, if the number field $F$ has $r_1$ real places and $r_2$ complex places, one can compute that 
\[
l_0 = \begin{cases} r_1\left(\frac{n}{2}\right) + r_2n - 1\ \mathrm{for}\ n\ \mathrm{even} \\
r_1\left(\frac{n+1}{2}\right) + r_2n-1\ \mathrm{for}\ n\ \mathrm{odd}.\end{cases} 
\]
These invariants were first introduced by Borel--Wallach in~\cite{BorelWallach} for a general connected reductive group $\mathrm{G}$ over $\Q$. There, they show up naturally in the computation of $(\mathfrak{g}, K'_\infty)$-cohomology of tempered representations of $\mathrm{G}(\R)$.  

We consider the range of cohomological degrees $[q_0, q_0+l_0]$, which is symmetric about the middle $\frac{1}{2}\dim_{\R}Y_K$ of the total range of cohomology. The following conjecture is formulated in~\cite{Emerton-ICM} -- see the discussion around Conjecture 3.3 in \emph{loc.~cit.} -- and in~\cite[Conj.~B]{CG18}.

\begin{conjecture}\label{conj:vanishing of cohomology GL_n integral}
    Let $\m \subset \mathbb{T}^S$ be a non-Eisenstein maximal ideal in the support of $H^*(Y_K, \F_{\ell})$. Then $H^i(Y_K, \F_{\ell})_{\m}\not = 0$ only if $i\in [q_0, q_0+l_0]$. 
\end{conjecture}

\noindent This conjecture is motivated by heuristics to do with Langlands reciprocity, since the invariant $l_0$ can be recovered naturally from a computation on the Galois side. We explain this more in \S~\ref{ss:CG sketch}. We also discuss the analogue of Conjecture~\ref{conj:vanishing of cohomology GL_n integral} for Shimura varieties in \S~\ref{s:vanishing}. 

\begin{example}\leavevmode
    \begin{enumerate}
        \item If $n=2$ and $F=\Q$, we have $l_0 = 0$ and $q_0 =1$. There is only one interesting degree for the Betti cohomology of modular curves, which is $H^1$. This is the setting in which the Taylor--Wiles method was first implemented.   Conjecture~\ref{conj:vanishing of cohomology GL_n integral} can be verified by hand in the case of modular curves. In fact, the maximal ideals $\m\subset \mathbb{T}^S$ in the support of $H^0(Y_K, \F_{\ell})$ can be shown to satisfy
        \[\bar{\rho}_{\m} \simeq \chi \oplus \chi_{\mathrm{cyclo}}\cdot \chi,\]
        where $\chi: \Gamma_{\Q}\to \overline{\F}_{\ell}^\times$ 
        is some continuous mod $\ell$ character and $\chi_{\mathrm{cyclo}}$ is the mod $\ell$ cyclotomic character. 
        
        The Taylor--Wiles method was later generalized to other settings in which $l_0 = 0$. For example, in the conjugate self-dual case, one can arrange to work with algebraic automorphic forms on a definite unitary group, so that $q_0 = l_0 = 0$. 
        
        \item If $n=2$ and $F$ is imaginary quadratic, we have $l_0 = 1$ and $q_0=1$. The interesting degrees of Betti cohomology for Bianchi manifolds are $H^1$ and $H^2$. One can again prove Conjecture~\ref{conj:vanishing of cohomology GL_n integral} by hand in this case. By Poincar\'e duality, which applies after localization at a non-Eisenstein ideal $\m$, it is enough to control $H^0$, which can be computed explicitly.  
        
        \item It should be possible to control another handful of low-dimensional cases by employing the congruence subgroup property to control $H^1$. If $n=2$ and $F$ is totally real, the locally symmetric spaces are closely related to Hilbert modular varieties and one should be able to deduce (most of) Conjecture~\ref{conj:vanishing of cohomology GL_n integral} from~\cite[Thm.~A]{CaraianiTamiozzo}. For arbitrary $n$ and an arbitrary number field $F$, Conjecture~\ref{conj:vanishing of cohomology GL_n integral} is still open, even in the case when $F$ is an imaginary CM field.
            \end{enumerate}
\end{example}

One can formulate an analogous conjecture with $\Q_{\ell}$ rather than $\F_{\ell}$ coefficients. This is a strictly weaker conjecture, but it can be proved for arbitrary $n$ when $F$ is a CM or totally real field and this turns out to be extremely important for applications. 

\begin{conjecture}\label{conj:vanishing of cohomology GL_n rational}
 Let $\m \subset \mathbb{T}^S$ be a non-Eisenstein maximal ideal in the support of $H^*(Y_K, \Z_{\ell})$. Then $H^i(Y_K, \Q_{\ell})_{\m}\not = 0$ only if $i\in [q_0, q_0+l_0]$. 
\end{conjecture}

\noindent To prove Conjecture~\ref{conj:vanishing of cohomology GL_n rational}, one needs to compute the Betti cohomology groups $H^*(Y_K, \Q_{\ell})$ in terms of automorphic representations of $\GL_n(\A_{F})$.  This goes back to a theorem of Franke~\cite{Franke} that uses harmonic analysis. If $\m$ is non-Eisenstein and such that $H^*(Y_K, \Z_{\ell})_{\m}[1/\ell]\not = 0$, then it is the reduction modulo $\ell$ of a system of Hecke eigenvalues occurring in a cuspidal automorphic representation of $\GL_n(\A_{F})$. In that case, one can eliminate all the terms in Franke's computation coming from proper parabolic subgroups of $\GL_n/F$, and one controls what is left by appealing to the computation of tempered cohomology in~\cite{BorelWallach}. See~\cite[Thm.~2.4.10]{10authors} for the detailed argument. We emphasise, 
however, that one needs Conjecture~\ref{conj:existence of Galois representations 1} in order to formulate the notion of non-Eisenstein maximal ideal, so one needs to work over a CM or totally real field to prove Conjecture~\ref{conj:vanishing of cohomology GL_n rational}.  

\begin{remark} When $l_0 = 0$ and $\m\subset \mathbb{T}^S$ is a maximal ideal in the support of $H^*(Y_K, \Z_{\ell})$, Conjecture~\ref{conj:vanishing of cohomology GL_n integral} implies that $H^{*}(Y_K, \Z_{\ell})_{\m}$ is concentrated in degree $q_0$ and $\ell$-torsion free. This can be shown by playing around with the long exact cohomology sequences associated to the short exact sequences of sheaves 
\[
0\to \Z/\ell^{m-1}\Z \stackrel{\cdot \ell}{\longrightarrow} \Z/\ell^m\Z \to \F_{\ell}\to 0
\]
on $Y_K$ for $m \in \Z_{\geq 2}$. 
\end{remark}

\begin{remark}
All the conjectures discussed in this section can be formulated for the faithful quotient of $\mathbb{T}^S$ that acts on the cohomology of $Y_K$ with compact support. A useful observation is that, after localizing at a non-Eisenstein maximal ideal, the natural map from cohomology with compact support to usual cohomology becomes an isomorphism. This follows from a now standard argument that involves a detailed study of strata in the boundary of the Borel--Serre compactification of $Y_K$. See~\cite[\S4]{newton-thorne} for more details. 
Moreover, one could formulate Conjecture~\ref{conj:vanishing of cohomology GL_n integral} instead by asking that $H^*(Y_K, \F_{\ell})_{\m}$ is only non-zero in some range of degrees of length $l_0$. By the observation just above, Poincar\'e duality shows that this apparently weaker conjecture is equivalent to the original one. 
\end{remark}

\subsection{A sketch of the Calegari--Geraghty method}\label{ss:CG sketch}

In this subsection, we briefly sketch the Taylor--Wiles patching method for proving modularity lifting theorems and then we discuss its  improvement due to Calegari--Geraghty. The goal is to motivate and illustrate the role of Conjectures~\ref{conj:existence of Galois representations 2} and \ref{conj:integral LGC} and especially of Conjecture \ref{conj:vanishing of cohomology GL_n integral} in proving modularity lifting theorems. We suppress many technical details in our account. 
For a more thorough expository account of the classical Taylor--Wiles method, see~\cite{GeeArizona}, and for the same on the Calegari--Geraghty enhancement of this method, see~\cite{ThorneNotes}.  

\medskip 

\subsubsection{The case $l_0 =0$} We let $G=\GL_2/\Q$, in which case the $Y_K$ are (open) modular curves, considered as Riemann surfaces. We let $\m \subset \mathbb{T}^S$ be a non-Eisenstein maximal ideal in the support of $H^1(Y_K, \cV_{\lambda})$. This cohomology group is torsion-free, as discussed above, and so $\bar{\rho}_{\m}$ is (up to a character twist) the reduction modulo $\ell$ of the Galois representation associated to some cuspidal eigenform (of weight determined by $\lambda$). 

In this setting, Conjectures~\ref{conj:existence of Galois representations 2}, \ref{conj:integral LGC} and \ref{conj:vanishing of cohomology GL_n integral}  are known. For simplicity, we assume that $\ell>2$, that the level $K_\ell$ is the hyperspecial maximal compact subgroup $\GL_2(\Z_{\ell})$ and that the weight $\lambda$ is trivial. This means that $\bar{\rho}_{\m}$ is the reduction modulo $\ell$ of the Galois representation attached to a cuspidal  eigenform $f$ of weight $2$ and level $K$ (prime to $\ell$). The lift 
\[
\rho_{\m}: \Gal(\overline{\Q}/\Q)\to \GL_2\left(\mathbb{T}(K,\lambda)^{\widehat{\ }}_{\m}\right). 
\]
can be constructed by interpolating the $2$-dimensional $\ell$-adic Galois representations attached to \emph{all} the cuspidal eigenforms of weight $2$ and level $K$ and which are congruent to $f$ modulo $\m$. Local-global compatibility for each of these Galois representations implies that $\rho_{\m}|_{\Gamma_{{\Q}_{\ell}}}$ is Barsotti--Tate, which after inverting $\ell$ means crystalline with Hodge--Tate weights equal to $\{0,1\}$.

We let $R^{\mathrm{BT}}_{\bar{\rho}_{\m}}$ denote the global deformation ring for $\bar{\rho}_{\m}$ which parametrizes deformations that are unramified outside a fixed finite set of primes containing $\ell$, and that are Barsotti--Tate at $\ell$. The corresponding local deformation problem at $\ell$ turns out to be smooth and this is necessary for the argument we present. (In practice, we might impose some additional local conditions at ramified primes other than $\ell$, but we suppress this in this sketch.) The existence and local-global compatibility for the automorphic Galois representation $\rho_{\m}$ give rise to a diagram
\begin{equation}\label{eq:Galois to automorphic}
R^{\mathrm{BT}}_{\bar{\rho}_{\m}}\to \mathbb{T}(K,\lambda)^{\widehat{\ }}_{\m} \circlearrowleft H^1(Y_K, \cV_{\lambda})_{\m}. 
\end{equation}
To prove a modularity lifting theorem, we wish to show that the cohomology group $H^1(Y_K, \cV_{\lambda})_{\m}$ has full support over the global deformation ring $R^{\mathrm{BT}}_{\bar{\rho}_{\m}}$. (Since $\mathbb{T}(K,\lambda)^{\widehat{\ }}_{\m} $ acts faithfully on $H^1(Y_K, \cV_{\lambda})_{\m}$, this is equivalent to the statement that the surjective map $R^{\mathrm{BT}}_{\bar{\rho}_{\m}}\to \mathbb{T}(K,\lambda)^{\widehat{\ }}_{\m}$ induces an isomorphism on the underlying reduced rings.)

Rather than study the diagram~\eqref{eq:Galois to automorphic} directly, we first put everything in $\ell$-adic families, in order to smooth out the underlying geometry. For this, we wish to ``thicken'' the global deformation ring $R^{\mathrm{BT}}_{\bar{\rho}_{\m}}$ into
a ring $R_\infty$ which we will then prove to be isomorphic to
a power series ring $\cO[\![ x_1,\dots,x_g]\!]$. The dimension of $R_\infty$ is determined via a Galois cohomology computation that identifies $g$ with the dimension over $\mathbb{F}$ of a certain Selmer group $H^1_{\mathcal{L}}\left(\Gamma_{\Q},\mathrm{ad}\ \bar{\rho}_{\m}\right)$.  

In order to thicken the diagram~\eqref{eq:Galois to automorphic} in a controlled fashion, we ``patch'' (co)homology groups\footnote{Technically, patching amounts to taking a projective limit as the level varies. For this reason, it is better to switch from cohomology to homology at this point.} $H^1(Y_{K_{Q_N}}, \cV_{\lambda}/\ell^N)_{\m_{N}}$ with additional tame level at sets $Q_N$ of so-called Taylor--Wiles primes. Patching is a highly non-canonical process that involves varying the sets $Q_N$ of Taylor--Wiles primes as $N$ goes to infinity. However, the cardinality of these sets is a fixed integer $r$, which can be computed as the dimension over $\mathbb{F}$ of a certain dual Selmer group $H^1_{\mathcal{L}^\perp}\left(\Gamma_{\Q}, \mathrm{ad}\ \bar{\rho}_{\m}(1)\right)$. This fixed cardinality ensures that the tangent spaces of the deformation rings with additional tame ramification at the Taylor--Wiles primes stay bounded.  

The output of Taylor--Wiles patching is a diagram 
\begin{equation}\label{eq:patched Galois to automorphic}
S_{\infty} \to R_{\infty} \circlearrowleft M_{\infty}
\end{equation}
where $S_\infty = \cO[\![ z_1,\dots,z_{2r}]\!]$ and $R_\infty$ is a priori a complete local $\cO$-algebra which only admits a surjection from $\cO[\![ x_1,\dots,x_g]\!]$. The module $M_\infty$ is patched from (co)homology groups of the modular curve with additional tame level. The ring $R_\infty$ is patched from global deformation rings with additional tame ramification at Taylor--Wiles primes. The ring $S_\infty$ keeps track of the additional tame ramification coming from Taylor--Wiles primes. 

Because the growth of the (co)homology groups as we add tame level is strictly controlled, $M_\infty$ will be finite free over $S_\infty$; it is also finite over $R_\infty$. We have a sequence of (in)equalities
\begin{equation}\label{eq:inequalities 1}
2r+1 = \dim S_\infty = \mathrm{depth}_{S_\infty}(M_\infty) \leq \mathrm{depth}_{R_\infty}(M_\infty) \leq \dim R_{\infty} \leq g+1. 
\end{equation}
The equality $\dim S_\infty = \mathrm{depth}_{S_\infty}(M_\infty)$ follows from the freeness of $M_\infty$ over $S_\infty$. 
The first inequality comes from the fact that the action of $S_\infty$ on $M_\infty$ factors through $R_\infty$. The second inequality comes from the fact that the depth of a finitely generated module is less than or equal to the dimension of its support. 

At this point,  Wiles makes use of an Euler characteristic formula in Galois cohomology that gives $2r-g = l_0$, which is $0$ in this setting. This implies that all the above inequalities are in fact equalities! In particular, this means that the given surjection $\cO[\![ x_1,\dots,x_g]\!]\twoheadrightarrow R_\infty$ is in fact an isomorphism.  We now apply the Auslander--Buchsbaum formula for the finitely generated module $M_\infty$ over the regular local ring $R_\infty$:
\begin{equation}\label{eq:Auslander-Buchsbaum}
\mathrm{depth}_{R_\infty}(M_\infty) + \mathrm{proj.}\ \mathrm{dim.}_{R_\infty}(M_\infty) = \dim R_\infty = g+1. 
\end{equation}
Since we have also found that $\mathrm{depth}_{R_\infty}(M_\infty) = g+1$, we deduce $M_\infty$ is finite free over $R_\infty$. 

To recover the original diagram~\eqref{eq:Galois to automorphic}, we take a tensor product of~\eqref{eq:patched Galois to automorphic} with  $\otimes_{S_\infty}\cO$.
We have $R_\infty \otimes_{S_\infty} \cO = R^{\mathrm{BT}}_{\bar{\rho}_{\m}}$ and $M_\infty\otimes_{S_\infty}\cO$ is (the $\cO$-linear dual of) $H^1(Y_K, \cV_{\lambda})_{\m}$. This shows that the cohomology group $H^1(Y_K, \cV_{\lambda})_{\m}$ is finite free over the global deformation ring $R^{\mathrm{BT}}_{\bar{\rho}_{\m}}$ and, in particular, has full support. 

\begin{remark} In the above sketch, we can weaken the assumptions that the level $K_\ell$ is hyperspecial and that the weight $\lambda$ is trivial. The right perspective here is due to Kisin~\cite{kisin-PST} and involves constructing $R_\infty$ as a power series ring over a product of local (possibly framed) deformation rings. We can then show that $M_\infty$ is maximal Cohen--Macaulay over $R_\infty$ and thus supported on a union of irreducible components of $\Spec R_\infty$. Therefore, the geometry of the local deformation rings determines the kind of modularity lifting result we can prove. A key problem is how to extend the known support of $M_\infty$ from one irreducible component to another, e.g. as in Taylor's Ihara avoidance. See~\cite{GeeArizona} for more details. The same phenomenon arises in the case when $l_0>0$ and, in this case, it is even more subtle to move between different irreducible components. 
\end{remark}

\medskip 

\subsubsection{The case $l_0>0$} For simplicity and to parallel the case $l_0=0$ discussed above, we restrict to $G=\GL_2/F$, where $F$ is now an imaginary CM field. For a neat compact open subgroup $K\subset \GL_2(\A_{F}^\infty)$, the spaces $Y_K$ are smooth real manifolds. We assume that the prime $\ell>2$ is unramified in $F$, that the level $K_v$ at each $v\mid \ell$ is hyperspecial, and that the weight $\lambda$ is trivial. We let $\m \subset \mathbb{T}^S$ be a non-Eisenstein maximal ideal in the support of $H^*(Y_K, \cV_{\lambda})$.

We assume Conjectures~\ref{conj:existence of Galois representations 2},
\ref{conj:integral LGC} and~\ref{conj:vanishing of cohomology GL_n integral}. As explained above, Conjecture~\ref{conj:existence of Galois representations 2} is essentially known by~\cite{Sch15}, and there are many settings where Conjecture~\ref{conj:integral LGC} is known as well in the case at hand, see~\cite{CaraianiNewton} for the state of the art result. However, Conjecture~\ref{conj:vanishing of cohomology GL_n integral} is still a significant assumption. The main exception to this is the case when $F$ is imaginary quadratic, where Conjecture~\ref{conj:vanishing of cohomology GL_n integral} is known, but Conjecture~\ref{conj:integral LGC} is surprisingly subtle and not known even in the case at hand.  

We let $R^{\mathrm{BT}}_{\bar{\rho}_{\m}}$ denote the universal deformation ring of $\bar{\rho}_{\m}$ with appropriate local conditions: in particular, we require these deformations to be unramified outside a fixed finite set of finite primes of $F$ and Barsotti--Tate at each prime $v\mid \ell$ of $F$. As before, we obtain a diagram
\begin{equation}\label{eq:Galois to automorphic CG}
R^{\mathrm{BT}}_{\bar{\rho}_{\m}}\to \mathbb{T}(K,\lambda)^{\widehat{\ }}_{\m} \circlearrowleft \bigoplus_{i=q_0}^{q_0+l_0} H^i(Y_K, \cV_{\lambda})_{\m}. 
\end{equation}

To prove a modularity lifting theorem in this setting, we want to show that
the direct sum has full support in $R^{\mathrm{BT}}_{\bar{\rho}_{\m}}$. 
It is, of course, enough to establish this for any individual cohomological degree. 

We want to proceed as in the case $l_0 = 0$ discussed above. To thicken the diagram~\eqref{eq:Galois to automorphic CG}, we again need to consider a patched deformation ring $R_\infty$, with $\dim R_\infty = g+1$, where $g$ is the dimension over $\F$ of a Selmer group $H^1_{\mathcal{L}}(\Gamma_F, \mathrm{ad}\ \bar{\rho}_{\m})$. At the same time, the cardinality of the sets of Taylor--Wiles primes must equal $r$, the dimension over $\F$ of the dual Selmer group $H^1_{\mathcal{L}^\perp}(\Gamma_F, \mathrm{ad}\ \bar{\rho}_{\m}(1))$. We therefore set $S_\infty = \cO[\![z_1,\dots, z_{2r}]\!]$ once again. However, if we try to implement the classical Taylor--Wiles method as outlined above, we get stuck at the Euler characteristic formula: 
\begin{equation}\label{eq:Euler characteristic}
\dim S_\infty - \dim R_\infty = l_0 > 0. 
\end{equation}
This means that there are not enough Galois representations to fill out the entire space $\Spec S_\infty$. On the other hand, because the cohomology groups of $Y_K$ are likely supported in multiple degrees, if we add tame level at sets of Taylor--Wiles primes and patch (co)homology groups (in a fixed degree or even taking a direct sum of cohomological degrees), we will likely not obtain a finite free module $M_\infty$ over $S_\infty$. Thus, it seems like there are also not enough automorphic forms to fill out $\Spec S_\infty$. Fortunately, as Calegari and Geraghty have observed, these two problems precisely cancel each other out, when considered from the right point of view! 

The idea is to work with a patched complex $C_\infty$ instead of a single patched module $M_\infty$. The complex $C_\infty$ is obtained from patching chain complexes that compute the homology of $Y_{K_{Q_N}}$, with additional tame level at a cardinality $r$ set of Taylor--Wiles primes $Q_N$, and after localizing at a non-Eisenstein maximal ideal $\m_N$ obtained from $\m$. Because we are working with the homology of the $Y_{K_{Q_N}}$ as chain complexes, rather than with individual homology groups, we obtain after patching a complex $C_\infty$ of finite projective $S_\infty$-modules. Moreover, Conjecture~\ref{conj:vanishing of cohomology GL_n integral} ensures that, up to replacing $C_\infty$ by a quasi-isomorphic complex, we can assume that it is a complex of finite projective $S_\infty$-modules concentrated in a range of degrees of length $l_0$. The following lemma in commutative algebra is now crucial, see~\cite[Lem.~6.2]{CG18} for its proof. 

\begin{lemma}\label{lem:commutative algebra} Let $S$ be a regular Noetherian local ring of dimension $d$. Let $C$ be a chain complex\footnote{In~\cite[Lem.~6.2]{CG18} the result is stated and proved for a cochain complex.} of finite projective $S$-modules concentrated in the range of degrees $[0,l]$ for some non-negative integer $l\leq d$. Then the dimension of the support of $H^*(C)$ in $S$ satisfies $\dim_{S}(H_*(C))\leq d-l$.
Moreover, if equality holds, then $H_i(C)$ is non-zero only for $i=0$,
\[
\mathrm{proj.}\ \mathrm{dim.}_{S}(H_0(C)) = l,\ \mathrm{and}\  \mathrm{depth}_{S}(H_{0}(C))= d-l. 
\]
(In particular, $C$ is a projective resolution of $H_0(C)$ of minimal length.) 
\end{lemma}

\noindent In the Calegari--Geraghty extension of the Taylor--Wiles method, we apply this lemma to $S=S_\infty$ and $C=C_\infty[q_0]$ (we shift our patched complex $q_0$ degrees to the left in order to get the concentration in the range $[0,l_0]$). We have a diagram 
\[
S_\infty \to R_\infty \circlearrowleft C_\infty, 
\]
giving us the sequence of inequalities 
\begin{equation}\label{eq:inequalities 2}
2r+1-l_0 \leq \dim_{S_\infty}(H_*(C_\infty))\leq \dim_{R_\infty}(H_*(C_\infty))\leq g+1. 
\end{equation}
The Euler characteristic formula~\eqref{eq:Euler characteristic} implies that all the inequalities in~\eqref{eq:inequalities 2} are equalities. Lemma~\ref{lem:commutative algebra} implies that the homology of $C_\infty$ is non-zero only in degree $q_0$ and that 
\[
\mathrm{depth}_{S_\infty}(H_{q_0}(C_\infty)) = 2r+1-l_0 = g+1.
\] 
At the same time, because the action of $S_\infty$ on $C_\infty$ factors through $R_\infty$, we have the (in)equalities
\[
g+1 = \mathrm{depth}_{S_\infty}(H_{q_0}(C_\infty))\leq \mathrm{depth}_{R_\infty}(H_{q_0}(C_\infty)) \leq \dim R_\infty = g+1. 
\]
These are all equalities and the Auslander--Buchsbaum formula applied to $H_{q_0}(C_\infty)$ viewed as an $R_\infty$-module implies that $H_{q_0}(C_\infty)$ is free over $R_\infty$. To recover (part of) the original diagram~\eqref{eq:Galois to automorphic CG}, we take a tensor product with $\otimes_{S_\infty}\cO$\footnote{Here, one could ask what would happen if we took a derived tensor product $\otimes^{\mathrm{L}}_{S_\infty}\cO$ instead. This will recover the entire cohomology as an object in a derived category. See the article~\cite{FengHarrisIHES} of Feng--Harris in these proceedings for more details on derived aspects.} and an $\cO$-linear dual to go from homology to cohomology. We ultimately deduce that $H^{q_0}(Y_K, \cV_{\lambda})_{\m}$ has full support over $R^{\mathrm{BT}}_{\bar{\rho}_{\m}}$, as desired. 
\medskip

\section{Shimura varieties, $p$-adic geometry, and torsion classes}\label{s:vanishing}

In this section, we discuss the steps that led to an unconditional implementation of the Calegari--Geraghty method for $\GL_n$ over CM fields. A key input is an analogue of the vanishing Conjecture~\ref{conj:vanishing of cohomology GL_n integral} for Shimura varieties. This analogue predicts that the part of the cohomology of Shimura varieties outside the middle degree is somewhat degenerate. In \S~\ref{ss:geometry of pi_HT}, we discuss the geometry of the Hodge--Tate period morphism, which has allowed us to make significant progress towards this conjecture, and which is a beautiful and important topic in its own right. In \S~\ref{ss:vanishing}, we discuss several different approaches towards the conjecture in the case of Shimura varieties PEL type A and point towards generalizations. Finally, in \S~\ref{ss:boundary}, we briefly describe how to use these results to make progress on the local-global compatibility Conjecture~\ref{conj:integral LGC} and mention further applications. 

\subsection{The geometry of the Hodge--Tate period morphism}\label{ss:geometry of pi_HT}

The goal of this subsection is to discuss more in depth the geometry of the Hodge--Tate period morphism that has already played a key role in Section~\ref{ss:perfectoid-Shimura}. We will briefly explain how this geometry is illuminated by an infinite-level version of Mantovan's product formula~\cite{mantovan-duke}. This product formula describes Newton strata in Shimura varieties in terms of local analogues of Shimura varieties, known as Rapoport--Zink spaces or local Shimura varieties, and objects of a semi-global nature, known as Igusa varieties. We end by mentioning an even more general version of the product formula, which allows us to vary the Newton stratum, conjectured by Scholze and recently established in many cases by Mingjia Zhang in her Bonn PhD thesis~\cite{mingjia}.  

For simplicity, we work with a Shimura datum $(G,X)$ of PEL type. We assume that $p$ is an unramified prime for this Shimura datum and $\ell$ is a rational prime different from $p$. We fix a sufficiently small compact open subgroup $K^p\subset G(\A^{p,\infty})$. Consider the minimally compactified perfectoid Shimura variety $\cS^*_{K^p}$ over $\C_p$ introduced in Theorem~\ref{thm:perfectoid-Sh} and the associated Hodge--Tate period morphism
\[
\pi_{\mathrm{HT}}: \cS^*_{K^p}\to \Fl_G.
\]
For technical reasons, we restrict the source of this morphism from now on to the good reduction locus $\cS^\circ_{K^p}\subset \cS^{*}_{K^p}$. 

\medskip 

\subsubsection{The Newton stratification}
Recall the Kottwitz set $B(G_{\Q_p})$ classifying isocrystals with $G_{\Q_p}$-structure, cf.~\cite{rapoport-richartz} for example. This set is equipped with a partial order, known as the Bruhat order. The Hodge cocharacter $\mu$ defines a subset $B(G_{\Q_p}, \mu^{-1})\subset B(G_{\Q_p})$ of $\mu^{-1}$-admissible elements.  
The special fiber of the Shimura variety with hyperspecial level at $p$ admits a Newton stratification
\begin{equation}\label{eq:Newton stratification special}
\overline{\Sh}_{K^pK_p} = \bigsqcup_{b\in B(G_{\Q_p}, \mu^{-1})} \overline{\Sh}^b_{K^pK_p}
\end{equation}
into locally closed strata indexed by this subset. This stratification is defined moduli-theoretically and detects at geometric points the isogeny class of the universal $p$-divisible group $\mathscr{G} = \cA[p^\infty]$ with $G_{\Q_p}$-structure living over $\overline{\Sh}_{K^pK_p}$. The closure relations between strata are determined by the degenerations of isocrystals with $G_{\Q_p}$-structure: for example, the $\mu$-ordinary locus is the unique open and dense stratum and the basic locus is the unique Zariski closed stratum. The Newton stratification can be extended in a natural way to the minimal compactification $\overline{\Sh}^*_{K^pK_p}$ using the notion of a well-positioned subset of the special fiber, cf.~\cite{lan-stroh}. 

There is an analogue of the Newton stratification on $\mathscr{F}\ell_G$ that is defined using purely local ingredients. More precisely, we have a decomposition
\begin{equation}\label{eq:Newton stratification flag}
\Fl_G = \bigsqcup_{b\in B(G_{\Q_p}, \mu^{-1})} \Fl^b_G 
\end{equation}
into locally closed strata, cf.~\cite[\S3]{CaraianiScholzeCompact}. Morally, the argument is as follows (although this language did not exist at the time when~\cite{CaraianiScholzeCompact} was written): one first constructs a map of v-stacks $\mathscr{F}\ell_G\to \mathrm{Bun}_G$, where the latter is the v-stack
of $G$-bundles on the Fargues--Fontaine curve. To construct this map of v-stacks, it is convenient to 
notice that one can identify the diamond associated to $\mathscr{F}\ell_G$ with the minuscule Schubert cell 
defined by $\mu$ inside the $B^+_{\mathrm{dR}}$-Grassmannian for $G$. 
Once we have obtained the map $\Fl_G\to \mathrm{Bun}_G$, we can use Fargues's result that the points of $\mathrm{Bun}_G$ are in bijection with the 
Kottwitz set $B(G_{\Q_p})$, cf.~\cite{fargues-G-torsors} (see also~\cite{anschutz} for an alternative proof that also 
works in equal characteristic). 
Moreover, the Newton decomposition~\eqref{eq:Newton stratification flag} is a \emph{stratification},
in the sense that, for $b\in B(G_{\Q_p}, \mu^{-1})$, we have a set-theoretic decomposition
\[
\overline{\mathscr{F}\ell^b_G}= \bigsqcup_{b'\geq b} \mathscr{F}\ell^{b'}_G,  
\]
where $\geq$ denotes the Bruhat order.
The latter fact follows from the analogous result of Viehmann that applies to $\mathrm{Bun}_G$, see~\cite[Thm.~1.1]{viehmann}. 

\begin{remark}\label{rem:opposite closure relations} The closure relations between Newton strata on the special fiber $\overline{\Sh}^*_{K^pK_p}$ are \emph{opposite} to the closure relations between Newton strata on $\Fl_G$. This is not an inconsistency, because the Hodge--Tate period map respects the Newton stratification only on points of rank one. The behaviour of $\pi_{\mathrm{HT}}$ on higher rank points is more subtle and this accounts for the difference in closure relations. See also~\cite[Theorem 1.12]{gleason-ivanov} for a modern proof that the closure relations on the special and generic fibres are reversed in the context of the stack of isocrystals with $G$-structure $\mathfrak{B}(G)$ and the stack $\mathrm{Bun}_G$ respectively. 
\end{remark}

\subsubsection{Igusa varieties} For each $b\in B(G_{\Q_p}, \mu^{-1})$, we can choose a (completely slope divisible) 
$p$-divisible group with $G_{\Q_p}$-structure $\mathbb{X}_b/\overline{\F}_p$, which is in the isogeny class determined by $b$ and also compatible with $\mu$, and define the corresponding 
\emph{Oort central leaf}. This is a smooth closed subscheme $\mathscr{C}^{\mathbb{X}_b}$ 
of the Newton stratum $\overline{\Sh}^b_{K^pK_p}$, such that the isomorphism class 
of the universal $p$-divisible group $\mathscr{G}$ with $G_{\Q_p}$-structure over each geometric 
point of the leaf is constant and equal to that of $\mathbb{X}_b$: 
\[
\mathscr{C}^{\mathbb{X}_b}=\left\{x\in \overline{\Sh}^b_{K^pK_p} 
\mid \mathscr{G}\times \kappa(\bar{x})
\simeq \mathbb{X}_b\times \kappa(\bar{x})\right\}.
\]
In general, there can be infinitely many non-empty leaves inside a given Newton stratum. 
Over each central leaf, one has the \emph{perfect Igusa variety} $\mathrm{Ig}^b/\overline{\F}_p$, 
a pro-finite cover of $\mathscr{C}^{\mathbb{X}^b}$ 
which parametrizes trivializations of $\mathscr{G}$ together with its $G_{\Q_p}$-structure. 

Variants of Igusa varieties were introduced in~\cite{HarrisTaylor}
in the special case of Shimura varieties of Harris--Taylor type. They were 
defined more generally for Shimura varieties of PEL type by Mantovan~\cite{mantovan-thesis} and their $\ell$-adic 
cohomology was computed in many cases using a counting point formula; see~\cite[Ch.V]{HarrisTaylor} as well as~\cite{ShinIgusa,ShinStable,ShinGalois}. The overall strategy is a variation of the LKR method (\S~\ref{sss:LKR-method}). Steps (1)--(4) of \S~\ref{sss:LKR-method} have analogues for Igusa varieties. See \cite[\S7]{ShinICM} for more details and the generalization to the case of Hodge or abelian type.

In these references, Igusa varieties are defined as pro-finite \'etale covers of central leaves, which trivialize the graded pieces of the slope filtration, which $\mathscr{G}$ admits after restriction to $\mathscr{C}^{\mathbb{X}^b}$. Providing a splitting of the slope filtration amounts to taking perfection on the level of the Igusa variety, recovering the object $\mathrm{Ig}^b$ introduced above. From now on, we only consider the perfect Igusa varieties, since they have a more elegant moduli-theoretic interpretation, but still have the same $\ell$-adic cohomology as the original objects.

\begin{remark}
While the central leaf $\mathscr{C}^{\mathbb{X}_b}$ depends
on the choice of $\mathbb{X}_b$ in its isogeny class, one can show that, up to isomorphism, the perfect Igusa 
variety $\mathrm{Ig}^b$ only depends on the isogeny class: this follows from
the equivalent moduli-theoretic description in~\cite[Lem.~4.3.4]{CaraianiScholzeCompact} 
(see also~\cite[Lem.~4.2.2]{CaraianiTamiozzo}, which keeps track of the extra structures more carefully). 
In particular, the pair $(G,\mu)$ is not intrinsically attached to the Igusa variety $\mathrm{Ig}^b$ -- 
it can happen that Igusa varieties that are a priori obtained from different Shimura varieties
are isomorphic. See~\cite[Thm.~4.2.4]{CaraianiTamiozzo} for an example of this phenomenon.  
\end{remark}

Because $\mathrm{Ig}^b/\overline{\F}_p$ is perfect, the base change $\mathrm{Ig}^b\times_{\overline{\F}_p}\cO_{\C_p}/p$ 
admits a canonical lift to a flat formal scheme over $\mathrm{Spf}\ \cO_{\C_p}$. We let $\mathfrak{Ig}^b$ denote 
the adic generic fibre of this lift, which is a perfectoid space over $\mathrm{Spa}(\C_p,\cO_{\C_p})$. 
The spaces $\mathrm{Ig}^b$ and $\mathfrak{Ig}^b$ both have an action of a locally profinite group 
$G_b(\Q_p)$, where $G_b$ is an inner form of a Levi subgroup of $G$, and they have a Hecke action away from $p$. The $\ell$-adic cohomology groups of $\mathrm{Ig}^b$ and $\mathfrak{Ig}^b$ are isomorphic, equivariantly for these actions.     

\medskip 

\subsubsection{Rapoport--Zink spaces} For each $b\in B(G_{\Q_p}, \mu^{-1})$, one can also consider the associated 
\emph{Rapoport--Zink space}, a moduli space of $p$-divisible groups with $G_{\Q_p}$-structure
that is a local analogue of a Shimura variety. Concretely in the PEL case, 
one considers a moduli problem 
of $p$-divisible groups equipped with $G_{\Q_p}$-structure, satisfying the Kottwitz determinant 
condition with respect to $\mu$, and with a modulo $p$
quasi-isogeny to the fixed $p$-divisible group $\mathbb{X}_b$. This moduli problem
was shown by Rapoport--Zink~\cite{rapoport-zink} to be representable by a formal scheme. We let $\cM^b$ denote the adic generic fibre
of this formal scheme, base changed to $\mathrm{Spa}(\C_p,\cO_{\C_p})$, and let $\cM^b_{\infty}$ denote
the corresponding infinite-level Rapoport--Zink space\footnote{As a consequence of the comparison
with moduli spaces of local shtukas in~\cite{berkeley-lectures}, one obtains a group-theoretic characterisation of
Rapoport--Zink spaces as local Shimura varieties determined by the tuple $(G,b,\mu)$. We suppress $(G,\mu)$ from
the notation for simplicity.}. 
The latter object can be shown to be a perfectoid space using the techniques of~\cite{ScholzeWeinstein}. 
By \emph{loc.~cit.}, the infinite-level Rapoport--Zink space 
admits a local analogue of the Hodge--Tate period morphism 
\[
\pi^b_{\mathrm{HT}}: \cM^b_{\infty} \to \mathscr{F}\ell_G.  
\]
The geometry of $\pi_{\mathrm{HT}}$ is intricately 
tied up with the geometry of its local analogues $\pi^b_{\mathrm{HT}}$, as we shall now explain.  

\medskip 

\subsubsection{The product formula}\label{sss:product formula}

For each $b\in B(G,\mu^{-1})$, 
we define a locally closed subset $\cS^{\circ b}_{K^p}$ of the good reduction locus $\cS^\circ_{K^p}$ by intersecting the pre-image under the specialization map of the Newton stratum corresponding to $b$ on the special fiber with the pre-image under $\pi_{\mathrm{HT}}$ of $\Fl^b_{G}$.  

\begin{theorem}\label{thm:product formula} There exists a Cartesian 
diagram of diamonds over $\mathrm{Spd}(\C_p, \cO_{\C_p})$
\begin{equation}\label{eq:upper Cartesian}
\xymatrix{
\cM^b_{\infty}\times_{\mathrm{Spd}(\C_p, \cO_{\C_p})} \mathfrak{Ig}^b \ar[r]\ar[d] & 
\cM^b_{\infty}\ar[d]^{\pi^b_{\mathrm{HT}}} \\
 \cS^{\circ b}_{K^p}\ar[r]^{\pi_{\mathrm{HT}}} & \mathscr{F}\ell^b_G.}
\end{equation}
Moreover, each vertical map is a pro-\'etale torsor for the ``unipotent group diamond'' $\widetilde{G}_b$ of~\cite{fargues-scholze}
(identified with $\mathrm{Aut}_G(\widetilde{\mathbb{X}}_b)$, in the 
notation of~\cite[\S4]{CaraianiScholzeCompact}).  
\end{theorem}

Theorem~\ref{thm:product formula} is established in~\cite[\S4]{CaraianiScholzeCompact}, although it is not formulated in terms of diamonds in \emph{loc.~cit.}
It has been generalized to Shimura varieties of Hodge type, in the context of Kisin--Pappas integral models (and under some additional technical assumption), by Hamacher--Kim~\cite{HamacherKim}. 
At the time of writing, a general version of the Mantovan product formula for Shimura varieties of abelian type is not yet available. 

Assume that the Shimura varieties $\mathrm{Sh}_K$ are compact. 
We have the following consequence of Theorem~\ref{thm:product formula} for the fibres of $\pi_{\mathrm{HT}}$: let 
$\bar{x}: \mathrm{Spa} (C,C^+)\to \mathscr{F}\ell^b_G$
be a geometric point. Then there is an inclusion of $\mathfrak{Ig}^b$ into 
$\pi_{\mathrm{HT}}^{-1}(\bar{x})$, which identifies the target with the canonical compactification of the source,
in the sense of~\cite[Prop.~18.6]{scholze-diamonds}. 

This computation of the fibres can be extended to minimal and toroidal compactifications of (non-compact) Shimura varieties. For a model of the argument, see~\cite[Thm.~1.10]{CaraianiScholzeNonCompact}, where the computation of the fibres is extended to Shimura varieties 
attached to quasi-split unitary groups. In this case, the fibres are obtained from partial minimal and toroidal compactifications of Igusa varieties. The partial minimal compactifications of Igusa varieties are affine - this is reflected in the ``affinoid'' nature of $\pi_{\mathrm{HT}}$ for the minimal compactification $\cS^*_{K^p}$. This result has been generalized to Shimura varieties of PEL type AC in the PhD thesis of Mafalda Santos~\cite{MafaldaSantos}.  

\begin{example}\label{ex:geometry ht for modular curve} 
We make the geometry of $\pi_{\mathrm{HT}}$ explicit in the case of the modular curve, i.e. for $G=\GL_2/\Q$. In this case, we identify $\mathscr{F}\ell_G = \mathbb{P}^{1,\mathrm{ad}}$ and we have 
the decomposition into Newton strata (at least on points of rank one): 
\[
\xymatrix{
\cS^{\circ}_{K^p} \ar[d]^{\pi_{\mathrm{HT}}} & = & \cS^{\circ,\mathrm{ord}}_{K^p}\ar[d] & \sqcup & \cS^{\circ, \mathrm{ss}}_{K^p}\ar[d] \\
\mathbb{P}^{1,\mathrm{ad}} & = & \mathbb{P}^{1,\mathrm{ad}}(\mathbb{Q}_p) & \sqcup &\Omega. 
} 
\]
The ordinary locus inside $\mathbb{P}^{1,\mathrm{ad}}$
consists of the set of points defined over $\Q_p$ and the basic / supersingular 
locus is its complement $\Omega$, the Drinfeld upper half-plane. 

The fibres
of $\pi_{\mathrm{HT}}$ over the ordinary locus are ``perfectoid versions'' of Igusa curves. 
The infinite-level version of the product formula reduces, in this case, to 
the statement that the ordinary locus is parabolically induced from 
$\mathfrak{Ig}^{\mathrm{ord}}$, as in~\cite[\S6]{CaraianiTamiozzo}.   
The fibres of $\pi_{\mathrm{HT}}$ over the supersingular locus are profinite sets: 
the corresponding Igusa varieties can be identified with double cosets 
$D^\times \backslash D^\times(\A_f)/K^p$,
where $D/\Q$ is the quaternion algebra ramified precisely at $\infty$ and $p$. In some form, this result goes back to Deuring--Serre; see~\cite{Howe} for this precise formulation. This description of the supersingular / basic Igusa varieties together with the Mantovan product formula recover Rapoport--Zink uniformization in this special case. 
\end{example}

\begin{remark}\label{rem:Mingjia's PhD thesis} The Cartesian diagram~\eqref{eq:upper Cartesian} can be obtained from a Cartesian diagram of $v$-stacks 
\begin{equation}\label{eq:lower Cartesian}
\xymatrix{
\cS^{\circ b}_{K^p}\ar[r]^{\pi_{\mathrm{HT}}}\ar[d] & \Fl^b_G\ar[d] \\ 
[\mathfrak{Ig}^b/\widetilde{G}_b] \ar[r] & 
[*/\widetilde{G}_b]\ar[r]^{\simeq} & \mathrm{Bun}^b_G
}. 
\end{equation} 
Indeed, pullback of~\eqref{eq:lower Cartesian} under $*\to [*/\widetilde{G}_b]$ produces a cube with all faces Cartesian diagrams and one of the faces can be identified with~\eqref{eq:upper Cartesian}. 
In her PhD thesis~\cite{mingjia}, 
Zhang constructs a more general version of~\eqref{eq:lower Cartesian}, where the element $b$, and thus the Newton stratum, is allowed to vary. 

More precisely, for Shimura varieties of PEL type AC and when $p$ is a prime of good reduction, Zhang constructs a $v$-stack $\mathrm{Igs}^\circ_{K^p}$ living over $\mathrm{Bun}_G$, called an \emph{Igusa stack}. She then establishes a Cartesian diagram of $v$-stacks
\begin{equation}\label{eq:Mingjia}
\xymatrix{
\cS^\circ_{K^p}\ar[d]\ar[r]^{\pi_{\mathrm{HT}}} & \Fl_G\ar[d] \\ 
\mathrm{Igs}^\circ_{K^p} \ar[r]^{\overline{\pi}_{\mathrm{HT}}} & \mathrm{Bun}_G 
}
\end{equation}
which specializes to~\eqref{eq:lower Cartesian} for fixed $b\in B(G,\mu^{-1})$. Under some mild technical assumptions, Zhang also extends~\eqref{eq:Mingjia} to minimal compactifications. That such a construction should be possible for general Shimura varieties was conjectured by Scholze.

Morally, one should think of the diagram~\eqref{eq:Mingjia} as taking as an input a moduli space of global $G$-iso-shtukas with no legs, and creating legs via Beauville--Laszlo gluing. This intuition cannot be made precise, of course, since we do not have a good notion of global shtukas in the number field setting. However, this does suggest that Igusa varieties and, more generally, Igusa stacks are, in some ways, just as fundamental as Shimura varieties. 
\end{remark}

\subsection{Vanishing of cohomology with torsion coefficients}\label{ss:vanishing}

The goal of this subsection is to discuss local counterparts to Conjecture~\ref{conj:vanishing of cohomology GL_n integral} in the special case of Shimura varieties, and to describe the state of the art for results towards these conjectures. We especially focus on the references~\cite{CaraianiScholzeCompact},\cite{CaraianiScholzeNonCompact} and~\cite{koshikawa} which concern Shimura varieties of PEL type A and a completely split prime $p$. However, we end by mentioning work in progress of Hamann and Lee~\cite{hamann-lee} that goes significantly beyond this case. 

In the case of Shimura varieties, the analogue of Conjecture~\ref{conj:vanishing of cohomology GL_n integral} predicts that the non-Eisenstein part of the cohomology with $\F_{\ell}$-coefficients is concentrated in the middle degree $q_0 = d:= \dim_{\C}\mathrm{Sh}_K$. The initial progress on this conjecture in the Shimura variety setting had rather strong additional assumptions: for example one needed $\ell$ to be an unramified prime for the Shimura datum and $K_{\ell}$ to be hyperspecial, as in the work of Dimitrov~\cite{Dimitrov} and Lan--Suh~\cite{LanSuhCompact, LanSuhNonCompact}. The geometry of the Hodge--Tate period morphism, as described in \S~\ref{ss:geometry of pi_HT}, has been a game-changer in this area, allowing us to obtain results that are simultaneously more general and more precise. 

Assume that $(G,X)$ is a Shimura datum of PEL type A. Let $\ell$ be a prime number and $\m \subset \mathbb{T}^S$ be a maximal ideal in the support of $H^*_{(c)}(\mathrm{Sh}_K(\C), \F_{\ell})$. Scholze's argument with fake Hasse invariants explained in \S~\ref{ss:construction-torsion-Galois} also allows us to associate a global modulo $\ell$ continuous and semi-simple Galois representation $\bar{\rho}_{\m}$ to the system of Hecke eigenvalues determined by $\m$. The non-Eisenstein condition on $\m$ makes sense as a \emph{global} condition imposed on the Langlands parameter $\bar{\rho}_{\m}$. In what follows, we will choose an auxiliary prime $p\not = \ell$ that is unramified for the Shimura datum and impose \emph{local} conditions on the Langlands parameter at $p$. 

\begin{definition}\label{defn:generic} 
Let $k$ be a finite field of characteristic $\ell$. 
\begin{enumerate}
    \item Let $p\not = \ell$ be a prime, $L/\mathbb{Q}_p$ be a finite extension, and $\bar{\rho}: \Gamma_L \to \GL_n(k)$ be a continuous Galois representation. We say that $\bar{\rho}$ is  \emph{generic} if it is unramified and 
    the eigenvalues (with multiplicity) $\alpha_1, \dots , \alpha_n \in \overline{k}^\times$
    of $\bar{\rho}(\Frob_L)$ satisfy $\alpha_i / \alpha_j \not = |\cO_L/\m_L|$ for $i\not = j$. For a product $\bar{\rho}: \Gamma_L\to \prod_{i=1}^m\GL_{n_i}(k)$ of continuous Galois representations, we say that $\bar{\rho}$ is \emph{generic} if each factor in the product is generic.

    \item Let $F$ be a number field and $\bar{\rho}: \Gamma_F \to \GL_n(k)$
be a continuous representation. We say that a prime $p\not = \ell$ is 
\emph{decomposed generic} for $\bar{\rho}$ if $p$ splits completely in $F$
and, for every prime $\mathfrak{p}\mid p$ of $F$, $\bar{\rho}|_{\Gamma_{F_{\mathfrak{p}}}}$ is generic.  
We say that $\bar{\rho}$ is \emph{decomposed generic}
if there exists a prime $p\not = \ell$ which is decomposed generic for $\bar{\rho}$.
(If one such prime exists, then infinitely many do, by the Chebotarev density theorem.) These notions also extend to products $\bar{\rho}: \Gamma_F \to \prod_{i=1}^m \GL_{n_i}(k)$ as above.  

\end{enumerate}
\end{definition} 

\begin{remark} 
The condition for a local representation $\bar{\rho}: \Gamma_L \to \GL_n(k)$ 
to be generic implies that any lift to characteristic $0$ of $\bar{\rho}$ corresponds under the
local Langlands correspondence to an irreducible and generic principal series representation of $\GL_n(L)$. 
Such a representation can never arise from a non-split inner form of $\GL_n/L$ 
via the Jacquet--Langlands correspondence. For this reason, a generic $\bar{\rho}$
cannot be the modulo $\ell$ reduction of the 
$L$-parameter of a smooth representation of a non-split
inner form of $\GL_n/L$. 

The genericity condition can also be characterized in terms of the modulo $\ell$ semi-simple local Langlands correspondence of Vign\'eras~\cite{Vigneras}: the condition is that 
$\bar{\rho}$ corresponds to an irreducible (and thus generic) modulo $\ell$ principal series representation of $\GL_n(L)$. If this is the case, then the representation of $\GL_n(L)$ attached to $\bar{\rho}$ by the Emerton--Helm mod $\ell$ local Langlands in families~\cite{EmertonHelm} is the same irreducible and generic principal series representation. 
\end{remark}

\begin{remark}
    The condition for a global Galois representation $\bar{\rho}:\Gamma_F\to \GL_n(k)$ to be decomposed generic can be ensured in many cases when $\bar{\rho}$ has sufficiently large image, using the Chebotarev density theorem. For example, if $F$ is totally real, $n=2$, $\ell>2$ and $\bar{\rho}$ is totally odd and has non-solvable image, then it is decomposed generic, cf.~\cite[Lem.~7.1.8]{CaraianiTamiozzo}. A similar result holds if $F$ is CM and Galois over $\Q$, cf.~\cite[Lem.~2.3]{AllenNewton}. 
\end{remark}

As mentioned above, the idea is to impose some kind of genericity condition on a system of Hecke eigenvalues $\m$ in order to restrict the range of degrees in which $\m$ can contribute. The expectation is that, under such a condition, 
\begin{equation}\label{eq:cohomology inequalities}
H^i_c(\mathrm{Sh}_K(\C), \F_{\ell})_{\m} = 0\ \mathrm{if}\ i> d\ \mathrm{and}\ H^i(\mathrm{Sh}_K(\C), \F_{\ell})_{\m} = 0\ \mathrm{if}\ i < d. 
\end{equation}
If either $\m$ is non-Eisenstein or the Shimura varieties under consideration are compact, the above cohomology groups would be concentrated only in the middle degree. This would prove Conjecture~\ref{conj:vanishing of cohomology GL_n integral} in many cases for Shimura varieties of PEL type A. For this purpose, the decomposed generic condition works extremely well, as exhibited by the following result. 

\begin{theorem}[~\cite{CaraianiScholzeCompact,CaraianiScholzeNonCompact,koshikawa,MafaldaSantos}]\label{thm:decomposed generic}
Let $(G,X)$ be a Shimura datum of PEL type A. Let $\m\subset \mathbb{T}$ be a maximal ideal in the support of $H^*_{(c)}(\mathrm{Sh}_K(\C), \F_{\ell})$ with associated Galois representation $\bar{\rho}_{\m}$.  
Assume that $\bar{\rho}_{\m}$ is decomposed generic. Then the expectation~\eqref{eq:cohomology inequalities} holds true. 
\end{theorem}

\begin{remark}
    In the case when $(G,X)$ gives rise to a compact Shimura variety of Harris--Taylor type, Theorem~\ref{thm:decomposed generic} was 
first proved by Boyer~\cite{boyer}. His argument uses the canonical integral models of these Shimura
varieties. In this setting, Boyer also goes
\emph{beyond genericity}, in the following sense. Given the  eigenvalues (with multiplicity) 
$\alpha_1,\dots,\alpha_n$ of $\bar{\rho}_{\m}(\mathrm{Frob}_{\mathfrak{p}})$, with $\mathfrak{p}\mid p$ 
the relevant prime of $F$\footnote{In this special case, one does not have to 
impose the condition that $p$ splits completely in $F$, and it suffices
to have genericity at one prime $\mathfrak{p}\mid p$.}, one can define a ``defect'' that measures how far 
$\bar{\rho}_{\m}$ is from being generic at $\mathfrak{p}$. Concretely, set 
$\delta_{\p}(\m)$ to be equal to the length of the maximal chain of 
eigenvalues where the successive terms have ratio equal to $|\cO_{F_{\p}}/\m_{F_{\p}}|$. 
Boyer shows that the cohomology groups $H^i_{(c)}(\mathrm{Sh}_{K}(\C), \F_{\ell})_{\m}$ 
are non-zero at most in the range $[d-\delta_{\p}(\m), d+ \delta_{\p}(\m)]$. 
As noted by both Emerton and Koshikawa, such a result is 
consistent with Arthur's conjectures on the cohomology of Shimura varieties
with $\C$-coefficients and points towards a modulo $\ell$ analogue of these conjectures. 
\end{remark}

\begin{remark} Let $p\not = \ell$ be a prime that splits completely in $F$, that is unramified for the Shimura datum $(G,X)$, and such that $K_p\subset G(\Q_p)$ is the hyperspecial compact open subgroup. Then the genericity condition can also be formulated purely locally, using the spherical Hecke algebra at $p$. This is the approach taken in~\cite{koshikawa}. 
Indeed, the spherical Hecke algebra $\mathcal{H}_p$ for $G(\Q_p)$ is a product of spherical Hecke algebras for general linear groups. For each system of Hecke eigenvalues $\m_p\subset \mathcal{H}_p$, 
the associated local $L$-parameter is a product of unramified $L$-parameters valued in general linear groups. The notion of generic as in the first part of Definition~\ref{defn:generic} makes sense and Koshikawa establishes the inequalities in~\eqref{eq:cohomology inequalities} after localizing at a maximal ideal $\m_p$ with generic $L$-parameter. 
 
\end{remark}

For simplicity, let us first discuss the proof of Theorem~\ref{thm:decomposed generic} in the case of compact Shimura varieties. The two main approaches that work in this generality are due to~\cite{CaraianiScholzeCompact} and~\cite{koshikawa} and they are, in some sense, complementary. They both rely on the Hodge--Tate period morphism and on the intuition that the generic part of the cohomology of a PEL type A Shimura variety at a completely split prime should come only from the ordinary locus. The geometry of the Hodge--Tate period morphism is used to control the non-ordinary Newton strata, showing that the generic part of their cohomology vanishes. 

At this point, the approaches diverge. The best way to think about their difference is in terms of the infinite-level Mantovan product formula, which describes a given Newton stratum as a (twisted) product of an Igusa variety and a Rapoport--Zink space. The original approach due to~\cite{CaraianiScholzeCompact} controls the cohomology of Igusa varieties using global methods, particularly the computation of their cohomology via the trace formula. The later approach due to~\cite{koshikawa} controls the cohomology of Rapoport--Zink spaces via Fargues--Scholze~\cite{fargues-scholze} and its comparison with the classical construction of the local Langlands correspondence for $\GL_n$ due to Harris--Taylor~\cite{HarrisTaylor}. 

The original approach starts by analyzing the $\mathbb{T}^S$-equivariant diagram 

\begin{equation}\label{eq:hodge-tate diagram}
\xymatrix{\ & \cS_{K^p}\ar[dl]\ar[dr]^{\pi_{\mathrm{HT}}} & \ \\
\cS_{K^pK_p} & \ & \mathscr{F}\ell_G.}
\end{equation}

\noindent The standard comparison theorems between various cohomology 
theories allow us to identify $H^*(\mathrm{Sh}_K(\C), \F_{\ell})_{\m}$ with $H^*(\cS_{K}, \F_\ell)_{\m}$. 
The arrow on the left hand side of~\eqref{eq:hodge-tate diagram} 
is a $K_p$-torsor, so the Hochschild--Serre spectral sequence allows us to 
recover $H^*(\cS_{K}, \F_\ell)_{\m}$ from 
$H^*(\cS_{K^p}, \F_\ell)_{\m}$. The idea
is now to compute 
$H^*(\cS_{K^p}, \F_\ell)_{\m}$ in two stages: first understand 
the complex of sheaves $(R\pi_{\mathrm{HT}*}\F_{\ell})_{\m}$ on $\mathscr{F}\ell$, then compute the total cohomology using the Leray--Serre spectral sequence. 
In the decomposed generic case, the structure 
of $(R\pi_{\HT*}\F_{\ell})_{\m}$ is particularly simple.
\begin{enumerate}
\item Firstly, $(R\pi_{\HT*}\F_{\ell})_{\m}$
behaves like a perverse sheaf on $\mathscr{F}\ell_G$ (up to shift). This is because
$\pi_{\HT}$ is simultaneously \emph{affinoid}, as explained in Theorem~\ref{thm:perfectoid-Sh} and the subsequent discussion,  
and \emph{partially proper}, because the Shimura varieties were assumed to be compact. 
In particular, the restriction of 
$(R\pi_{\HT*}\F_{\ell})_{\m}$ to a highest-dimensional stratum in its support 
is concentrated in one degree. By the computation of the fibres of $\pi_{\HT}$, 
this implies that the cohomology groups $R\Gamma(\mathfrak{Ig}^b, \Z_{\ell})_{\m}$ 
are concentrated in one degree and torsion-free. 

\item Secondly, whenever the group
$G_b(\Q_p)$ acting on $\mathfrak{Ig}^b$ 
comes from a non-quasi-split inner form, 
the localisation $R\Gamma(\mathfrak{Ig}^b, \Q_{\ell})_{\m}$ vanishes. Recall from \S~\ref{ss:geometry of pi_HT} that we have an isomorphism $R\Gamma(\mathfrak{Ig}^b, \Q_{\ell})_{\m}\simeq R\Gamma(\mathrm{Ig}^b, \Q_{\ell})_{\m}$. The computation of $R\Gamma(\mathrm{Ig}^b, \Q_{\ell})_{\m}$,
at least at the level of the Grothendieck group, can be done 
using the trace formula method pioneered in~\cite[Ch.V]{HarrisTaylor} and \cite{ShinIgusa}, 
which is particularly well-suited to Shimura varieties of PEL type A (at least when working with a unitary group defined over a totally real field $F^+\not= \Q$).   

\item Finally, the condition that $p$ splits completely in $F$ guarantees that the only 
Newton stratum for which $G_b$ is quasi-split is the ordinary one. 

\end{enumerate}

\noindent Therefore, the hypotheses of Theorem~\ref{thm:decomposed generic} guarantee 
that $(R\pi_{\HT*}\F_{\ell})_{\m}$ is as simple as possible - it is 
supported in one degree on a zero-dimensional stratum! The degree in which it is concentrated is in fact equal to $d$. The Leray--Serre spectral sequence computing $H^*(\cS_{K^p}, \F_{\ell})_{\m}$ degenerates right away and the resulting cohomology is also concentrated in degree $d$. The same result can then also be deduced for $H^*(\cS_{K}, \F_{\ell})_{\m}$ by analyzing the Hochschild--Serre spectral sequence. 

\begin{remark}\label{rem:perversity} The intuition that $(R\pi_{\mathrm{HT}*}\F_{\ell})_{\m}$ behaves like a perverse sheaf (up to shift) can be made precise if we instead work with the morphism of small $v$-stacks 
\[
\pi_{\mathrm{HT}/K_p}: \cS_{K^pK_p} \to \mathscr{F}\ell_G/\underline{K_p}. 
\]
One can define an ad hoc perverse $t$-structure on $\mathscr{F}\ell_G/\underline{K}_p$ with respect to its Newton stratification and show that $(R\pi_{\mathrm{HT}/K_p*}\F_{\ell})_{\m}[d]$ is in the heart of this $t$-structure. The perverse $t$-structure on $\mathscr{F}\ell_G/\underline{K}_p$ is defined by glueing the shifted standard $t$-structures $(D^{\leq -\delta_b},D^{\geq -\delta_b})$ on each Newton stratum $\mathscr{F}\ell^b_G/\underline{K}_p$, with $\delta_b$ equal to the dimension of the stratum.  In fact, this perverse $t$-structure is pulled back from a perverse $t$-structure on $\mathrm{Bun}_G$ under the $\ell$-cohomologically smooth morphism $\mathscr{F}\ell_{G}/\underline{K_p}\to \mathrm{Bun}_G$ (see~\cite{hamann} for the definition of the $t$-structure on $\mathrm{Bun}_G$ and~\cite[Prop.~2.10]{Hansen} for $\ell$-cohomological smoothness). One then obtains semi-perversity of $(R\pi_{\mathrm{HT}/K_p*}\F_{\ell})_{\m}[d]$
in one direction from the affineness of Igusa varieties and Artin vanishing and in the other direction from Verdier duality, which uses properness of the Shimura varieties. We thank Matteo Tamiozzo for explaining the details of this argument to us.  
\end{remark}

We now explain the more recent approach due to~\cite{koshikawa}. In this case, one first proves a vanishing result for the generic part of the cohomology of local Shimura varieties. Let $L/\Q_p$ be a finite extension. Let $(G,b,\mu)$ be a local Shimura datum, with $G/L$ a group of the form $\prod_{i\in I}\GL_{n_i}$ for some finite index set $I$. We set $K_L:=\prod_{i\in I}\GL_{n_i}(\cO_L)\subset G(L)$ and we consider the associated Rapoport--Zink space $\cM_{K_L}$ at level $K_L$. We consider the spherical $\F_{\ell}$-Hecke algebra $\mathcal{H}_L$ acting on $H^*_c(\cM_{K_L}, \F_{\ell})$ with $\ell\not = p$. To a maximal ideal $\m_L\subset \mathcal{H}_L$ we have an associated unramified Galois representation $\bar{\rho}_{\m_L}\to \prod_{i\in I}\GL_{n_i}(\overline{\F}_{\ell})$. We also let $G_b$ denote the twisted Levi subgroup of $G$ determined by the element $b$.  

\begin{theorem}[\cite{koshikawa}]\label{thm:local generic} Let $\m_L\subset \mathcal{H}_L$ be a maximal ideal with associated Galois representation $\bar{\rho}_{\m_L}$. Assume that $\bar{\rho}_{\m_L}$ is generic  as in Definition~\ref{defn:generic} and that $G_b$ is not quasi-split. Then 
$H^*_c(\cM_{K_L}, \F_{\ell})_{\m_L} = 0$ 
\end{theorem}

The proof of Theorem~\ref{thm:local generic} relies on two steps. Firstly, Koshikawa proves that, if $G_b$ is not quasi-split and if $\pi$ is an irreducible smooth $\overline{\F}_{\ell}$-representation of $G_b(L)$, then the semi-simple $L$-parameter $\varphi_{\pi}$ associated to it by Fargues--Scholze cannot be generic. 
This uses the compatibility between the $L$-parameters constructed by Fargues--Scholze and the ones constructed via classical Jacquet--Langlands and local Langlands; this compatibility holds up to semi-simplification in the case of $\GL_n$ by~\cite{HKW}. Secondly, Koshikawa uses the action of the spectral Bernstein center to show that, for any irreducible subquotient $\pi$ of $H^*_c(\cM_{K_L}, \overline{\F}_{\ell})_{\m_L}$, the Fargues--Scholze parameter satisfies $\varphi_{\pi} = \bar{\rho}_{\m_L}$. 

Let us return to the global perspective of a Shimura variety of PEL type A and assume that $p$ is an unramified prime for the Shimura datum that splits completely in the underlying CM field $F$. Then the only Newton stratum for which the group $G_b$ is quasi-split is the ordinary one. 
Assume that $\m\subset \mathbb{T}$ is a maximal ideal in the support of $H^*(\mathrm{Sh}_K(\C), \F_{\ell})$ such that $\bar{\rho}_{\m}$ is decomposed generic. By the cohomological version of the Mantovan product formula, cf.~\cite[Thm.~7.1]{koshikawa}, and by Theorem~\ref{thm:local generic}, the only contribution to $H^*(\cS_{K^p}, \F_{\ell})_{\m}$ comes from the ordinary locus. At this point, Theorem~\ref{thm:decomposed generic} follows from the perversity of $(R\pi_{\HT*}\F_{\ell})_{\m}[d]$ (in the sense used originally in~\cite{CaraianiScholzeCompact}, or in the sense of Remark~\ref{rem:perversity}). 

For the non-compact case of Theorem~\ref{thm:decomposed generic}, the key conceptual difference is that one must weaken perversity to semi-perversity. For example, if one considers the restriction of the Hodge--Tate period morphism to the good reduction locus
\[
\pi^\circ_{\HT}: \cS^{\circ}_{K^p}\to \mathscr{F}\ell_G,
\] 
it is shown in~\cite[\S4.6]{CaraianiScholzeNonCompact} that $(R\pi^\circ_{\HT*}\F_{\ell})_{\m}[d]$ is, in some sense, semi-perverse. This uses the computation of the fibers of the Hodge--Tate period morphism for a toroidal compactification of the Shimura variety. Such a result can again be made precise as in Remark~\ref{rem:perversity} by showing that $(R\pi^\circ_{\HT/K_p*}\F_{\ell})_{\m}[d]$ is in $^pD^{\geq 0}$ for the perverse $t$-structure on $\mathscr{F}\ell_G/\underline{K}_p$. (There is also a dual result, which applies to $(R\pi^\circ_{\HT/K_p!}\F_{\ell})_{\m}[d]$ and shows that this is in $^pD^{\leq 0}$. This dual result can be obtained from the computation of the fibers of the Hodge--Tate period morphism for the minimal compactification of the Shimura variety. These fibers are related to partial minimal compactifications of Igusa varieties, which are affine and therefore satisfy Artin vanishing.) 

In the original strategy of~\cite{CaraianiScholzeNonCompact}, once one has semi-perversity, one still needs to control the boundary of Shimura variety. This is done via an intricate induction argument; along the way, one constructs Galois representations attached to torsion classes occurring in the cohomology of Igusa varieties. In the more recent strategy of~\cite{koshikawa}, semi-perversity together with Theorem~\ref{thm:local generic} are essentially all that is needed. The proof of Theorem~\ref{thm:decomposed generic} in the non-compact case is completed in~\cite{MafaldaSantos} by extending semi-perversity to more general PEL cases. 

Several different generalizations of Theorem~\ref{thm:decomposed generic} are now emerging. Let us mention the upcoming work of Hamann--Lee~\cite{hamann-lee}, who work with a prime that does not necessarily split completely in $F$
and who adapt their method to more general reductive groups under an axiomatic set-up. The key assumption in their axiomatic set-up is the compatibility between the work of Fargues--Scholze and more classical approaches to the local Langlands correspondence. This is known, for example, for $\mathrm{GSp}_4$ and for odd unitary groups. The key new ingredient in their work is the theory of geometric Eisenstein series over the Fargues--Fontaine curve developed in~\cite{hamann} in the case of the Borel subgroup. In the cases considered by Hamann--Lee, only the so-called unramified elements in $B(G_{\Q_p})_{\mathrm{un}}\subset B(G_{\Q_p})$ contribute to the cohomology of the Shimura variety and, moreover, the contribution of each $b\in B(G_{\Q_p})_{\mathrm{un}}$ can be controlled with the theory of geometric Eisenstein series. 

\begin{remark} Recall the Igusa stack $\overline{\pi}_{\HT}: \mathrm{Igs}^\circ_{K^p}\to \mathrm{Bun}_G$ constructed in~\cite{mingjia} and discussed in Remark~\ref{rem:Mingjia's PhD thesis}. 
One advantage of the original approach to Theorem~\ref{thm:decomposed generic} via Igusa varieties is that it describes the sheaves $(R\bar{\pi}_{\HT!}\F_{\ell})_{\m}$ on $\mathrm{Bun}_G$, at least under the assumption that $p$ splits completely in $F$. These sheaves are expected to satisfy a form of local-global compatibility for the conjectural equivalence of categories of Fargues--Scholze. In fact, with the methods introduced in~\cite{hamann}, one can investigate this question in a slightly more general setting, cf.~\cite[\S 8.3]{mingjia}. The cohomology of Igusa varieties is admissible, therefore the sheaves $R\bar{\pi}_{\HT*}\F_{\ell}$ are ULA objects in $D(\mathrm{Bun}_G, \F_{\ell})$, cf.~\cite[Thm.~1.5.1]{fargues-scholze}. Given a semi-simple $L$-parameter $\phi$ which satisfies an appropriate genericity condition, one can use the spectral action of Fargues--Scholze to define the ``localisation'' $(R\bar{\pi}_{\mathrm{HT}*}\overline{\F}_{\ell})_{\phi}$.
Assume the compatibility between Fargues--Scholze and more classical approaches to local Langlands. The techniques of Hamann and Hamann--Lee (specifically~\cite[Prop.~4.15]{hamann-lee}) show that 
\begin{equation}\label{eq:unramified elements}
(R\bar{\pi}_{\mathrm{HT}*}\overline{\F}_{\ell})_{\phi}
\simeq \bigoplus_{b\in B(G_{\Q_p})_{\mathrm{un}}} R\Gamma(\mathrm{Ig}^b, \overline{\F}_{\ell})_{\phi},
\end{equation}
where the terms on the RHS of the decomposition are supported on $[*/\widetilde{G}_b]$ for $b\in B(G_{\Q_p})_{\mathrm{un}}$. 
Furthermore, in the case when the Igusa stack agrees with its minimal compactification, these terms are all supported in one degree by the perversity result discussed in Remark~\ref{rem:perversity}. 

\end{remark}
 
\subsection{Applications to locally symmetric spaces}\label{ss:boundary}

In this subsection, we return to the case when $F$ is an imaginary CM field and $G=\GL_n/F$. For $K\subset G(\A_F^\infty)$ a sufficiently small compact open subgroup, recall that $Y_K$ is the associated locally symmetric space. In this subsection, we explain how Theorem~\ref{thm:decomposed generic} can be used to approach Conjecture~\ref{conj:integral LGC} for the cohomology groups of $Y_K$. We then explain how to implement the Calegari--Geraghty method unconditionally over CM fields.  

As we have seen in \S~\ref{sss:Borel-Serre}, $G$ can be realized as a Levi subgroup of a maximal parabolic subgroup of a quasi-split unitary group $\widetilde{G}=U_{2n}/F^+$ and $\mathrm{Res}_{F^+/\Q}\widetilde{G}$ admits a Shimura datum. For a neat compact open subgroup $\widetilde{K}\subset \widetilde{G}(\A_{F^+}^\infty)$, we let $\Sh_{\widetilde{K}}$ denote the corresponding Shimura variety for (the restriction of scalars of) $\widetilde{G}$, considered in this subsection as a complex manifold. We let $\Sh_{\widetilde{K}}^\mathrm{BS}$ denote its Borel--Serre compactification, a real manifold with corners, and  $\partial \Sh_{\widetilde{K}}^\mathrm{BS}$ denote the boundary of the Borel--Serre compactification. 

We let $S$ be a finite set of finite places of $F$ that is stable under complex conjugation. This determines a finite set of finite places $\bar{S}$ of $F^+$. 
We have the usual abstract Hecke algebra $\mathbb{T}^S$ for $G$. We consider an analogous Hecke algebra $\widetilde{\mathbb{T}}^{\bar{S}}$ for $\widetilde{G}$ (although we only need to take the product of local spherical Hecke algebras at places in $F^+$ that split in $F$). There is an unnormalized Satake transform map $\cS: \widetilde{\mathbb{T}}^{\bar{S}} \to \mathbb{T}^S$; see~\cite[\S2]{newton-thorne} for its definition.

Recall from \S~\ref{ss:construction-torsion-Galois} that the excision sequence for the boundary of the Borel--Serre compactification induces the $\widetilde{\mathbb{T}}^{\bar{S}}$-equivariant long exact sequence

\begin{multline}\label{eq:excision les}
  \cdots \to H^i_c(\Sh_{\widetilde{K}}, \Z/\ell^m\Z)\to H^i(\Sh_{\widetilde{K}}, \Z/\ell^m\Z)\to  \\ 
  H^i(\partial\Sh_{\widetilde{K}}^{\mathrm{BS}},\Z/\ell^m\Z)\to 
  H^{i+1}_c(\Sh_{\widetilde{K}}, \Z/\ell^m\Z)\to \cdots. 
\end{multline}
 Assume that a maximal ideal $\widetilde{\m}\subset \widetilde{\mathbb{T}}^{\bar{S}}$ satisfies an appropriate genericity condition at an auxiliary prime $p\not = \ell$. Theorem~\ref{thm:decomposed generic} applies to the Shimura varieties $\Sh_{\widetilde{K}}$\footnote{More precisely, one considers a variant of $\widetilde{G}$ that is a unitary similitude group, giving rise to Shimura varieties of PEL type A. The difference is only on the level of geometric connected components.} and can be used to simplify the long exact sequence~\eqref{eq:excision les} after localisation at $\widetilde{\m}$. 
 In degree $d = \dim_{\C}\Sh_{\widetilde{K}}$, we obtain the following diagram: 
\begin{equation}\label{eq:middle degree diagram}
H^d(\mathrm{Sh}_{\widetilde{K}}, \Z_{\ell})_{\widetilde{\m}} \otimes_{\Z_{\ell}}\Q_{\ell}\hookleftarrow  H^d(\mathrm{Sh}_{\widetilde{K}}, \Z_{\ell})_{\widetilde{\m}} \twoheadrightarrow
H^d(\partial\Sh_{\widetilde{K}}^{\mathrm{BS}}, \Z_{\ell})_{\widetilde{\m}},
\end{equation} 
which is used to approach Conjecture~\ref{conj:integral LGC}. 
Very roughly, the cohomology groups on the LHS can be described in terms of (essentially) conjugate self-dual, regular $C$-algebraic cuspidal automorphic representations of $\GL_{r}/F$ with various $r\le 2n$\footnote{In order to make this work in reality, one needs to introduce a non-trivial local system $\cV_{\tilde{\lambda}}$ on $\Sh_{\widetilde{K}}$, which can eliminate contributions from non-conjugate self-dual automorphic forms by imposing a central character constraint.}. 
Their Galois representations are known to satisfy local-global compatibility up to Frobenius semi-simplification at all places, including at the places $v\mid \ell$ of $F$. This ultimately gives control over the cohomology group on the RHS of~\eqref{eq:middle degree diagram}. 

On the other hand, if $\widetilde{\m} = \cS^*(\m)$, the cohomology groups $H^i(Y_K, \Z/\ell^m\Z)_{\m}$ for $i=0,\dots d-1$ can be shown to contribute to the cohomology group on the RHS of~\eqref{eq:middle degree diagram}. 
This is the basic strategy employed in~\cite{10authors} and later refined via $P$-ordinary parts in~\cite{CaraianiNewton}. The strategy requires there to be enough primes of $F^+$ dividing $\ell$: we fix one such prime $\bar{v}\mid \ell$ where we want to prove local-global compatibility, and we use the primes $\bar{v}'\mid \ell$ with $\bar{v}'\not = \bar{v}$ to create congruences and shift to cohomological degree $d$ on the boundary of the Borel--Serre compactification. (This is the reason why Conjecture~\ref{conj:integral LGC} seems to be most subtle when $F$ is an imaginary quadratic field -- there are not enough auxiliary primes!) 

\begin{remark} There can be genuine torsion classes in the cohomology group $H^d(\partial\Sh_{\widetilde{K}}^{\mathrm{BS}}, \Z_{\ell})_{\widetilde{\m}}$. The diagram~\eqref{eq:middle degree diagram} implies that these lift to characteristic $0$ classes in $H^d(\mathrm{Sh}_{\widetilde{K}}, \Z_{\ell})_{\widetilde{\m}} \otimes_{\Z_{\ell}}\Q_{\ell}$. By design, the proofs of Conjecture~\ref{conj:integral LGC} using the method outlined above go through these torsion classes.  
\end{remark}

At this point, we have access to Conjectures~\ref{conj:existence of Galois representations 2} and~\ref{conj:integral LGC}, at least under some technical assumptions. In order to implement the Calegari--Geraghty method unconditionally as sketched in \S~\ref{ss:CG sketch}, we also need to access Conjecture~\ref{conj:vanishing of cohomology GL_n integral}. This seems extremely hard in this generality. However, a modification of the Calegari--Geraghty method due to Khare--Thorne~\cite{khare-thorne} works in certain settings by assuming only Conjecture~\ref{conj:vanishing of cohomology GL_n rational}. This predicts the vanishing of cohomology with $\Q_{\ell}$-coefficients, a result that can be proved when $F$ is a CM field. One of the main challenges in~\cite{10authors} was to make this insight of Khare--Thorne compatible with other techniques in automorphy lifting, which rely on reduction modulo $\ell$. This challenge was resolved  by considering reduction modulo $\ell$ from a derived perspective, leading to an unconditional implementation of the Calegari--Geraghty method in arbitrary dimensions.  

We end by briefly mentioning some applications of this circle of ideas: a proof of the Sato--Tate conjecture for elliptic curves over CM fields, cf.~\cite{10authors}, a proof that all elliptic curves over many imaginary fields are modular, cf.~\cite{AKT, CaraianiNewton}, and proofs of many new cases of the Ramanujan conjecture for $\mathrm{GL}_2$ over CM fields, cf.~\cite{10authors, 5authors}. 

\medskip


\bibliographystyle{amsplain}


\bibliography{bib}

\end{document}